\documentclass[mnsc,nonblindrev]{informs3-vvm}

\OneAndAHalfSpacedXI %

\usepackage{subcaption}

\usepackage{amsmath}
\usepackage{amssymb}
\usepackage{natbib}
 \bibpunct[, ]{(}{)}{,}{a}{}{,}%
 \def\bibfont{\small}%

\usepackage{color}
 \usepackage[table]{xcolor} 
 \usepackage{ctable}
 \usepackage{multirow}
\usepackage{xspace}

\usepackage{algorithm}
\usepackage{algorithmic}

\def\ab{\mathbf{a}}
\def\bb{\mathbf{b}}
\def\cb{\mathbf{c}}
\def\db{\mathbf{d}}
\def\hb{\mathbf{h}}
\def\qb{\mathbf{q}}
\def\ub{\mathbf{u}}

\def\vb{\mathbf{v}}
\def\wb{\mathbf{w}}
\def\xb{\mathbf{x}}
\def\yb{\mathbf{y}}
\def\zb{\mathbf{z}}

\def\zerob{\mathbf{0}}
\def\oneb{\mathbf{1}}

\def\betab{\boldsymbol \beta}

\def\lambdab{\boldsymbol \lambda}
\def\mub{\boldsymbol \mu}

\def\pib{\boldsymbol \pi}
\def\taub{\boldsymbol \tau}
\def\xib{\boldsymbol \xi}
\def\varphib{\boldsymbol \varphi}

\def\Sigmab{\boldsymbol \Sigma}
\def\Xib{\boldsymbol \Xi}

\def\Ab{\mathbf{A}}
\def\Cb{\mathbf{C}}

\def\zerob{\mathbf{0}}

\def\Acal{\mathcal{A}}
\def\Kcal{\mathcal{K}}
\def\Lcal{\mathcal{L}}
\def\Ncal{\mathcal{N}}
\def\Pcal{\mathcal{P}}
\def\Scal{\mathcal{S}}
\def\Tcal{\mathcal{T}}
\def\Ucal{\mathcal{U}}
\def\Vcal{\mathcal{V}}
\def\Xcal{\mathcal{X}}

\def\Ibb{\mathbb{I}}
\def\Rbb{\mathbb{R}}

\def\leftnode{\mathbf{left}}
\def\rightnode{\mathbf{right}}

\def\modelP{\textsf{P}\xspace}
\def\modelPRobustS{\textsf{P-Robust}\xspace}
\def\modelPRobustSBudget{\textsf{P-Robust-Budget}\xspace}
\def\modelRA{\textsf{RA}\xspace}
\def\modelPRPT{\textsf{P-RPT}\xspace}
\def\modelSOC{\textsf{SOC}\xspace}

\def\KKDP{\textsf{KKDP}\xspace}
\def\Greedy{\textsf{Greedy}\xspace}
\def\LS{\textsf{LS}\xspace}
\def\GM{\textsf{GM}\xspace}

\def\WCL{\mathrm{WCL}}
\def\RI{\mathrm{RI}}

\def\Halmos{$\square$}

\newenvironment{proofvvm}{\paragraph{Proof:}}{\\[1em]}

\def\Halmos{$\square$}

\def\toubia{{\tt timbuk2}\xspace}

\def\bank{{\tt bank}\xspace}
\def\candidate{{\tt candidate}\xspace}
\def\immigrant{{\tt immigrant}\xspace}

\TheoremsNumberedThrough     %
\ECRepeatTheorems

\EquationsNumberedThrough    %
\begin{document}
\RUNTITLE{Exact Logit-Based Product Design}

\TITLE{Exact Logit-Based Product Design}

\ARTICLEAUTHORS{%
	\AUTHOR{\.{I}rem Akchen}
	\AFF{UCLA Anderson School of Management, University of California, Los Angeles, California 90095, United States, \EMAIL{\tt irem.akchen@gmail.com}}
	\AUTHOR{Velibor V. Mi\v{s}i\'{c}}
	\AFF{UCLA Anderson School of Management, University of California, Los Angeles, California 90095, United States, \EMAIL{\tt velibor.misic@anderson.ucla.edu}} %
}
\ABSTRACT{

The share-of-choice product design (SOCPD) problem is to find the product, as defined by its attributes, that maximizes market share arising from a collection of customer types or segments. When customers follow a logit model of choice, the market share is given by a weighted sum of logistic probabilities, leading to the logit-based share-of-choice product design problem. In this paper, we develop a methodology for solving this problem to provable optimality. We first analyze the complexity of this problem, and show that this problem is theoretically intractable: it is NP-Hard to solve exactly, even when there are only two customer types, and it is furthermore NP-Hard to approximate to within a non-trivial factor. Motivated by the difficulty of this problem, we propose three different mixed-integer exponential cone programs of increasing strength for solving the problem exactly, which allow us to leverage modern integer conic program solvers such as Mosek. Using both synthetic problem instances and instances derived from real conjoint data sets, we show that our methodology can solve large instances to provable optimality or near optimality in operationally feasible time frames and yields solutions that generally achieve higher market share than previously proposed heuristics.

}%


\maketitle

\allowdisplaybreaks

\section{Introduction}
\label{sec:introduction}

Consider the following canonical marketing problem. A firm has to design a product, which has a collection of attributes, and each attribute can be set to one of a finite set of levels. The product will be offered to a collection of customers, which differ in their preferences and specifically, in the utility that they obtain from different levels of different attributes. What product should the firm offer -- that is, to what level should each attribute be set -- so as to maximize the share of customers who choose to purchase the product? This problem is referred to as the share-of-choice product design (SOCPD) problem, and has received a significant amount of attention in the marketing science research literature. 

The problem of designing a single product is of significant practical interest; we name a couple of notable historical examples. The paper of \cite{wind1989courtyard} details an academic consulting engagement initiated by Marriott to design a new hotel chain, Courtyard by Marriott, using conjoint analysis. The design of this new hotel chain involved specifying 50 attributes with a total of 167 levels. As another example, the paper of \cite{vavra1999evaluating} details another academic consulting engagement that resulted in the design of the EZPass toll system by the state of New Jersey and New York. The problem of designing a single product is particularly relevant in settings where launching the product is sufficiently capital intensive (e.g., a new hotel chain) that it is not possible to launch multiple versions of the same product (i.e., a product line). 

The SOCPD problem is a challenging problem for several reasons. First, since a product corresponds to a combination of attribute levels, the number of candidate products scales exponentially with the number of attributes, and can be enormous for even a modest number of attributes. This, in turn, renders solution approaches based on brute force enumeration computationally cumbersome. Second, it is common to represent customers using discrete choice models that are built on the multinomial logit model. Under this assumption of customer behavior, the problem becomes more complex, because the purchase probability under a logit choice model is a nonlinear function of the product design's utility that is neither convex nor concave. Finally, product design problems in real life settings may also often involve constraints, arising from engineering or other considerations, which can further constrain the set of candidate products.

In this paper, we consider the logit-based SOCPD problem. In this problem, the firm must design a product that maximizes the expected number of customers who choose to purchase a product, where customers are assumed to follow logit models of choice, and the probability of a customer purchasing a product is given by a logistic response function (i.e., the function $f(u) = e^u / (1 + e^u)$). We propose an exact solution methodology for this problem that is based on modern integer, convex and conic optimization. To the best of our knowledge, this is the first exact solution methodology for the logit-based SOCPD problem. 

We make the following specific contributions:

\begin{enumerate}
\item We formally define the logit-based SOCPD problem. We show that this problem is NP-Hard, even when there are two customer types. We further show that the problem is NP-Hard to approximate the problem to within a factor $O(1 / n^{1-\epsilon})$, where $n$ is the number of product attributes and $\epsilon > 0$. Lastly, we also establish that when each customer type has up to three non-zero partworths, the problem is APX-Hard, and when each customer type has exactly three non-zero partworths, is NP-Hard to approximate to a factor better than $7/8 + \epsilon$ for any $\epsilon > 0$.
\item We develop three different mixed-integer exponential cone program formulations of the problem of increasing tightness. The first formulation relies on a characterization of logit probabilities as being the optimal solutions to a representative agent problem, in which an agent chooses the probability of selecting two alternatives so as to maximize a regularized expected utility. The second formulation is derived by applying a perspective function-based convexification of the logit probability expression. The third formulation relies on applying the reformulation-perspectification technique (RPT) of \cite{zhen2021extension} to the second formulation. We additionally show how our second formulation can be extended to account for uncertainty in partworths. Lastly, we also show that the problem can be formulated as a mixed-integer second order cone program.
\item We demonstrate the practical tractability of our approach using synthetic problem instances, as well as a set of problem instances derived from real conjoint data sets. Using synthetic problem instances with up to $n = 70$ attributes and up to $K = 30$ customer types, we show that our approach can solve the logit-based SOCPD problem to within an optimality gap of 10\% within two hours, and solutions obtained by our approach outperform heuristic solutions, in some cases by as much as 30\%. On our real problem instances, we are able to solve the logit-based SOCPD problem to provable optimality in all cases within ten minutes, and we again find that our solutions outperform those obtained by heuristics.
\end{enumerate}

The rest of this paper is organized as follows. Section~\ref{sec:literature_review} provides a review of the related literature. Section~\ref{sec:model} provides a formal definition of the logit-based SOCPD problem, and all of our complexity results. Section~\ref{sec:micp} presents our three mixed-integer exponential cone formulations of the logit-based SOCPD problem. Section~\ref{sec:numerical_experiments} presents the results of our numerical experiments. Lastly, in Section~\ref{sec:conclusions}, we conclude. All proofs are relegated to the electronic companion. %

\section{Literature Review}
\label{sec:literature_review}

We divide our literature review according to four subsets: the single product design literature; the product line design literature; the representative agent model literature; and lastly, the broader optimization literature.

\paragraph{Single product design:} Product design has received significant attention in the marketing science community; we refer readers to \cite{schmalensee1988perceptual} and \cite{green2004buyer} for overviews of this topic. The majority of papers on the SOCPD problem assume that customers follow a deterministic, first-choice model, i.e., they purchase the product if the utility exceeds a ``hurdle'' utility, and do not purchase it otherwise. Many papers have proposed heuristic approaches for this problem; examples include \cite{kohli1987heuristic,kohli1989optimal} and \cite{balakrishnan1996genetic}. Other papers have also considered exact approaches based on branch-and-bound \citep{camm2006conjoint} and nested partitions \citep{shi2001optimization}. %

The main difference between our work and the majority of the prior work on the product design problem is the use of a logit-based share-of-choice objective function. As stated earlier, when the share-of-choice is defined as the sum of logit probabilities, the SOCPD problem becomes a discrete nonlinear optimization problem, and becomes significantly more difficult than the SOCPD problem when customers follow first-choice/max-utility models. To the best of our knowledge, our approach is the first approach for obtaining provably optimal solutions to the SOCPD problem when customers follow a logit model.

\paragraph{Product line design/assortment optimization:} Besides the product design problem, a more general problem is the product line design (PLD) problem, where one must select several products so as to either maximize the share-of-choice, the expected profit or some other criterion. A number of papers have considered the PLD problem under a first-choice model of customer behavior, where customers deterministically select the product with the highest utility; examples of such papers include \cite{mcbride1988integer}, \cite{kohli1990heuristics}, \cite{wang2009branch}, \cite{belloni2008optimizing} and \cite{bertsimas2019exact}. Besides the first-choice model, several papers have also considered the PLD problem under the (single-class) multinomial logit model \citep[see, e.g.,][]{chen2000mathematical,schon2010optimal}. Other work has also considered objective functions corresponding to a worst-case expectation over a family of choice models \citep{bertsimas2017robust}.

Outside of the marketing literature, the PLD problem is closely related to the problem of assortment optimization which appears in the operations management literature. In this problem, one must select a set of products from a larger universe of products so as to maximize expected revenue. The difference between PLD and assortment optimization arises from where the choice model comes: in PLD, the choice model usually comes from conjoint survey data and the task is to select a set of new products, whereas in assortment optimization, typically the set of candidate products consists of products that have been offered in the past, and the choice model is estimated from historical transactions. There is an extensive literature on solving this problem under a variety of choice models, such as the single class MNL model \citep{talluri2004revenue}, the nested logit model \citep{davis2014assortment} and the Markov chain choice model \citep{feldman2017revenue}; we refer readers to \cite{gallego2019assortment} for a recent overview of the literature.  %

Our paper differs from the product line and assortment optimization literatures in that we focus on the selection of a single product, and the decision variables of our optimization problem are the attributes of the product. In contrast, virtually all mathematical programming-based approaches to PLD/assortment optimization require one to input a set of candidate products, and the main decision variable is a set of products from the overall set of candidate products. The attributes of the products are only relevant in specifying problem data (e.g., in an MNL assortment problem, one would determine the utilities of the candidate products from their attributes), but do not directly appear as decision variables. %

Lastly, we comment on a couple of papers in the assortment optimization literature that use conic optimization, and which are related to ours. The first is the working paper of \cite{shao2021tractable}, which considers a product line design and pricing problem where the goal is to select continuous attributes for a collection of products under the single-class multinomial logit model, the Markov chain choice model and the nested logit model, so as to maximize profit. For the MNL model in particular, \cite{shao2021tractable} develop a continuous nonlinear program for the profit maximization problem which can also be cast as an exponential cone program. The main technique of formulating the MNL probability involves effectively converting a non-convex inequality constraint into a convex one by multiplying both sides by a decision variable corresponding to a choice probability, which is similar to the perspectification technique used to arrive at one of our formulations (model~\modelP). Despite this similarity, there are a number of important differences between that paper and our paper. First, \cite{shao2021tractable} consider a single-class MNL model where attributes are continuous, whereas in our paper, we consider a mixture MNL model where attributes are discrete. Second, our paper proposes three other formulations for the logit-based SOCPD problem, one based on the representative agent model, one based on the very recently proposed reformulation-perspectification technique \citep{zhen2021extension}, and one based on a rotated second order cone representation. Besides these differences in the problem setting and formulations, our paper also differs in several other critical aspects: we develop a number of important intractability results for the logit-based SOCPD problem (Section~\ref{subsec:model_complexity}); we extend our approach to address uncertainty, using robust optimization (Section~\ref{sec:robust} of the ecompanion); and we perform extensive computational experiments using both synthetic and real problem instances. Our numerical experiments are also closely connected to the marketing research literature, in that (1) we compare our methods against commonly used heuristics from the marketing literature, and (2) for our real data instances we estimate real mixture MNL models from real conjoint datasets, using both the EM algorithm and the MCMC method applied to a hierarchical Bayesian formulation of the choice model. 

The second is the paper of \cite{sen2018conic}, which develops a mixed-integer second order cone program for solving the mixture of MNL assortment optimization problem. Our three main formulations in Section~\ref{sec:micp} (\modelRA, \modelP and \modelPRPT) differ from this formulation in that they are based on the exponential cone, whereas the formulation in \cite{sen2018conic} is based on rotated second-order cone constraints. A fourth formulation that we present in the ecompanion, dubbed \modelSOC (see Section~\ref{sec:SOC}), shares a similarity to the formulation of \cite{sen2018conic}, in that a constraint of the form $x(1+e^u) \geq 1$ is modeled using a rotated second order cone constraint, which allows the formulation to lower bound the no-purchase probability. However, the main idea of our formulation \modelSOC, which is based on representing a term of the form $e^{ \beta_{k,1} a_1 + \dots + \beta_{k,n} a_n}$ using a collection of rotated second order cone constraints, is novel and distinct from \cite{sen2018conic}.

\paragraph{Representative agent model:} One of our formulations (formulation~\modelRA) is based on a characterization of logit probabilities as solutions of a concave maximization problem where the decision variables correspond to the choice probabilities and the objective function is the entropy-regularized expected utility. This concave maximization problem is an example of a representative agent model, and has been studied in a number of papers in the economics and operations management literatures \citep{anderson1988representative,hofbauer2002global,feng2017relation}. The goal of our paper is not to contribute directly to this literature, but rather to leverage one such result so as to obtain an exact and computationally tractable reformulation of the logit-based SOCPD problem. To the best of our knowledge, the representative agent-based characterization of logit probabilities has not been previously used in optimization models arising in marketing or operations; we believe that this characterization could potentially be useful in other contexts outside of product design.

\paragraph{Optimization literature:} Lastly, we comment on the relation of our paper to the general optimization literature. Our paper contributes to the growing literature on mixed-integer convex and mixed-integer conic programming. In the optimization community, there has been an increasing interest in developing general solution methods for this class of problems \citep[see, for example,][]{lubin2018polyhedral,coey2020outer} as well as understanding what types of optimization problems can and cannot be modeled as mixed-integer convex programs \citep{lubin2017mixed}. At the same time, mixed-integer convex and mixed-integer conic programming have been used in a variety of applications, such as power flow optimization \citep{lubin2019chance}, robotics \citep{liu2020new}, portfolio optimization \citep{benson2013mixed}, joint inventory-location problems \citep{atamturk2012conic} and designing battery swap networks for electric vehicles \citep{mak2013infrastructure}. All of our mathematical programming formulations rely on the exponential cone, and thus our paper contributes to a growing set of applications of exponential cone programming, which include scheduling charging of electric vehicles \citep{chen2021exponential}, robust optimization with uncertainty sets motivated by estimation objectives \citep{zhu2021joint} and manpower planning \citep{jaillet2018strategic}.

Outside of this literature, we note that a couple of prior papers have considered the problem of designing a product to maximize the share-of-choice under a mixture of logit models. The first is the paper of \cite{udell2013maximizing} that considers the sum of sigmoids optimization problem, which is an optimization problem where the objective function is a sum of sigmoid (S-shaped) functions; the logistic response function $f(u) = e^u / (1+e^u)$ is a specific type of sigmoid function. The paper of \cite{udell2013maximizing} develops a general purpose branch-and-bound algorithm for solving this problem when the decision variables are continuous. Our paper differs in that we focus specifically on an objective that corresponds to a sum of logistic response functions, and the main decision variables in our formulation are binary variables, indicating the presence or absence of certain attributes. In addition, our formulation is an exact reformulation of the problem into a mixed-integer convex problem, which can then be solved directly using a commercial mixed-integer conic solver (such as Mosek). In contrast, the approach of \cite{udell2013maximizing} requires one to solve the problem using a custom branch-and-bound algorithm. 

The second is the paper of \cite{huchette2017nonconvex}, which develops a mixed-integer linear programming formulation for general nonconvex piecewise linear functions. As an example of the application of the framework, the paper applies the framework to the problem of deciding on continuous product attributes to maximize a logit-based share-of-choice objective, which involves approximating the logistic response function $f(u) = e^u / (1+e^u)$ using a piecewise linear function. As with our discussion of \cite{udell2013maximizing}, our formulation differs in that it is exact, and that the attributes are discrete rather than continuous.

Lastly, we note that our paper also contributes to a growing literature on optimization models where the objective function to be optimized is obtained from a predictive model or a machine learning model. Some examples include work on optimizing objective functions obtained from tree ensemble models (such as random forests; see \citealt{ferreira2016analytics}, \citealt{mivsic2020optimization}) and neural networks \citep{anderson2020strong}.

\section{Model}
\label{sec:model}

We begin by formally defining our model in Section~\ref{subsec:model_problem_definition}. We then discuss the computational complexity of this model in Section~\ref{subsec:model_complexity}. %

\subsection{Problem definition}
\label{subsec:model_problem_definition}

We assume that there are $n$ binary attributes. We assume the product design is described by a binary vector $\ab = (a_1, \dots, a_n)$, where $a_i$ denotes the presence of attribute $i$. We let $\Acal \subseteq \{0,1\}^n$ denote the set of feasible attribute vectors. While we formulate the problem in terms of binary attributes, we note that this is without loss of generality, as one can represent an attribute with $M$ levels using $M-1$ binary attributes, and one can specify $\Acal$ to include a constraint that requires at most one of the new $M-1$ binary attributes to be selected.

We assume that there are $K$ different segments or \emph{customer types}. Each customer type is associated with a nonnegative weight $\lambda_k$, which is the fraction of customers who belong to that type/segment, or alternatively the probability that a customer belongs to that type/segment; note that we always have that $\sum_{k=1}^K \lambda_k = 1$. Each customer type is also associated with a partworth vector $\betab_k = (\beta_{k,1}, \dots, \beta_{k,n}) \in \mathbb{R}^n$, where $\beta_{k,i}$ is the partworth of attribute $i$. In addition, we let $\beta_{k,0} \in \mathbb{R}$ denote the constant part of the customer's utility. Given a candidate design $\ab \in \Acal$, the customer's utility for the product is given by 
\begin{equation*}
u_k(\ab) = \beta_{k,0} + \sum_{i=1}^n \beta_{k,i} a_i.
\end{equation*}
We assume that each customer type is choosing between our product design corresponding to the vector $\ab$, and an outside/no-purchase option. Without loss of generality, we fix the utility of the outside option to zero. This assumption is not restrictive, as choice probabilities under the logit model are unaffected when all of the utilities are adjusted by a constant. In particular, an equivalent representation (one that would lead to the same choice probabilities) is to specify the utility of the product as $\sum_{i=1}^n \beta_{k,i} a_i$ and the utility of the no-purchase option as $- \beta_{k,0}$. As a result, the constant term $\beta_{k,0}$ effectively captures the utility of the no-purchase option.

We assume that each customer type chooses to buy or not buy the product according to a multinomial logit model. Thus, given $\ab \in \Acal$, the customer chooses to purchase the product with probability $ \exp( u_k(\ab) ) / (1 + \exp( u_k(\ab)))$ and chooses the outside option with probability $1 / (1 + \exp(u_k(\ab)))$.

With these definitions, the logit-based share-of-choice product design problem can then be defined as 
\begin{equation}
\underset{ \ab \in \Acal}{\text{maximize}} \sum_{k=1}^K \lambda_k \cdot \frac{ \exp( u_k(\ab)) }{1 + \exp(u_k(\ab))}. \label{prob:SOCPD_abstract}
\end{equation}
The objective function of this problem can be thought of as the share or fraction of all customers who choose to purchase the product, or the (unconditional) probability that a random customer chooses to purchase the product.

\subsection{Complexity}
\label{subsec:model_complexity}

In this section, we characterize the complexity of the logit-based SOCPD problem. 
Our first major theoretical result is that problem~\eqref{prob:SOCPD_abstract} is intractable even in the case that there are two customer types.
\begin{theorem}
The logit-based SOCPD problem~\eqref{prob:SOCPD_abstract} with $\Acal = \{0,1\}^n$ and with $K = 2$ customer types is NP-Hard. 
\label{theorem:NPHard_Keq2}
\end{theorem}
We establish this result by considering the decision form of problem~\eqref{prob:SOCPD_abstract}, which asks whether there exists a product attribute vector $\ab$ that achieves a share-of-choice of at least $\theta$, where $\theta$ is a user specified parameter. We show that the decision form of problem~\eqref{prob:SOCPD_abstract} is equivalent to the partition problem, which is another well-known NP-Hard problem \citep{garey2002computers}. We note that our result was inspired by and shares some similarity with another result from the paper of \cite{rusmevichientong2014assortment}. In particular, Theorem~3.2 of the paper of \cite{rusmevichientong2014assortment} shows that the problem of assortment optimization under the mixture of multinomial logits (MMNL) model is NP-Hard when the number of customer classes/types is equal to 2. This result is similar to ours in that our underlying choice model is also a mixture of multinomial logits/latent-class multinomial logit model, and that paper also establishes this result using the partition problem. Although there is a similarity in the choice model used and the use of the partition problem, the problem of assortment optimization under the MMNL model is quite different from the share-of-choice product design problem, because the objective function of the MMNL assortment optimization problem is a sum of weighted linear fractional functions of a binary vector, i.e., it is a problem of the form $\max_{\xb \in \{0,1\}^n} \sum_{k=1}^K \lambda_k \frac{ \sum_{i=1}^n r_i w_{k,i} x_i}{ 1 + \sum_{i=1}^n w_{k,i} x_i}$ (where $n$ is the number of candidate products, $\xb$ is a binary vector that encodes for each product $i \in \{1,\dots,n\}$ whether it is offered or not, $r_i$ is the marginal profit/revenue of product $i$, and $w_i$ is the preference weight of product $i$ for customer type $k$). The nonlinearity in this problem arises from the ratios of linear functions. Note that the exponential function does not appear, because it is ``baked into'' the $w_{k,i}$'s: each $w_{k,i}$ can be thought of as $w_{k,i} = e^{u_{k,i}}$, where $u_{k,i}$ is the mean utility of product $i$ for customer type $k$. In contrast, in our logit-based SOCPD problem the decision variable is also a binary vector $\ab \in \{0,1\}^n$, but the objective function has a more complicated dependence on this vector $\ab$ through the exponential function. As a result the proof of Theorem~\ref{theorem:NPHard_Keq2} is quite different from that in \cite{rusmevichientong2014assortment}, and is not a straightforward adaptation of the proof of Theorem~3.2 from \cite{rusmevichientong2014assortment}.

Our second result in this section concerns the ability to approximate problem~\eqref{prob:SOCPD_abstract}. Letting $F^*$ denote the optimal objective value of problem~\eqref{prob:SOCPD_abstract}, we say that an algorithm achieves an approximation factor of $C \in [0,1]$ if it is guaranteed to produce a solution $\ab$ whose objective value is at least $C \cdot F^*$. We then have the following result.
\begin{theorem}
The logit-based SOCPD problem~\eqref{prob:SOCPD_abstract} with $\Acal = \{0,1\}^n$ is NP-Hard to approximate to within a factor of $O(1 / n^{1-\epsilon})$ for any $\epsilon > 0$. \label{theorem:APXHard}
\end{theorem}

We establish this result by designing an approximation-preserving reduction between the logit-based SOCPD problem and the maximum independent set problem, for which the same inapproximability result holds \citep{hastad1996clique}. We note that this result also shares some similarities with known results in the assortment optimization literature. In particular, the excellent papers of \cite{aouad2018approximability} and \cite{desir2022capacitated} respectively showed that the assortment optimization problem under ranking preferences and the MMNL assortment optimization problem are NP-Hard to approximate to within a factor better than $O( n^{1-\epsilon})$ for $\epsilon > 0$, where $n$ is the number of candidate products, also using the maximum independent set problem (see Theorem~1 of \citealt{aouad2018approximability} and Theorem~2 of \citealt{desir2022capacitated}). As with our discussion of our Theorem~\ref{theorem:NPHard_Keq2}, the proof of our Theorem~\ref{theorem:APXHard} differs from these existing results because the dependence of the objective function of problem~\eqref{prob:SOCPD_abstract} on the binary vector $\ab$ of product attributes is completely different from how the expected revenue under the ranking-based model or under the MMNL model depend on the binary vector $\xb$ that encodes which products are included/excluded from the assortment. As a result, the proof of Theorem~\ref{theorem:APXHard} is not a direct adaptation of the proofs of these prior results. 

Lastly, we establish that problem~\eqref{prob:SOCPD_abstract} is APX-Hard, even in the case that each customer type has up to three non-zero partworths, i.e., $|\{i \in [n] \mid \beta_{k,i} \neq 0\}| \leq 3$ for all $k \in \{1,\dots,K\}$. 

\begin{theorem}
The logit-based SOCPD problem~\eqref{prob:SOCPD_abstract} with $\Acal = \{0,1\}^n$ and with each customer type having up to three non-zero partworths is APX-Hard. In addition, when $\Acal = \{0,1\}^n$ and with exactly three non-zero partworths for each customer type, problem~\eqref{prob:SOCPD_abstract} is NP-Hard to approximate to within factor $7/8 + \epsilon$ for any $\epsilon > 0$. \label{theorem:APXHard_MAX3SAT}
\end{theorem}

The proof of Theorem~\ref{theorem:APXHard_MAX3SAT} follows by establishing an approximation preserving reduction between the logit-based SOCPD problem and the MAX-3SAT problem, which is known to be APX-Complete \citep{papadimitriou1991optimization}. Additionally, using the known result that MAX-E3SAT (the special case of MAX-3SAT where each clause has exactly three literals) is NP-Hard to approximate to within factor $7/8 + \epsilon$ for any positive $\epsilon$ \citep{haastad2001some}, we can establish a similar inapproximability result for problem~\eqref{prob:SOCPD_abstract}. Note that in our two previous results, Theorems~\ref{theorem:NPHard_Keq2} and \ref{theorem:APXHard}, the reductions rely on setting up the instance of problem~\eqref{prob:SOCPD_abstract} in such a way that the number of non-zero partworths for each customer type can be as large as $n$. Theorem~\ref{theorem:APXHard_MAX3SAT} establishes that even if each customer type's purchase probability is affected by even as few as three attributes, the problem remains difficult to even approximate. 

In terms of approximation algorithms, Theorem~\ref{theorem:APXHard} implies that an upper bound on the best possible approximation factor is $O(1 / n)$ in the most general case (i.e., with no restrictions on the problem parameters), while Theorem~\ref{theorem:APXHard_MAX3SAT} implies that, when each customer type has up to three non-zero partworths, the best possible approximation algorithm one can hope for is a constant factor approximation algorithm. As a complement to these results, we also develop an approximation algorithm (see Section~\ref{sec:approximation_algorithm} in the ecompanion) that runs in polynomial time when one assumes that the number of customer types $K$ and the input size are both treated as constant. The development of more powerful approximation algorithms (e.g., a fully polynomial time approximation scheme for the case where $K = O(1)$) is an interesting direction for future research.

The main takeaway from these results is that problem~\eqref{prob:SOCPD_abstract} is fundamentally a very difficult problem to solve. In the next section, we develop exact mathematical programming formulations for this problem.

\section{Mixed-Integer Convex Programming Formulations}
\label{sec:micp}

Motivated by the hardness results that we established for problem~\eqref{prob:SOCPD_abstract} in Section~\ref{subsec:model_complexity}, in this section, we develop exact mathematical programming formulations of problem~\eqref{prob:SOCPD_abstract}. 

The formulations that we will develop belong to the general class of mixed-integer convex programming (MICONVP) problems, and specifically, the class of mixed-integer conic programming (MICP) problems; we provide a brief overview of MICP here. An MICP problem has the following general form:
\begin{subequations}
\begin{alignat}{2}
& \underset{\xb \in \mathbb{R}^n}{\text{minimize}} & & \cb^T \xb \\ 
& \text{subject to} & \quad & \Ab \xb - \bb \in \Kcal, \\
& & & x_i \in \mathbb{Z}, \quad \forall \ i \in I, \label{prob:MICP_canonical_integrality}
\end{alignat}%
\label{prob:MICP_canonical}%
\end{subequations}
where $\xb$ is an $n$-dimensional vector of decision variables, $\cb$ is an $n$-dimensional vector, $\bb$ is an $m$-dimensional vector, $\Ab$ is an $m$-by-$n$ matrix, $I \subseteq [n]$ is the set of integer variables in the problem and finally, $\Kcal$ is a closed convex cone. A closed convex cone $\Kcal$ is a closed subset of $\Rbb^n$ that contains all nonnegative combinations of its elements, i.e., a set $\Kcal$ satisfying the following property:
\begin{equation}
\yb_1, \yb_2 \in \Kcal \ \Rightarrow \ \alpha_1 \yb_1 + \alpha_2 \yb_2 \in \Kcal\ \text{for any}\ \alpha_1, \alpha_2 \geq 0.  \label{eq:cone_definition}
\end{equation}
While the cone $\Kcal$ can in theory be specified as any set that satisfies \eqref{eq:cone_definition}, in practice, it is common to model $\Kcal$ as a Cartesian product of a collection of cones drawn from the set of standard cones. An example of a standard cone is the the cone $\Kcal_{\geq} = \{ \yb \in \Rbb^m \mid \yb \geq \zerob \}$, where $\zerob$ is an $m$-dimensional vector of zeros. This cone is known as the nonnegative cone, as it corresponds to the nonnegative orthant of $\Rbb^n$. The constraint $\Ab \xb - \bb \in \Kcal_{\geq}$ is equivalent to the constraint $\Ab \xb \geq \bb$, which is just a linear constraint. Other standard cones include the zero cone, the second order cone, the exponential cone and the positive semidefinite cone. The three formulations presented here will involve functions that are representable using the exponential cone, which is defined as 
\begin{equation}
\Kcal_{\exp} = \{ (r, 0, t) \in \mathbb{R}^3 \mid r \geq 0, t \leq 0\} \cup \{ (r,s,t) \in \mathbb{R}^3 \mid s > 0, r \geq s \exp( t / s) \}.
\end{equation}

We refer readers to \cite{aps2021mosekcookbook} for an overview of other standard cones. 

Having described mixed-integer conic programming in generality, we now elaborate on why this representation is valuable. Many mixed-integer convex programs involve complicated nonlinear functions. Until recently, the method of choice for tackling such problems has been to use mixed-integer nonlinear programming solvers, which treat these nonlinear functions in a ``black-box'' fashion and rely on evaluating these functions and their derivatives to solve the problem. Often, it turns out that constraints involving these nonlinear functions can be re-written through additional variables and additional conic constraints involving standard cones.\footnote{A notable recent example of this is the paper of \cite{lubin2018polyhedral}, which found that all 194 mixed-integer convex programming problems in the MINLPLIB2 (\url{http://www.gamsworld.org/minlp/minlplib2/html/}) benchmark library could be represented as mixed-integer conic programs using standard cones.} In so doing, one obtains a mixed-integer conic program, which is then amenable to solution methods for such problems. This is important because conic programs -- problems of the same form as \eqref{prob:MICP_canonical}, without the integrality constraint~\eqref{prob:MICP_canonical_integrality} -- are considered to be among the easiest continuous nonlinear programs to solve: the theory of numerical algorithms for solving these problems is quite developed, there are numerous software packages for solving these problems at practical scale, and there continues to be active development both in the theory and in software implementations. Solution algorithms for mixed-integer conic programs are built on top of algorithms for (continuous) conic programs and can exploit the conic structure. By formulating the problem as a mixed-integer conic program, one is able to use state-of-the-art commercial solvers such as Mosek \citep{aps2021moseksoftware}, as well as new open-source solvers such as Pajarito \citep{coey2020outer} to solve the problem. (In our numerical experiments in Section~\ref{sec:numerical_experiments} we will indeed solve all of our formulations using Mosek.)

The challenge in developing an MICP formulation of problem~\eqref{prob:SOCPD_abstract} is how to model probabilities of the form $e^{u}/(1+e^{u})$; as we have already discussed, such probabilities have a non-convex dependence on $u$. The first formulation we will present in Section~\ref{subsec:micp_RA}, formulation \modelRA, is based on leveraging the fact that logit probabilities arise as optimal solutions of a concave maximization problem called the representative agent model. The second formulation we will present in Section~\ref{subsec:micp_P}, formulation \modelP, is based on a technique known as perspectification. The last formulation we will present in Section~\ref{subsec:micp_PRPT}, formulation \modelPRPT, is based on applying the reformulation-perspectification technique (RPT; \citealt{zhen2021extension}) to formulation \modelP. In all three of our formulations, the representation of the logit probabilities will critically depend on certain convex functions that are representable using the exponential cone, and our formulations will thus be mixed-integer exponential cone programs, which are supported by Mosek as of 2019.  Finally, in Section~\ref{subsec:micp_extensions}, we briefly describe a couple of extensions of our models here (on handling a profit objective and handling uncertainty in problem data), as well as an alternate approach to exactly modeling logit probabilities using second order cones; we discuss these extensions in more detail in the ecompanion.

Before we continue, we make the following assumption on the structure of the set $\Acal$.
\begin{assumption}
The set $\Acal$ can be written as $\Acal = \{ \ab \in \{0,1\}^n \mid \Cb \ab \leq \db \}$ for some choice of $\Cb \in \mathbb{R}^{m \times n}$ and $\db \in \mathbb{R}^m$, where $m \in \mathbb{Z}_+$. \label{assumption:Acal_polyhedron}
\end{assumption}
Assumption~\ref{assumption:Acal_polyhedron} just requires that the set of candidate products $\Acal$ can be represented as the set of binary vectors satisfying a finite collection of linear inequality constraints. 
This assumption is necessary in order to ensure that our problem can be formulated as a mixed-integer convex program of finite size. 
We note that this assumption is not too restrictive, as many natural constraints can be expressed in this way; we discuss some examples in Section~\ref{subsec:extra_Acal} of the ecompanion.

\subsection{Formulation \modelRA}
\label{subsec:micp_RA}
 
Our first formulation relies on the representative agent model, which we now briefly review. In the representative agent model, an agent is faced with $M$ alternatives. Each alternative $m \in [M]$ is associated with a utility $\pi_m$, where we use the notation $[N] = \{1,\dots,N\}$ for any positive integer $N$. The agent must choose the probability $x_m$ with which each alternative will be selected; we let $\xb = (x_1,\dots, x_M)$ be the probability distribution over the $M$ alternatives. The agent seeks to maximize his adjusted expected utility, where the adjustment is achieved through a convex regularization function $V(\xb)$. The representative agent model is then the following optimization problem:
\begin{subequations}
\begin{alignat}{2}
& \underset{\xb}{\text{maximize}} & \quad  & \sum_{i=1}^M \pi_m x_m - V(\xb) \\
& \text{subject to} & & \sum_{m=1}^M x_m = 1, \\
& & & x_m \geq 0, \quad \forall \ m \in [M]. 
\end{alignat}
\end{subequations}
Since the function $V(\cdot)$ is a convex function, the above problem is a concave maximization problem. By carefully choosing the function $V$, the optimal solution of this problem -- the probability distribution $\xb$ -- can be made to coincide with choice probabilities under different choice models. In particular, it is known that the function $V(\xb) = \sum_{i=1}^M x_i \log (x_i)$ gives choice probabilities corresponding to a multinomial logit model \citep{anderson1988representative}. We refer the reader to the excellent paper of \cite{feng2017relation} for a complete characterization of which discrete choice models can be captured by the representative agent model. 

For our problem, the specific instantiation of the representative agent model that we are interested in is one corresponding to the choice of the $k$th customer type between our product and the no-purchase option. We let $x_{k,1}$ denote the probability of choosing our product with attribute vector $\ab$, and $x_{k,0}$ denote the probability of choosing the no-purchase option. Recall that the utility of our product is $u_k(\ab)$, and the utility of the no-purchase option is 0. The representative agent model for this customer type can thus be formulated as
\begin{subequations}
\begin{alignat}{2}
& \underset{x_{k,0}, x_{k,1} } { \text{maximize} } & \quad & u_k(\ab) \cdot x_{k,1} + 0 \cdot x_{k,0} - x_{k,1} \log(x_{k,1}) - x_{k,0} \log(x_{k,0}) \\
& \text{subject to} & & x_{k,1} + x_{k,0} = 1, \\
& & & x_{k,1}, x_{k,0} \geq 0. 
\end{alignat}
\label{prob:representative_agent_Meq2}%
\end{subequations}
The following theoretical result establishes  two key properties of this problem. First, the unique optimal solution $(x^*_{k,1}, x^*_{k,0})$ is exactly the pair of logit choice probabilities for the two alternatives. Second, the optimal objective value can be found in closed form. The proof of this proposition follows straightforwardly by analyzing the Lagrangean of problem~\eqref{prob:representative_agent_Meq2}, and is thus omitted.
\begin{proposition}
The unique optimal solution $(x^*_{k,1}, x^*_{k,0})$ of problem~\eqref{prob:representative_agent_Meq2} is given by 
\begin{align*}
x^*_{k,1} & =   \frac{e^{u_k(\ab)}}{1 + e^{u_k(\ab)}}, \\
x^*_{k,0} & = \frac{1}{ 1+ e^{u_k(\ab)} }.
\end{align*} 
In addition, the optimal objective value is $\log(1 + e^{ u_k(\ab)})$. 
\label{proposition:representative_agent_Meq2_optimality}
\end{proposition}

Using this result, we can now proceed with the formulation of our SOCPD problem. With a slight abuse of notation, let $u_k$ be a decision variable that denotes the utility of the candidate product $\ab$ for customer type $k$. As before, let $x_{k,1}$ and $x_{k,0}$ denote the probability of customer type $k$ purchasing and not purchasing the product, respectively. Then, the optimization problem can be formulated as
\begin{subequations}
\begin{alignat}{2}
& \underset{\ab, \ub, \xb}{\text{maximize}} & \quad  & \sum_{k=1}^K \lambda_k \cdot x_{k,1} \\
& \text{subject to} & & x_{k,1} + x_{k,0} = 1, \quad \forall k \in \{1,\dots, K\}, \label{prob:bilinear_x_unitsum} \\ 
& & & u_k x_{k,1}  - x_{k,1} \log(x_{k,1}) - x_{k,0} \log(x_{k,0}) \geq \log(1 + \exp(u_k)), \quad \forall k \in \{1,\dots, K\}, \label{prob:bilinear_mainconstraint} \\ 
& & & u_k = \beta_{k,0} + \sum_{i=1}^n \beta_{k,i} a_i, \quad \forall k \in \{1,\dots, K\}, \label{prob:bilinear_utility} \\ 
& & & \Cb \ab \leq \db, \\
& & & a_i \in \{0,1\}, \quad \forall i \in \{1,\dots,n\}, \\
& & & x_{k,1}, x_{k,0} \geq 0, \quad \forall k \in \{1,\dots, K\}. \label{prob:bilinear_x_nonnegative}
\end{alignat}
\label{prob:bilinear}
\end{subequations}
Observe that this formulation is a valid formulation of problem~\eqref{prob:SOCPD_abstract}. To see this, observe that any solution $(\ab,\ub, \xb)$ that is feasible has the property that for every $k \in [K]$, $(x_{k,0}, x_{k,1})$ is an optimal solution to problem~\eqref{prob:representative_agent_Meq2}. This is because constraint~\eqref{prob:bilinear_mainconstraint} enforces that the objective function of \eqref{prob:representative_agent_Meq2} is at least as good as $\log(1+e^{u_k})$, which is equal to $\log(1 +e^{u_k(\ab)})$ (because constraint~\eqref{prob:bilinear_utility} ensures that the decision variable $u_k$ is exactly equal to $u_k(\ab)$). Since we know that $\log(1+e^{u_k(\ab)})$ is the optimal objective value of problem~\eqref{prob:representative_agent_Meq2}, it follows that $(x_{k,0}, x_{k,1})$ is an optimal solution to problem~\eqref{prob:representative_agent_Meq2}. Additionally, since the solution of the representative agent problem~\eqref{prob:representative_agent_Meq2} is unique, it must be that $x_{k,0} = 1/(1+e^{u_k(\ab)})$, $x_{k,1} = e^{u_k(\ab)}/(1 + e^{u_k(\ab)})$. 
The key feature of this formulation is that it no longer explicitly involves the logit choice probabilities, which are a nonconvex function of $u_k$. We note that this problem is \emph{almost} a mixed-integer convex program. In the main nonlinear constraint~\eqref{prob:bilinear_mainconstraint}, the functions $- x_{k,1} \log(x_{k,1})$ and $-x_{k,0} \log(x_{k,0})$ appearing on the left hand side are instances of the \emph{entropy function} $- x \log(x)$ \citep{boyd2004convex} and are concave in $x_{k,1}$ and $x_{k,0}$. Similarly, the function $\log(1 + \exp(u_k))$ appearing on the right hand side, which is known as the \emph{softplus function} \citep{aps2021mosekcookbook}, is a convex function of $u_k$. Thus, this constraint can almost be written in the form $f(u_k, \xb_k) \leq 0$, where $f$ is a convex function. The main obstacle that prevents us from doing this is the bilinear term $u_k x_{k,1}$, which is not jointly concave in $u_k$ and $x_{k,1}$. 

Fortunately, we can use the fact that $u_k = \beta_{k,0} + \sum_{i=1}^n \beta_{k,i} a_i$ to re-write this bilinear term as 
\begin{align*}
u_k x_{k,1} & = (\beta_{k,0} + \sum_{i=1}^n \beta_{k,i} a_i) x_{k,1} \\
& = \beta_{k,0} x_{k,1} + \sum_{i=1}^n \beta_{k,i} \cdot a_i x_{k,1}.
\end{align*}
Using the fact that each $a_i \in \{0,1\}$ and that $0 \leq x_{k,1} \leq 1$, we can now linearize the bilinear terms $a_i x_{k,1}$ using a standard modeling technique from integer programming. In particular, we introduce a continuous decision variable $y_{k,i}$ for each $k$ and $i$ which corresponds to the product $a_i x_{k,1}$, and a continuous decision variable $w_k$ which corresponds to the product $u_k x_{k,1}$. This leads to the following equivalent formulation, which we denote by \modelRA: 
{\allowdisplaybreaks
\begin{subequations}
\begin{alignat}{2}
\modelRA: \quad & \underset{\ab, \ub, \wb, \xb, \yb}{\text{maximize}} & \quad & \sum_{k=1}^K \lambda_k \cdot x_{k,1} \\
& \text{subject to} & & x_{k,1} + x_{k,0} = 1, \quad \forall k \in [K], \label{prob:MICONVP_x_unitsum} \\ 
& & & w_k  - x_{k,1} \log(x_{k,1}) - x_{k,0} \log(x_{k,0}) \geq \log(1 + \exp(u_k)), \quad \forall k \in [K],  \label{prob:MICONVP_main} \\ 
& & & u_k = \beta_{k,0} + \sum_{i=1}^n \beta_{k,i} a_i, \quad \forall k \in [K], \\
& & & w_k = \beta_{k,0} x_{k,1} + \sum_{i=1}^n \beta_{k,i} y_{k,i}, \quad \forall k \in [K], \\ 
& & & y_{k,i} \leq x_{k,1}, \quad \forall k \in [K], \ i \in [n], \\ 
& & & y_{k,i} \leq a_i, \quad \forall k \in [K], \ i \in [n], \\ 
& & & y_{k,i} \geq a_i - 1 + x_{k,1}, \quad \forall k \in [K], \ i \in [n], \\ 
& & & y_{k,i} \geq 0, \quad \forall k \in [K], \ i \in [n], \\ 
& & & \Cb \ab \leq \db, \\
& & & a_i \in \{0,1\}, \quad \forall i \in [n], \\
& & & x_{k,1}, x_{k,0} \geq 0, \quad \forall k \in [K]. \label{prob:MICONVP_x_nonnegative}
\end{alignat}
\label{prob:MICONVP}
\end{subequations}
}
Note that in the formulation above, when $\ab \in \Acal$, $y_{k,i}$ is forced to take the value of $a_i \cdot x_{k,1}$, and $w_k$ takes the value of $u_k \cdot x_{k,1}$, which ensures the correctness of the formulation. We note that the nonlinear functions that appear in formulation~\modelRA, which are the entropy functions $x_{k,1} \log(x_{k,1})$ and $x_{k,0} \log(x_{k,0})$ and the softplus function, are representable using the exponential cone; we refer readers to \cite{aps2021mosekcookbook} for more details. As a result, formulation~\modelRA is a mixed-integer exponential cone program that can be solved using Mosek. %

\subsection{Formulation \modelP}
\label{subsec:micp_P}

In this section, we develop an alternate formulation of the logit-based SOCPD problem using a trick that is colloquially known in the nonlinear programming/convex optimization literature as \emph{perspectification} \citep{zhen2021extension}. We briefly review the idea of perspectification, and then demonstrate how it can be applied in our setting. 

In nonconvex optimization problems, one sometimes encounters a constraint of the form 
\begin{equation*}
g(x, \yb) + x f( \yb) \leq 0
\end{equation*}
where $x \geq 0$ is a scalar decision variable, $\yb$ is a vector of decision variables, $f$ is convex in $\yb$ and $g$ is jointly convex in $x$ and $\yb$. This constraint is not a convex constraint because although the term $x f(\yb)$ is marginally convex in each of $x$ and $\yb$ (specifically, it is linear in $x$ for any fixed $\yb$ and convex in $\yb$ for a fixed $x \geq 0$), it is not jointly convex in $(x,\yb)$. However, one way in which one can turn this constraint into a convex one is as follows. First, we multiply and divide $\yb$ by $x$ inside the function $f$:
\begin{equation*}
g(x,\yb) + x f\left( \frac{ x \yb}{x} \right)  \leq 0,
\end{equation*}
where we assume that $x f( \ub / x) = 0$ when $x = 0$. We now replace $x \yb$ with a new decision variable vector $\ub$, which serves as the linearization of the vector of bilinear terms $x \yb$. This leads to the constraint
\begin{equation*}
g(x,\yb) + x f \left( \frac{ \ub }{x} \right) \leq 0,
\end{equation*}
This last constraint now is in fact a convex constraint, because the new term $x f( \ub / x)$, is exactly the \emph{perspective function} of $f$. Recall that the perspective function of $f$ is the function $\tilde{f}(\yb, t) = t \cdot f( \yb / t)$, defined for $t \geq 0$. The perspective function is significant because whenever $f$ is a convex function of $\yb$, then the perspective function $\tilde{f}$ is a convex function of $\yb$ and $t$ \citep{boyd2004convex}. By including additional constraints that appropriately constrain the $\ub$ variable vector, we can potentially obtain a good relaxation of the original problem. This idea has been identified and exploited in a number of recent papers that consider challenging optimization problems \citep[for example][]{zhen2022disjoint,gorissen2022hidden}.

Let us now see how this idea applies in the context of our problem. Recall that in our problem, ideally we wish to enforce the constraint 
\begin{equation}
x_{k,1} = \frac{1}{1+ e^{-u_k}}
\end{equation}
where $x_{k,1}$ is the decision variable that represents the purchase probability of customer type $k$, and $u_k$ is the utility of the product for customer type $k$. Note that since the objective function is $\sum_{k=1}^K \lambda_k x_{k,1}$, which is a nonnegative weighted combination of the $x_{k,1}$ variables, we can safely relax the equality to an inequality, to obtain the constraint
\begin{equation}
x_{k,1} \leq \frac{1}{1+ e^{-u_k}}.
\end{equation}
If we now move the denominator of the right-hand side to the left, we get
\begin{equation}
x_{k,1} + x_{k,1} e^{-u_k} \leq 1.
\end{equation}
Here, recall that $f(y) = e^y$ is a convex function of $y$. We can thus apply the perspectification idea by multiplying and dividing the argument of $e^{\cdot}$ by $x_{k,1}$, which yields
\begin{equation}
x_{k,1} + x_{k,1} e^{ \frac{- x_{k,1} u_k}{x_{k,1}}} \leq 1.
\end{equation}
Finally, as in formulation~\modelRA, let us use $w_k$ to denote the linearization of $u_k \cdot x_{k,1}$. We now replace $x_{k,1} u_k$ with $w_k$, to arrive at the following convex constraint:
\begin{equation}
x_{k,1} + x_{k,1} e^{-w_k / x_{k,1} } \leq 1.
\end{equation}
This suggests the following formulation of the logit-based SOCPD problem, which we refer to as formulation~\modelP. Note that the definitions of the decision variables are the same as in formulation~\modelRA, and that constraints~\eqref{prob:P_w_definition} - \eqref{prob:P_x_definition} are the same as in formulation~\modelRA. 
\begin{subequations}
\begin{alignat}{2}
\modelP: \quad & \underset{\ab, \wb, \xb, \yb}{\text{maximize}} & \quad & \sum_{k=1}^K \lambda_k x_{k,1} \\
& \text{subject to} & & x_{k,1} + x_{k,1} e^{-w_k / x_{k,1} } \leq 1, \quad \forall \ k \in [K],  \label{prob:P_perspective}\\
& & & w_k = \beta_{k,0} x_{k,1} + \sum_{i=1}^n \beta_{k,i} y_{k,i}, \quad \forall \ k \in [K], \label{prob:P_w_definition} \\
& & & y_{k,i} \leq a_i, \quad \forall \ k \in [K], \ i \in [n], \\
& & & y_{k,i} \leq x_{k,1}, \quad \forall \ k \in [K], \ i \in [n], \\
& & & y_{k,i} \geq x_{k,1} + a_i - 1, \quad \forall \ k \in [K], \ i \in [n], \\
& & & y_{k,i} \geq 0, \quad \forall \ k \in [K], \ i \in [n], \\
& & & x_{k,0}, x_{k,1} \geq 0, \quad \forall \ k \in [K], \label{prob:P_x_definition} \\
& & & \Cb \ab \leq \db, \\
& & & \ab \in \{0,1\}^n. \label{prob:P_a_definition}
\end{alignat}
\label{prob:P}%
\end{subequations}
Note that this formulation is a valid formulation of problem~\eqref{prob:SOCPD_abstract}. To see this, observe that when $\ab \in \Acal$, $y_{k,i}$ will take the value $a_i \cdot x_{k,1}$ and $w_k$ will take the value $u_k \cdot x_{k,1}$, where $u_k = \beta_{k,0} + \sum_{i=1}^n \beta_{k,i} a_i$ is the utility of the product for customer type $k$. Thus, constraint~\eqref{prob:P_perspective} essentially enforces that $x_{k,1} + x_{k,1} e^{ - u_k \cdot x_{k,1} / x_{k,1}} \leq 1$, or equivalently, that $x_{k,1} \leq 1/ (1 + e^{-u_k})$. Since the objective function is a nonnegative weighted combination of the $x_{k,1}$ variables, at optimality we will have that $x_{k,1}$ will be equal to this upper bound, i.e., $x_{k,1} = 1/(1+e^{-u_k})$. Note that unlike formulation~\modelRA, it is no longer necessary to include a decision variable $u_k$ to represent the utility of the product for customer $k$. 

We make a couple of important remarks about this formulation. First, we note that formulation~\modelP, like formulation~\modelRA, can also be represented as a mixed-integer exponential cone program and solved directly using Mosek. However, one advantage that formulation~\modelP has over \modelRA is that \modelP requires only a single exponential cone constraint per customer type to represent the perspective function $x_{k,1} e^{ -w_k / x_{k,1}}$, whereas \modelRA requires four (one for each of the entropy functions, $x_{k,1} \log x_{k,1}$ and $x_{k,0} \log x_{k,0}$, and two for the softplus function $\log(1 + e^{u_k})$). Thus, formulation~\modelP should be easier to solve in general. On this point, another advantage that formulation~\modelP has over \modelRA is in regard to numerical stability. In formulation~\modelRA, constraint~\eqref{prob:MICONVP_main} by design must hold at equality at integer solutions, which can lead to ill-posedness issues \citep{aps2021mosekcookbook}. On the other hand, constraint~\eqref{prob:P_perspective} does not have to hold at equality for integer solutions, and so formulation~\modelP does not have this same issue of ill-posedness.

Second, a natural question is how formulation~\modelP and \modelRA compare. Here, the continuous relaxations of formulations \modelP and \modelRA are the convex optimization problems that one obtains when the constraint $\ab \in \{0,1\}^n$ is replaced with the constraint
\begin{equation}
0 \leq a_i \leq 1, \quad \forall \ i \in [n].
\end{equation}
The relaxation bound is the objective value that is obtained when one solves the continuous relaxation. Let $Z^*_{\modelP}$ and $Z^*_{\modelRA}$ denote the optimal values of the continuous relaxations of formulations P and RA respectively. The following result (see Section~\ref{proof:proposition_P_leq_RA} of the ecompanion for the proof) characterizes the relation between the relaxation bounds of the two formulations.
\begin{proposition}
$Z^*_{\modelP} \leq Z^*_{\modelRA}$.  \label{proposition:P_leq_RA}
\end{proposition}
This result implies that the relaxation bound of formulation~\modelP is always at least as tight as that of formulation RA. This is important because a tighter relaxation bound generally implies that the integer problem can be solved more quickly via branch-and-bound. We will see in Section~\ref{subsec:numerical_experiments_synthetic} that the relaxation bound of \modelP can be substantially tighter than that of \modelRA. In the next section, we discuss one way in which formulation~\modelP can be modified to obtain an even stronger (albeit less tractable) formulation.

\subsection{Formulation \modelPRPT}
\label{subsec:micp_PRPT}

In this section, we propose a modified version of formulation~\modelP that produces a tighter relaxation bound. The key idea in this new formulation is to leverage a recently proposed technique from the paper of \cite{zhen2021extension} called the reformulation-perspectification technique (RPT). RPT is a generalization of the well-known reformulation-linearization technique (RLT) originally proposed by \cite{sherali1990hierarchy}. The basic idea of RPT is to multiply a pair of constraints together, where one constraint is a linear constraint and one is a convex constraint, to generate new constraints that are valid but potentially intractable. These new constraints are then converted into tractable constraints by applying perspectification.

To illustrate the idea, suppose that we are given the following two constraints: 
\begin{align*}
\cb^T \yb \geq d \\
f(x) \leq 0
\end{align*}
where $\yb$ is a vector of decision variables, $x$ is a scalar decision variable, $\cb$ is a vector of the same size as $\yb$, $d$ is a scalar, and $f$ is a convex function. Observe that from the first constraint, we know that $\cb^T \yb - d$ must be nonnegative. Therefore, we can obtain a valid new constraint by multiplying the left and right hand sides of $f(x) \leq 0$ by $\cb^T \yb - d$: 
\begin{equation}
(\cb^T \yb - d) f(x) \leq 0.
\end{equation}
This new constraint is no longer convex. However, we can now apply the perspectification trick to obtain a tractable convex constraint. We multiply and divide the argument of $f$ by $\cb^T \yb - d$:
\begin{equation}
(\cb^T \yb - d) f \left( \frac{ x(\cb^T \yb - d) }{ \cb^T \yb - d} \right) \leq 0.
\end{equation}
Now, letting $\ub$ denote the linearization of $x \cdot \yb$, we can reformulate this as 
\begin{equation}
(\cb^T \yb - d) f \left( \frac{ \cb^T \ub - d x }{ \cb^T \yb - d} \right) \leq 0.
\end{equation}
As in the example in Section~\ref{subsec:micp_P}, this new constraint is a convex constraint, because again we can write the left hand side of the constraint as $\tilde{f}( \cb^T \ub - d x, \cb^T \yb - d)$, where $\tilde{f}(x,t) = t f( x / t)$ for $t \geq 0$ is the perspective function of $f$. This example is a simple example of the procedure; in the paper of \cite{zhen2021extension}, there are a number of more complicated instances shown (for example, where $f$ is a multivariate function).

To apply the RPT technique to formulation~\modelP, let us use as a starting point the constraints 
\begin{align*}
a_i & \geq 0, \\
x_{k,1} + x_{k,1} e^{ -u_k } & \leq 1.
\end{align*}
Note that the first constraint is a valid constraint that must be satisfied by $a_i$, whether it is binary or relaxed to be continuous, while the second constraint is the main constraint from formulation~\modelP prior to perspectification. If we now multiply the second constraint on the left and right by $a_i$, we obtain
\begin{align*}
a_i x_{k,1} + a_i x_{k,1} e^{ - u_k } & \leq a_i
\end{align*}
We now perspectify the second term by multiplying and dividing by $a_i x_{k,1}$ inside the $e$: 
\begin{align*}
a_i x_{k,1} + a_i x_{k,1} e^{ \frac{-a_i u_k x_{k,1} }{a_i x_{k,1}} } & \leq a_i.
\end{align*}
Finally, introducing the new variable $\varphi_{k,i}$ to represent the linearization of $a_i \cdot x_{k,1} \cdot u_k$, and recalling that we had previously introduced $y_{k,i}$ to denote the linearization of $a_i x_{k,1}$, we can re-write this constraint as 
\begin{align*}
y_{k,i} + y_{k,i} e^{ \frac{- \varphi_{k,i}  }{ y_{k,i} }} & \leq a_i,
\end{align*}
which is a convex constraint. We can apply similar steps using the two constraints 
\begin{align*}
1 - a_i & \geq 0, \\
x_{k,1} + x_{k,1} e^{ -u_k } & \leq 1,
\end{align*}
where the first constraint is just a re-arrangement of $a_i \leq 1$. By multiplying the left and right hand side of the second constraint by $1 - a_i$ and following the same steps, we obtain 
\begin{align*}
(x_{k,1} - y_{k,i}) + (x_{k,1} - y_{k,i}) e^{ \frac{ - (w_k - \varphi_{k,i})}{x_{k,1} - y_{k,i}}} \leq 1 - a_i.
\end{align*}
This leads us to the following formulation, which we refer to as formulation~\modelPRPT:
\begin{subequations}
\begin{alignat}{2}
& \underset{ \substack{\ab, \bb, \wb, \xb, \\ \yb, \zb, \varphib} }{\text{maximize}} & \quad & \sum_{k=1}^K \lambda_k x_{k,1} \\
& \text{subject to} & & y_{k,i} + y_{k,i} e^{ \frac{- \varphi_{k,i}  }{ y_{k,i} } } \leq a_i, \quad \forall \ k \in [K], \ i \in [n], \label{prob:PRPT_apersp} \\
& & & (x_{k,1} - y_{k,i}) + (x_{k,1} - y_{k,i}) e^{ \frac{ - (w_k - \varphi_{k,i})}{x_{k,1} - y_{k,i}}} \leq 1 - a_i, \quad \forall \ k \in [K], \ i \in [n], \label{prob:PRPT_oneminusapersp} \\
& & & \varphi_{k,i} = \beta_{k,0} y_{k,i} + \sum_{j=1}^n \beta_{k,j} z_{k,i,j}, \quad \forall \ k \in [K], \ i \in [n], \\
& & & z_{k,j,i} = z_{k,i,j}, \quad \forall \ k \in [K],\ i, j \in [n], i < j, \label{prob:PRPT_z_linearization_1}\\
& & & z_{k,i,i} = y_{k,i}, \quad \forall \ k \in [K], \ i \in [n], \\
& & & z_{k,i,j} \leq y_{k,i}, \quad \forall \ k \in [K], \ i, j \in [n], i < j, \\
& & & z_{k,i,j} \leq y_{k,j}, \quad \forall \ k \in [K], \ i, j \in [n], i < j, \\
& & & z_{k,i,j} \geq y_{k,i} + y_{k,j} - x_{k,1}, \quad \forall \ k \in [K], \ i, j \in [n], i < j, \\
& & & z_{k,i,j} \leq b_{i,j}, \quad \forall \ k \in [K], \ i, j \in [n], i < j, \\
& & & z_{k,i,j} \geq b_{i,j} + y_{k,j} - a_j, \quad \forall \ k \in [K], \ i, j \in [n], i < j, \\
& & & z_{k,i,j} \geq b_{i,j} + y_{k,i} - a_i, \quad \forall \ k \in [K], \ i, j \in [n], i < j, \\
& & & z_{k,i,j} \leq 1 - x_{k,1} - a_i + y_{k,i} - a_j + y_{k,j} + b_{i,j}, \quad \forall \ k \in [K], \ i, j \in [n], i < j, \\
& & & z_{k,i,j} \geq 0, \quad \forall \ k \in [K], \ i, j \in [n], i < j, \label{prob:PRPT_z_linearization_10} \\
& & & b_{j,i} = b_{i,j}, \quad \forall \ i, j \in [n], i < j, \label{prob:PRPT_b_linearization_1} \\
& & & b_{i,i} = a_i, \quad \forall \ i \in [n],\\
& & & b_{i,j} \leq a_i, \quad \forall \ i, j \in [n], i < j, \\
& & & b_{i,j} \leq a_j, \quad \forall \ i, j \in [n], i < j, \\
& & & b_{i,j} \geq a_i + a_j - 1, \quad \forall \ i, j \in [n], i < j, \\
& & & b_{i,j} \geq 0, \quad \forall \ i, j \in [n], i < j, \label{prob:PRPT_b_linearization_6} \\
& & & \text{constraints~\eqref{prob:P_perspective} - \eqref{prob:P_a_definition}}.
\end{alignat}
\label{prob:PRPT}%
\end{subequations}
In the above formulation, constraints~\eqref{prob:PRPT_apersp} and \eqref{prob:PRPT_oneminusapersp} are those obtained by applying RPT. In addition to the new decision variable $\varphi_{k,i}$, the formulation also includes the decision variables $z_{k,i,j}$, which represents the linearization of $a_i \cdot a_j \cdot x_{k,1}$, and $b_{i,j}$, which represents the linearization of $a_i \cdot a_j$. Constraints~\eqref{prob:PRPT_z_linearization_1} - \eqref{prob:PRPT_z_linearization_10} are standard constraints to linearize $a_i a_j x_{k,1}$, while constraints~\eqref{prob:PRPT_b_linearization_1} - \eqref{prob:PRPT_b_linearization_6} are similar constraints to linearize $a_i a_j$. 

We make a couple of observations regarding formulation~\modelPRPT. First, just like formulation~\modelP, this formulation can be represented as a mixed-integer exponential cone program. Second, in terms of strength, observe that the projection of the feasible region of the continuous relaxation of \modelPRPT is contained in the feasible region of the relaxation of P, since the constraints of \modelP are a superset of those in \modelPRPT. Therefore, it follows straightforwardly that the relaxation bound of formulation \modelPRPT, which we denote by $Z^*_{\modelPRPT}$, is at least as tight as that of $Z^*_{\modelP}$. 
\begin{proposition}
$Z^*_{\modelPRPT} \leq Z^*_{\modelP}$. \label{proposition:PRPT_leq_P}
\end{proposition}
We will see in Section~\ref{subsec:numerical_experiments_synthetic} that the relaxation bound of \modelPRPT can be substantially tighter than that of \modelP.

Third, in terms of tractability, formulation~\modelPRPT is significantly more complex than formulation~\modelP. In particular, whereas \modelP is representable using $K$ exponential cone constraints and $O(Kn)$ linear constraints, \modelPRPT requires $O(Kn)$ exponential cone constraints, and $O(Kn^2)$ linear constraints. In our experience with this formulation, it is generally much slower to solve with integrality constraints than formulation~\modelP, and we did not have success with solving this formulation in a reasonable amount of time for the synthetic instances considered in Section~\ref{subsec:numerical_experiments_synthetic}. However, the continuous relaxation of formulation~\modelPRPT can be solved relatively quickly. Thus, this formulation can be useful in some cases in allowing one to quickly obtain a better upper bound than \modelP.

\subsection{Extensions and other formulations}
\label{subsec:micp_extensions}

Before we conclude this section, we comment on a couple of extensions of the models we consider here. First, all of our formulations focus on optimizing the market share of the product. A firm may instead be interested in maximizing expected profit. In the case that profit is a linear function of $\ab$, then all of our formulations can be straightforwardly modified to represent this new objective; we provide more details in Section~\ref{subsec:extra_profit} of the ecompanion. 

Second, an interesting question is whether our formulation techniques in this section can be generalized to handle the case of designing multiple products, i.e., a product line. Indeed, it is possible to generalize formulations \modelRA and \modelP to the problem of selecting attributes for multiple products; we provide the details in Section~\ref{subsec:extra_PLD}.

Third, an important consideration in practice is the reliability of the problem data. In the situation where there are estimation errors in the parameters (the $\lambda_k$ and $\beta_{k,i}$ values), it is possible that a product designed under the assumption of one set of values for these parameters will perform poorly if another set of parameter values is realized. In Section~\ref{sec:robust}, we present two different variants of the logit-based SOCPD problem, based on robust optimization, for addressing uncertainty in the problem data. 

Lastly, we note that all three of our models - \modelRA, \modelP and \modelPRPT - are based on the exponential cone. An interesting question is whether it may be possible to design an exact formulation that is not based on the exponential cone. Surprisingly, it turns out that it is possible to formulate the logit-based SOCPD problem exactly as a mixed-integer second order cone program; we describe this formulation in Section~\ref{sec:SOC} of the ecompanion.

\section{Numerical Experiments}
\label{sec:numerical_experiments}

In this section, we present the results of our numerical experiments. Section~\ref{subsec:numerical_experiments_synthetic} presents the results of our experiments with synthetically generated problem instances, while Section~\ref{subsec:numerical_experiments_real} presents the results of our experiments with instances derived from real conjoint datasets. All of our numerical experiments are implemented in the Julia technical computing language, version 1.5 \citep{bezanson2017julia} using the JuMP package (Julia for Mathematical Programming; see \citealt{dunning2017jump}). All mixed-integer exponential cone programs are solved using Mosek version 10 \citep{aps2021moseksoftware} with a maximum of 8 threads. All of our experiments are conducted on Amazon Elastic Compute Cloud (EC2), on a single instance of type {\tt m6a.48xlarge} (AMD EPYC 7R13 processor, with 192 virtual CPUs and 768 GB of memory).

\subsection{Experiments with synthetic instances}
\label{subsec:numerical_experiments_synthetic}

In our first collection of numerical experiments, we test our approaches on synthetically generated problem instances. We generate these instances as follows. For a fixed number of binary attributes $n$ and number of customer types $K$, we draw an independent uniformly distributed random number $v_{k,i}$ in the interval $[-1, +1]$ for each customer type $k$ and attribute $i$. We then set the partworth $\beta_{k,i}$ of attribute $i$ for customer type $k$ as $\beta_{k,i} = c \cdot v_{k,i}$, where $c$ is a positive constant. For each customer type, we set the utility $\beta_{k,0}$ of the no-purchase option as $\beta_{k,0} = -3$. This choice of the no-purchase option can be interpreted as assuming that when offered a product with no attributes, i.e., $\ab = \zerob$, then the utility of the product for a segment $k$ is $u_k(\ab) = \beta_{k,0} = -3$, which corresponds to a purchase probability of $\sigma(-3) \approx 0.0474$, i.e., a roughly 5\% chance of the customer buying the product. We assume that the probability $\lambda_k$ of each customer type $k$ is set to $1 / K$. 

We vary $n \in \{30,40,50,60,70\}$, $K \in \{10,20,30\}$, and $c \in \{5,10,20\}$. For each combination of $n$ and $K$, we generate 20 collections of values $v_{k,i}$ for each $k \in [K]$ and $i \in [n]$. For each such collection, we vary $c \in \{5, 10, 20\}$. Consequently, this gives rise to $5 \times 3 \times 3 = 45$ sets of 20 problem instances, for a total of $45 \times 20 = 900$ problem instances.

In our first experiment, we solve formulation \modelP on each instance with a time limit of 2 hours. We record the computation time and the optimality gap, which is defined as 
\begin{equation*}
O_{\modelP} = 100\% \times \frac{Z_{\modelP,UB} - Z_{\modelP,LB}}{Z_{\modelP,UB}},
\end{equation*}
where $Z_{\modelP,UB}$ is the best upper bound obtained upon termination of formulation~\modelP and $Z_{\modelP,LB}$ is the best lower bound at termination of \modelP (which corresponds to the best possible integer solution). We compute the average of the computation time, denoted by $T_{\modelP}$, and $O_{\modelP}$ over the twenty instances corresponding to each combination of $n$, $K$ and $c$. 

Table~\ref{table:R1_synthetic_gap_time} displays the results. Overall, we can see that in a large number of cases one can solve formulation \modelP to provable optimality (i.e., the average gap $O_{\modelP}$ is zero) within two hours. In other cases, it is not possible to solve it to provable optimality, but the resulting solution comes with a relatively low optimality gap. For example, for $n = 70$, $K = 30$, $c = 5$, the average optimality gap is 4.95\%. Across all combinations of $n, K, c$ where there is a non-zero average optimality gap, the average optimality gap ranges from 0.19\% to 9.60\%. To put these results into perspective, we attempted to solve some of our instances using exhaustive enumeration in Julia. For instances with $n = 30$, the average time required for exhaustive enumeration was on the order of one hour. For instances with $n = 40$, the time exceeded two hours; given that $\Acal = \{0,1\}^n$, we should expect this time to be about $2^{40 - 30} = 2^{10}$ times larger, i.e., on the order of approximately 1000 hours. In comparison, formulation \modelP can be solved to provable optimality in no more than 2333 seconds (approximately 39 minutes) on average for instances with $n = 30$ and to an average optimality gap of below 7\% within 7026 seconds (just under two hours) for instances with $n = 40$.

With regard to the three instance parameters, $n$, $K$ and $c$, we observe the following behavior. The computation time tends to be decreasing in $c$. There are potentially a couple of reasons for this. First, a larger value of $c$ means that the difference in utility of two distinct attribute vectors $\ab$ and $\ab'$ will be larger. This is potentially beneficial in branch-and-bound, because branching on a single $a_i$ would lead to two child nodes in the branch-and-bound tree where there is a larger difference in the relaxation bounds of the two child nodes, which should generally induce more pruning in the branch-and-bound tree. Second, the constant part of the utility is set to $\beta_{k,0} = -3$ for all customer types. Thus, a smaller $c$ would mean that for more customer types and more product attribute vectors, the utility $u_k(\ab)$ will be negative and reside in the convex (i.e., non-concave) part of the logistic curve. Given that the relaxation of \modelP is a concave maximization problem subject to convex constraints, the relaxation bound is more likely to be loose for problem with small $c$. With regard to $K$, we generally see that computation times are increasing in $K$; this is to be expected, as \modelP becomes larger (i.e., has more variables and constraints) as $K$ increases. With regard to $n$, we generally see that the computation time increases, but then plateaus or decreases slightly when $n$ increases beyond 50. The reason for this is due to the random generation of the problem instances. Recall that the $v_{k,i}$ values are generated in an IID fashion from a uniform distribution over $[-1,+1]$. For a fixed $K$, as $n$ increases, there is a higher probability of there existing at least one attribute where most customer types have a positive partworth, i.e., $v_{k,i} > 0$ for many $k$'s. Setting such attributes to 1 allows the solver to quickly obtain nearly optimal solutions to the problem and close the gap.

\begin{table}
\centering
\begin{tabular}{lllrr} \toprule
$c$ & $n$ & $K$ & $O_{\modelP}$ & $T_{\modelP}$ \\ \midrule
  5 &  30 &  10 & 0.00 & 15.54 \\ 
    5 &  30 &  20 & 0.00 & 202.62 \\ 
    5 &  30 &  30 & 0.00 & 2333.12 \\ 
    5 &  40 &  10 & 0.00 & 44.99 \\ 
    5 &  40 &  20 & 0.37 & 2393.76 \\ 
    5 &  40 &  30 & 6.91 & 7026.27 \\ 
    5 &  50 &  10 & 0.00 & 21.62 \\ 
    5 &  50 &  20 & 0.70 & 2809.89 \\ 
    5 &  50 &  30 & 9.95 & 7213.79 \\ 
    5 &  60 &  10 & 0.00 & 47.54 \\ 
    5 &  60 &  20 & 1.03 & 2317.13 \\ 
    5 &  60 &  30 & 9.60 & 7213.76 \\ 
    5 &  70 &  10 & 0.00 & 18.33 \\ 
    5 &  70 &  20 & 0.45 & 1267.83 \\ 
    5 &  70 &  30 & 4.95 & 6410.63 \\ \bottomrule
\end{tabular}
\qquad
\begin{tabular}{lllrr} \toprule
$c$ & $n$ & $K$ & $O_{\modelP}$ & $T_{\modelP}$ \\ \midrule
   10 &  30 &  10 & 0.00 & 14.96 \\ 
   10 &  30 &  20 & 0.00 & 80.05 \\ 
   10 &  30 &  30 & 0.00 & 1107.93 \\ 
   10 &  40 &  10 & 0.00 & 17.78 \\ 
   10 &  40 &  20 & 0.00 & 399.22 \\ 
   10 &  40 &  30 & 2.50 & 5316.60 \\ 
   10 &  50 &  10 & 0.00 & 17.92 \\ 
   10 &  50 &  20 & 0.19 & 1351.77 \\ 
   10 &  50 &  30 & 6.16 & 7161.57 \\ 
   10 &  60 &  10 & 0.00 & 22.60 \\ 
   10 &  60 &  20 & 0.03 & 709.22 \\ 
   10 &  60 &  30 & 4.88 & 6692.13 \\ 
   10 &  70 &  10 & 0.00 & 16.69 \\ 
   10 &  70 &  20 & 0.00 & 87.16 \\ 
   10 &  70 &  30 & 1.64 & 3310.76 \\ \bottomrule
\end{tabular}
\qquad 
\begin{tabular}{lllrr} \toprule
$c$ & $n$ & $K$ & $O_{\modelP}$ & $T_{\modelP}$ \\ \midrule
   20 &  30 &  10 & 0.00 & 12.23 \\ 
   20 &  30 &  20 & 0.00 & 66.61 \\ 
   20 &  30 &  30 & 0.00 & 614.31 \\ 
   20 &  40 &  10 & 0.00 & 18.42 \\ 
   20 &  40 &  20 & 0.00 & 496.22 \\ 
   20 &  40 &  30 & 1.19 & 4346.31 \\ 
   20 &  50 &  10 & 0.00 & 24.78 \\ 
   20 &  50 &  20 & 0.22 & 817.61 \\ 
   20 &  50 &  30 & 3.40 & 5550.03 \\ 
   20 &  60 &  10 & 0.00 & 38.61 \\ 
   20 &  60 &  20 & 0.00 & 114.48 \\ 
   20 &  60 &  30 & 2.18 & 5115.71 \\ 
   20 &  70 &  10 & 0.00 & 21.44 \\ 
   20 &  70 &  20 & 0.00 & 202.43 \\ 
   20 &  70 &  30 & 0.23 & 2610.99 \\ \bottomrule
\end{tabular}

\caption{Optimality gap and computation time of formulation \modelP as $c$, $n$ and $K$ vary. \label{table:R1_synthetic_gap_time} }
\end{table}

We note that in this experiment we focus exclusively on \modelP, as it generally gave the best performance. In Section~\ref{subsec:R2_PRPT_gap_time}, we present analogous results for models~\modelRA and \modelPRPT. We find that both \modelRA and \modelPRPT are slower to solve than \modelP and result in higher optimality gaps. However, as we will see in the next experiment, the relaxation bound of \modelPRPT can substantially improve on that of \modelP.

In our second experiment, we compare the three formulations -- \modelRA, \modelP and \modelPRPT -- in terms of their continuous relaxations. We solve the continuous relaxation of each of the three formulations for each instance, and then compute the integrality gap, denoted by $I_m$, as 
\begin{equation}
I_m = 100\% \times \frac{Z_{m,rlx} - Z'}{Z'},
\end{equation}
where $Z'$ denotes the objective value of the best integer solution obtained from \modelP after two hours of computation. In addition to $I_m$, we also calculate $T_{rlx,m}$, which is the time required to solve the continuous relaxation of model $m$. We focus on only those instances with $c = 5$, as these instances resulted in the largest separation between the integrality gaps. We note that in a small but not negligible number of cases (79 instance-method pairs out of 900), the value of $Z_{m,rlx}$ obtained was lower than $Z'$ because of numerical precision issues arising from Mosek. Over these 79 instance-method pairs, the error was extremely small, with the relaxation gap ranging between -0.091\% and -0.000016\%. Of the 79 instance-method pairs, 5 corresponded to \modelRA, 23 to \modelP and 51 to \modelPRPT. Due to the relatively small magnitudes of the errors, we treated these negative values as zeros in the calculation of our result metrics. 

Table~\ref{table:R1_synthetic_relaxations} below displays the average integrality gap over the 20 instances for each combination of $(c, n, K)$. From this table, we can see that \modelPRPT generally has the smallest integrality gap, followed by \modelP, and finally \modelRA. In some cases, the difference can be quite large (for example, for $c = 5$, $n = 30$ and $K = 30$, the average integrality gap is about 30\% for \modelRA, compared to about 25\% for \modelP and 18\% for \modelPRPT). This agrees with our theoretical results on the objective values of the continuous relaxations of these formulations (Proposition~\ref{proposition:P_leq_RA} and \ref{proposition:PRPT_leq_P}). In terms of computation time, we note that the relaxation of \modelP is fastest to solve, while \modelRA is about two orders of magnitude slower, and \modelPRPT is an additional 1-2 orders of magnitude slower and can take up to 10 minutes to solve. The edge of \modelP in time is to be expected, as \modelP requires only $K$ exponential cones, while \modelRA requires $4K$ and \modelP requires $2nK + K$. As noted earlier, the integer version of \modelPRPT can be much slower to solve than \modelP; however, these results suggest that \modelPRPT's relaxation could still be useful when one needs to quickly obtain a good upper bound as a complement to a heuristic solution.

\begin{table}
\centering
\begin{tabular}{lllrrrrrr} \toprule
$c$ & $n$ & $K$ & $I_{\modelRA}$ & $I_{\modelP}$ & $I_{\modelPRPT}$ & $T_{rlx,\modelRA}$ & $T_{rlx,\modelP}$ & $T_{rlx,\modelPRPT}$  \\ \midrule
 5 &  30 &  10 & 4.84 & 3.46 & 2.28 & 9.05 & 0.02 & 1.06 \\ 
    5 &  30 &  20 & 17.80 & 13.66 & 9.46 & 9.72 & 0.04 & 15.36 \\ 
    5 &  30 &  30 & 29.57 & 24.55 & 17.99 & 9.89 & 0.08 & 29.95 \\ 
    5 &  40 &  10 & 2.40 & 1.74 & 1.24 & 9.21 & 0.02 & 2.35 \\ 
    5 &  40 &  20 & 10.88 & 8.36 & 6.06 & 10.17 & 0.08 & 41.71 \\ 
    5 &  40 &  30 & 19.83 & 16.33 & 12.51 & 10.68 & 0.17 & 81.65 \\ 
    5 &  50 &  10 & 0.05 & 0.03 & 0.02 & 9.18 & 0.08 & 3.84 \\ 
    5 &  50 &  20 & 5.19 & 4.02 & 3.09 & 10.33 & 0.08 & 91.80 \\ 
    5 &  50 &  30 & 19.32 & 16.58 & 13.79 & 10.69 & 0.17 & 162.83 \\ 
    5 &  60 &  10 & 0.03 & 0.02 & 0.01 & 10.61 & 0.10 & 6.20 \\ 
    5 &  60 &  20 & 2.54 & 2.12 & 1.80 & 10.65 & 0.10 & 224.54 \\ 
    5 &  60 &  30 & 15.71 & 13.73 & 12.04 & 12.47 & 0.23 & 379.66 \\ 
    5 &  70 &  10 & 0.01 & 0.01 & 0.00 & 9.37 & 0.04 & 8.46 \\ 
    5 &  70 &  20 & 0.98 & 0.83 & 0.74 & 10.75 & 0.11 & 388.53 \\ 
    5 &  70 &  30 & 7.93 & 6.98 & 6.19 & 11.97 & 0.39 & 615.05 \\  \bottomrule
\end{tabular}
\caption{Comparison of integrality gaps for continuous relaxations of \modelRA, \modelP and \modelPRPT using synthetic instances with $c = 5$. \label{table:R1_synthetic_relaxations}}
\end{table}

In our final experiment, we test formulation \modelP and compare the quality of its solutions against those of several heuristic approaches. We compare it against several different heuristic approaches, which we summarize below:
\begin{enumerate}
\item \KKDP: This is the dynamic programming (DP) heuristic of \cite{kohli1987heuristic}, which sequentially fixes the elements of $\ab$. Although this heuristic was originally proposed for the problem of product design under a first-choice/max-utility model, with some minor modifications it can also be applied to the logit-based SOCPD problem. 
\item \Greedy: This is the greedy heuristic described in \cite{shi2001optimization}, which involves simply picking the attribute vector $\ab$ that maximizes the weighted average of the customer utilities:
\begin{equation}
\max_{\ab \in \Acal} \sum_{k=1}^K \lambda_k \cdot u_k(\ab).
\end{equation}
When $\Acal = \{0,1\}^n$, this problem can be solved by setting each attribute independently of the others based on the sign of the quantity $\sum_{k=1}^K \lambda_k \cdot \beta_{k,i}$. More generally, it can be solved by formulating a mixed-integer linear program. 
\item \LS: This is a local search heuristic. This heuristic involves starting from a random attribute vector $\ab \in \Acal$, and then moving to a new attribute vector $\ab'$ in the neighborhood $\Ncal(\ab)$ of $\ab$ that leads to the greatest improvement in the objective value; we then repeat this at the new attribute vector, and continue in this way until there is no attribute vector in the neighborhood of the current one that leads to an improvement. We consider a neighborhood $\Ncal(\ab)$ defined as $\Ncal(\ab) = \{ \ab' \in \Acal \mid \|\ab' - \ab\|_1 \leq 1 \}$, in other words, all attribute vectors $\ab'$ that differ from $\ab$ in exactly one coordinate $i \in [n]$. We apply the local search heuristic from ten uniformly randomly generated starting points, and retain the best solution of those ten repetitions. 
\item \GM: This is a novel benchmark heuristic we developed that is based on approximating the objective function of the logit-based SOCPD problem using a weighted geometric mean. We describe this approach in more detail in Section~\ref{sec:extensions_geometric_mean} of the ecompanion.
\end{enumerate}
For formulation~\modelP, we solve it using Mosek with a time limit of 2 hours. 

We execute each of the four heuristics and formulation~\modelP on each of the 900 problem instances. For each problem instance, we identify the solution with the highest objective value and denote its objective function value by $Z'$. For each approach $m$ (one of \KKDP, \Greedy, \LS, \GM, or formulation \modelP), we then compute the post-hoc gap of its solution relative to the best solution for that instance:
\begin{equation}
G_m = 100\% \times \frac{Z' - Z_m}{Z'},
\end{equation}
where $Z_m$ is the objective value attained by approach $m$. We compute the average of $G_m$ over all instances with the same values of $n$, $K$ and $c$, for each approach. The post-hoc gap measures how far the solution of a method is from the best solution found by any method, with a value of 0\% implying that that method's solution is the best (or tied for best) out of all the solutions obtained for that instance.

Table~\ref{table:R1_synthetic_heuristics_gap} below displays the average post-hoc gap of each approach. From this table, we can see that in general, across all the values of $n$, $K$ and $c$, the solution obtained using our mixed-integer exponential cone formulation \modelP tends to be the best one, as it has the lowest average gap, and this gap is very close to zero. Out of the heuristic approaches, \KKDP tends to deliver very poor solutions (average gap over all 900 instances of 16.3\%), followed by \Greedy (average gap of 16.2\%), followed by \LS (average gap of 8\%). We note that the geometric mean heuristic \GM generally tends to give better solutions than all three heuristic in the aggregate, with an average gap of 5.7\% over all 900 instances, although there are many cases where \GM is the weakest of the four heuristic approaches (for example $c = 5$, $n = 30$ $K = 30$); generally, \GM seems to perform best for large $n$ and low $K$. Overall, the main takeaway from these results is that our exact solution approach can lead to solutions that are significantly better than heuristic approaches that do not guarantee global optimality. 

\begin{table}
\centering
\begin{tabular}{rrrrrrrr} \toprule
$c$ & $n$ & $K$ & $G_{\Greedy}$ & $G_{\LS}$ & $G_{\KKDP}$ & $G_{\GM}$ & $G_{\modelP}$ \\ \midrule
 5 &  30 &  10 & 15.74 & 3.91 & 19.25 & 2.37 & 0.00 \\ 
    5 &  30 &  20 & 16.93 & 5.20 & 21.86 & 17.78 & 0.00 \\ 
    5 &  30 &  30 & 17.25 & 4.35 & 17.43 & 20.01 & 0.00 \\ 
    5 &  40 &  10 & 15.64 & 6.17 & 15.73 & 1.03 & 0.00 \\ 
    5 &  40 &  20 & 21.07 & 7.51 & 21.16 & 12.91 & 0.00 \\ 
    5 &  40 &  30 & 19.80 & 7.72 & 18.46 & 18.42 & 0.00 \\ 
    5 &  50 &  10 & 12.21 & 3.56 & 16.48 & 0.00 & 0.00 \\ 
    5 &  50 &  20 & 19.18 & 8.49 & 21.80 & 4.99 & 0.00 \\ 
    5 &  50 &  30 & 16.31 & 6.75 & 15.72 & 15.34 & 1.20 \\ 
    5 &  60 &  10 & 9.43 & 3.49 & 10.45 & 0.00 & 0.00 \\ 
    5 &  60 &  20 & 16.24 & 9.42 & 20.38 & 0.63 & 0.00 \\ 
    5 &  60 &  30 & 17.71 & 7.84 & 18.60 & 9.35 & 0.49 \\ 
    5 &  70 &  10 & 7.74 & 2.49 & 8.52 & 0.00 & 0.00 \\ 
    5 &  70 &  20 & 22.17 & 13.47 & 20.72 & 0.00 & 0.14 \\ 
    5 &  70 &  30 & 19.25 & 12.09 & 22.14 & 3.20 & 0.19 \\ \midrule
   10 &  30 &  10 & 14.14 & 3.54 & 15.73 & 1.74 & 0.00 \\ 
   10 &  30 &  20 & 18.75 & 7.01 & 18.34 & 12.23 & 0.00 \\ 
   10 &  30 &  30 & 18.93 & 7.91 & 18.52 & 20.12 & 0.00 \\ 
   10 &  40 &  10 & 14.05 & 4.30 & 13.69 & 0.00 & 0.00 \\ 
   10 &  40 &  20 & 21.49 & 10.63 & 19.72 & 7.22 & 0.00 \\ 
   10 &  40 &  30 & 20.81 & 11.59 & 17.98 & 13.26 & 0.00 \\ 
   10 &  50 &  10 & 9.97 & 2.50 & 10.72 & 0.00 & 0.00 \\ 
   10 &  50 &  20 & 19.43 & 11.07 & 19.87 & 1.20 & 0.00 \\ 
   10 &  50 &  30 & 17.53 & 10.35 & 17.88 & 10.86 & 0.26 \\ 
   10 &  60 &  10 & 7.47 & 3.51 & 8.05 & 0.00 & 0.00 \\ 
   10 &  60 &  20 & 15.57 & 9.59 & 16.58 & 0.00 & 0.01 \\ 
   10 &  60 &  30 & 19.13 & 12.82 & 18.63 & 2.17 & 0.00 \\ 
   10 &  70 &  10 & 6.49 & 1.50 & 4.56 & 0.00 & 0.00 \\ 
   10 &  70 &  20 & 19.36 & 13.30 & 17.04 & 0.00 & 0.01 \\ 
   10 &  70 &  30 & 20.66 & 16.67 & 22.53 & 0.40 & 0.26 \\ \midrule
   20 &  30 &  10 & 12.95 & 3.28 & 11.15 & 1.58 & 0.00 \\ 
   20 &  30 &  20 & 19.37 & 8.90 & 20.49 & 10.22 & 0.00 \\ 
   20 &  30 &  30 & 19.51 & 10.38 & 18.80 & 42.97 & 0.00 \\ 
   20 &  40 &  10 & 13.23 & 4.50 & 11.12 & 0.00 & 0.00 \\ 
   20 &  40 &  20 & 20.41 & 11.29 & 15.37 & 5.46 & 0.00 \\ 
   20 &  40 &  30 & 21.06 & 13.76 & 19.20 & 10.23 & 0.00 \\ 
   20 &  50 &  10 & 8.94 & 2.50 & 5.50 & 0.00 & 0.00 \\ 
   20 &  50 &  20 & 18.47 & 10.51 & 17.83 & 0.11 & 0.01 \\ 
   20 &  50 &  30 & 19.50 & 13.39 & 19.70 & 7.92 & 0.00 \\ 
   20 &  60 &  10 & 6.83 & 4.50 & 6.60 & 0.00 & 0.00 \\ 
   20 &  60 &  20 & 14.98 & 7.76 & 16.25 & 0.00 & 0.02 \\ 
   20 &  60 &  30 & 20.21 & 14.72 & 20.61 & 0.30 & 0.01 \\ 
   20 &  70 &  10 & 6.05 & 3.50 & 5.76 & 0.00 & 0.00 \\ 
   20 &  70 &  20 & 17.78 & 11.50 & 17.91 & 0.00 & 0.01 \\ 
   20 &  70 &  30 & 20.62 & 16.65 & 19.72 & 0.00 & 0.08 \\ \midrule
   \multicolumn{3}{c}{(Mean)} & 16.23 & 8.13 & 16.32 & 5.65 & 0.06 \\
   \multicolumn{3}{c}{(Median)} & 17.40 & 8.85 & 16.81 & 0.00 & 0.00 \\ \bottomrule
   \end{tabular}
   \caption{Comparison of post-hoc gap of heuristic approaches and exact approach (from solving \modelP) on synthetic instances. \label{table:R1_synthetic_heuristics_gap} }
\end{table}

Besides the post-hoc gap, it is also helpful to compare the exact solution of formulation \modelP with the heuristics in terms of computation time. Due to space considerations, these results are reported in Table~\ref{table:R1_synthetic_time} in Section~\ref{subsec:R1_synthetic_time} of the ecompanion, which compares the approaches in terms of average computation time, where the average is taken over the twenty instances for a fixed $n$, $K$, $c$ combination. From this table, our formulation~\modelP requires the most time, while the \KKDP, \Greedy and \LS heuristics are extremely fast, requiring no more than a second in all cases. Although solving \modelP requires more time than the heuristics, we believe that the additional runtime is justified in light of the fact that \modelP produces solutions for which the level of suboptimality (i.e., the optimality gap) is known, which is not the case for \KKDP, \Greedy or \LS. The geometric mean approach \GM requires significantly less time compared to formulation~\modelP; across all of the instances, we were able to solve the geometric mean formulation~\eqref{prob:GM_MICONVP} to provable optimality in under two minutes on average; across all 900 instances, the largest  time we observed was 996 seconds (just over 16 minutes). %

\subsection{Experiments with instances based on real conjoint datasets}
\label{subsec:numerical_experiments_real}

In our second set of numerical experiments, we test our approaches using instances built with logit models estimated from real conjoint datasets. We use four different data sets: \toubia, a dataset on preferences for laptop bags produced by Timbuk2 from \cite{toubia2003fast} (see also \citealt{belloni2008optimizing}, \citealt{bertsimas2017robust,bertsimas2019exact}, which also use this data set for profit-based product line design); \bank, a dataset on preferences for credit cards from \cite{allenby1995using} (accessed through the \texttt{bayesm} package for R; see \citealt{rossi2019bayesm}); \candidate, a dataset on preferences for a hypothetical presidential candidate from \cite{hainmueller2014causal}; \immigrant, a dataset on preferences for a hypothetical immigrant from \cite{hainmueller2014causal}. The high-level characteristics of each dataset are summarized in Table~\ref{table:real_data_summary} below, and we provide additional details on the datasets in Section~\ref{subsec:R1_real_data_details}.

We note that for some of these datasets, the product design problem is of a more hypothetical nature. For example, for \candidate, the problem is to ``design'' a political candidate maximizing the share of voters who would vote for that candidate. Similarly, for \immigrant, the problem is to ``design'' an ideal immigrant that would maximize the fraction of people who would support granting admission to such an immigrant. Clearly, it is not possible to create a political candidate or immigrant with certain characteristics. Despite this, we believe that identifying what an optimal ``product'' would be for these data sets, and what share-of-choice such a product would achieve, would still be insightful. Notwithstanding these concerns, these datasets are still valuable from the perspective of verifying that our optimization methodology can solve problem instances derived from real data. 

\begin{table}
\centering 
\begin{tabular}{lllll} \toprule
Dataset & Respondents & Attributes & Attribute Levels & $n$ \\ \midrule
\bank & 946 & 7 & $4 \times 4 \times 3 \times 3 \times 3 \times 2 \times 2$ & 14 \\
\candidate & 311 & 8 & $6 \times 2 \times 6 \times 6 \times 6 \times 6 \times 6 \times 2$   & 32 \\
\immigrant & 1396 & 9 & $7 \times 2 \times 10 \times 3 \times 11 \times 4 \times 4 \times 5 \times 4$ & 41 \\
\toubia & 330 & 10 &  $7 \times 2 \times 2 \times 2 \times 2 \times 2 \times 2 \times 2 \times 2 \times 2$ & 15 \\ \bottomrule
\end{tabular}
\caption{Summary of real conjoint datasets used in Section~\ref{subsec:numerical_experiments_real}. The column ``Attributes'' indicates the number of attributes, and ``Attribute Levels'' indicates the structure of each attribute (e.g., $2 \times 3 \times 5$ indicates that the product has one attribute with two levels, followed by one with three levels, followed by one with five levels). The column labelled $n$ indicates the resulting number of binary attributes when the dataset is used to formulate the logit-based SOCPD problem.  \label{table:real_data_summary} }
\end{table}

For each data set, we develop two different types of logit models, which we summarize below. 
\begin{enumerate}
\item \emph{Latent-class logit}: For each dataset, we estimate a latent-class (LC) multinomial logit with a finite number of classes $K$. We estimate each model using a custom implementation of the expectation-maximization (EM) algorithm \citep{train2009discrete}. For each dataset, we run the EM algorithm from five randomly chosen starting points, and retain the model with the lowest log likelihood. To ensure numerical stability, we impose the constraint $-10 \leq \beta_{k,i} \leq 10$ for each $i$ in the M step of the algorithm. We vary the number of classes $K$ in the set $\{5, 10, 15, 20, 30, 40, 50\}$. Thus, in the associated logit-based SOCPD instance, each customer class corresponds to one of the customer types and the customer type probability $\lambda_k$ is the class $k$ probability estimated via EM. Although we do not focus on tuning $K$ here, we provide additional details on which value of $K$ would be preferred from a model selection standpoint using standard metrics (such as the Akaike information criterion and the Bayesian information criterion) in Section~\ref{appendix:numerics_tuning_K} of the ecompanion.
\item \emph{Hierarchical Bayes}: For each dataset, we estimate a mixture multinomial logit (MMNL) model with a multivariate normal mixture distribution using the hierarchical Bayesian (HB) approach; we use a standard specification with normal-inverse Wishart second stage priors (see Section~\ref{appendix:numerics_HB_specification} of the ecompanion for more details). We estimate this model using Markov chain Monte Carlo (MCMC) via the \texttt{bayesm} package in R \citep{rossi2019bayesm}. We simulate 50,000 draws from the posterior distribution of $(\beta_{r,1}, \dots, \beta_{r,n})$ for each respondent $r$, and thin the draws to retain every 100th draw. Of those draws, we retain the last $J = 100$ draws, which we denote as $(\beta^j_{r,1}, \dots, \beta^j_{r,n})$, where $j \in \{1,\dots, J\}$, and we compute the average partworth vector $(\beta_{k,1},\dots,\beta_{k,n})$ as 
\begin{equation}
(\beta_{r,1}, \dots, \beta_{r,n}) = ( \frac{1}{J} \sum_{j=1}^{J} \beta^j_{r,1}, \dots, \frac{1}{J} \sum_{j=1}^{J} \beta^j_{r,n}).
\end{equation}
This approach leads to an estimate of the partworths for each of the respondents. In the corresponding logit-based SOCPD instance, the number of customer types $K$ corresponds to the number of respondents, and the probability of each customer type $k$ is $1/K$.
\end{enumerate}
We additionally visualize the distributions of the partworths for all of the logit models that we consider in Section~\ref{appendix:partworth_distributions} of the ecompanion.

Before continuing, we note that there may be other approaches for defining a mixture of logits model. (For example, given an estimate of the mean and covariance matrix of a normal mixture distribution defining a mixture logit model, one could sample a set of $K$ partworth vectors and use those as the set of customer types, with each $\lambda_k = 1/K$.) We emphasize that our goal is not to advocate for one approach over another. The estimation approaches described here are simply for the purpose of obtaining problem instances that are of a realistic scale and correspond to real data. We note that our optimization approach is agnostic to how the customer choice model is constructed and is compatible with any estimation approach, so long as it results in a finite set of customer types that each follow a logit model of choice.

For each dataset, we define the set $\Acal$ to be the set of all binary vectors of size $n$ that respect the attribute structure of the dataset; in particular, for attributes that are not binary, we introduce constraints of the form $\sum_{i \in S} a_i \leq 1$ as appropriate (cf. constraints~\eqref{eq:attribute1_constraint} and \eqref{eq:attribute2_constraint} in Section~\ref{subsec:extra_Acal}). For \immigrant, we also follow \cite{hainmueller2014causal} in not allowing certain combinations of attributes (for example, it is not possible for a hypothetical immigrant to be a doctor and have only a high school education). We briefly describe the constraints for \immigrant in Section~\ref{appendix:numerics_immigrant_constraints} of the ecompanion.

With regard to the no-purchase option, recall from Section~\ref{subsec:model_problem_definition} that the constant part of each customer's utility function, $\beta_{k,0}$, can be thought of as the negative of the utility of the no-purchase option. None of the four data sets include explicit information on the no-purchase option, and they did not include any tasks where respondents were able to choose between the no-purchase option and a hypothetical product. Thus, to define the utility of the no-purchase option, we take a different approach, where we assume that in each problem instance, each customer can choose from three different competitive offerings which are defined using the same attributes as the product that is to be designed. This is a standard assumption in the product design and product line design literature \citep{belloni2008optimizing,bertsimas2017robust,bertsimas2019exact}. More details on this calculation are provided in Section~\ref{subsec:R1_real_data_details} and we provide the specific details of the competitive offerings for each data set in Section~\ref{appendix:numerics_competitive_offerings} of the ecompanion. 

For each real data instance, we test the \Greedy, \LS and \GM heuristics. We solve formulation~\modelP using Mosek, with a time limit of two hours. 

Table~\ref{table:real} shows the computation time and the objective value of all of the different methods for all four datasets. From this table, we can see that all of the LC (latent class logit) instances can be solved to complete optimality using formulation~\modelP within 16 seconds, while all of the HB instances are solved within ten minutes. As in our synthetic experiments, the \Greedy and \LS heuristics are the fastest, requiring under a second to execute in all cases. Although \Greedy and \LS sometimes obtain the optimal solution, this is not always the case, and in some cases there can be a large gap between these heuristic solutions and the optimal solution; to focus on one example, for \immigrant with LC and $K = 30$, \Greedy and \LS are about 12\% and 8\% suboptimal, respectively). With regard to the geometric mean approach, we find that Mosek is able to solve all of the instances very quickly (within 6 seconds in all cases), but the solutions obtained from \GM perform worse in these datasets than in the synthetic datasets considered in the previous section and exhibit higher suboptimality gaps than \Greedy and \LS. These results again highlight the value of a provably optimal approach to the logit-based SOCPD problem.

\begin{table}
\centering
\begin{tabular}{lllcccccccc} \toprule
& & & \multicolumn{4}{c}{Objective Value} & \multicolumn{4}{c}{Computation Time (s)} \\
Dataset & Model & $K$ & \Greedy & \LS & \GM & \modelP &  \Greedy & \LS & \GM & \modelP \\ \midrule
  \bank & LC & 5 & 0.737 & 0.742 & 0.675 & 0.742 & 0.01 & 0.00 & 0.12 & 0.19 \\ 
   & LC & 10 & 0.743 & 0.749 & 0.685 & 0.749 & 0.00 & 0.00 & 0.05 & 0.29 \\ 
   & LC & 15 & 0.752 & 0.752 & 0.682 & 0.752 & 0.00 & 0.00 & 0.06 & 0.54 \\ 
   & LC & 20 & 0.719 & 0.719 & 0.687 & 0.719 & 0.00 & 0.00 & 0.08 & 0.49 \\ 
   & LC & 30 & 0.736 & 0.764 & 0.637 & 0.764 & 0.00 & 0.00 & 0.14 & 0.97 \\ 
   & LC & 40 & 0.757 & 0.757 & 0.665 & 0.757 & 0.00 & 0.00 & 0.17 & 1.27 \\ 
   & LC & 50 & 0.746 & 0.749 & 0.642 & 0.749 & 0.00 & 0.00 & 0.13 & 1.70 \\ 
  & HB & 946 & 0.817 & 0.817 & 0.812 & 0.817 & 0.01 & 0.01 & 2.23 & 45.52 \\ \midrule
  \candidate & LC & 5 & 0.504 & 0.471 & 0.509 & 0.626 & 0.00 & 0.00 & 0.05 & 1.14 \\ 
   & LC & 10 & 0.573 & 0.597 & 0.637 & 0.694 & 0.00 & 0.00 & 0.09 & 1.49 \\ 
   & LC & 15 & 0.637 & 0.616 & 0.651 & 0.670 & 0.00 & 0.00 & 0.07 & 2.13 \\ 
   & LC & 20 & 0.563 & 0.628 & 0.574 & 0.705 & 0.00 & 0.00 & 0.11 & 2.55 \\ 
   & LC & 30 & 0.549 & 0.567 & 0.534 & 0.627 & 0.00 & 0.00 & 0.16 & 6.12 \\ 
   & LC & 40 & 0.585 & 0.671 & 0.537 & 0.671 & 0.00 & 0.00 & 0.22 & 9.48 \\ 
   & LC & 50 & 0.710 & 0.669 & 0.680 & 0.710 & 0.01 & 0.00 & 0.34 & 9.83 \\ 
   & HB & 311 & 0.829 & 0.851 & 0.851 & 0.852 & 0.01 & 0.03 & 1.05 & 42.65 \\ \midrule
  \immigrant & LC & 5 & 0.555 & 0.670 & 0.687 & 0.689 & 0.00 & 0.00 & 0.05 & 0.63 \\ 
   & LC & 10 & 0.688 & 0.718 & 0.688 & 0.738 & 0.01 & 0.00 & 0.08 & 1.91 \\ 
   & LC & 15 & 0.683 & 0.631 & 0.393 & 0.726 & 0.00 & 0.00 & 0.08 & 2.66 \\ 
   & LC & 20 & 0.706 & 0.696 & 0.552 & 0.756 & 0.01 & 0.00 & 0.18 & 4.25 \\ 
   & LC & 30 & 0.595 & 0.622 & 0.467 & 0.675 & 0.01 & 0.00 & 0.20 & 9.02 \\ 
   & LC & 40 & 0.724 & 0.692 & 0.344 & 0.724 & 0.01 & 0.00 & 0.29 & 9.84 \\ 
   & LC & 50 & 0.713 & 0.689 & 0.628 & 0.731 & 0.01 & 0.00 & 0.29 & 15.81 \\ 
   & HB & 1396 & 0.828 & 0.851 & 0.846 & 0.865 & 0.01 & 0.22 & 6.34 & 549.31 \\ \midrule
  \toubia & LC & 5 & 0.519 & 0.519 & 0.510 & 0.519 & 0.00 & 0.00 & 0.04 & 0.14 \\ 
   & LC & 10 & 0.543 & 0.543 & 0.536 & 0.543 & 0.00 & 0.00 & 0.06 & 0.34 \\ 
   & LC & 15 & 0.551 & 0.567 & 0.430 & 0.567 & 0.00 & 0.00 & 0.08 & 0.48 \\ 
   & LC & 20 & 0.556 & 0.557 & 0.556 & 0.557 & 0.00 & 0.00 & 0.11 & 0.80 \\ 
   & LC & 30 & 0.596 & 0.620 & 0.436 & 0.620 & 0.01 & 0.00 & 0.16 & 1.21 \\ 
   & LC & 40 & 0.579 & 0.579 & 0.560 & 0.579 & 0.00 & 0.00 & 0.24 & 1.78 \\ 
   & LC & 50 & 0.628 & 0.628 & 0.446 & 0.628 & 0.00 & 0.00 & 0.25 & 2.19 \\ 
   & HB & 330 & 0.644 & 0.644 & 0.644 & 0.644 & 0.00 & 0.01 & 1.54 & 16.75 \\ \bottomrule
\end{tabular}
\caption{Results for numerical experiment with real data. \label{table:real} }
\end{table}

In addition to the performance of the different methods, it is also interesting to examine the optimal solutions. Table~\ref{table:candidate_solutions} visualizes the optimal solution for the \candidate dataset for the LC model with $K = 20$ segments. The table also shows the three outside options/competitive offerings that were defined for this dataset. In addition, the table also shows the structure of the solution obtained by \Greedy, which finds the vector $\ab$ in $\Acal$ that maximizes $\sum_{k=1}^K \lambda_k u_k(\ab)$. 

From this table, we can see that the optimal solution matches some of the outside options on certain attributes (such as income and profession), but differs on some (for example, age). In addition, while the optimal solution does match the heuristic on many attributes, it differs on a couple of key attributes, namely race/ethnicity (the optimal candidate is Black, while the heuristic candidate is Asian American) and gender (the optimal candidate is male, while the heuristic candidate is female). While this may appear to be a minor difference, it results in a substantial difference in market share: the heuristic candidate attracts a share of 0.563, while the optimal candidate attracts a share of 0.705, which is an improvement of 25\%. This illustrates that intuitive solutions to the logit-based product design problem can be suboptimal, and demonstrates the value of a principled optimization-based approach to this problem.

\begin{table}

\centering
\scriptsize

\begin{tabular}{|l|c|c|c|c|c|} \hline
Attribute & Outside  & Outside  &  Outside  & Optimal  & \Greedy  \\ 
& Option 1 & Option 2 & Option 3 & Solution & Solution \\ \hline
Age: 36 & \cellcolor{black!25} &  &  &  &  \\ \hline
Age: 45 &  &  &  & \cellcolor{black!25} & \cellcolor{black!25} \\ \hline
Age: 52 &  & \cellcolor{black!25} &  &  &  \\ \hline
Age: 60 &  &  &  &  &  \\ \hline
Age: 68 &  &  & \cellcolor{black!25} &  &  \\ \hline
Age: 75 &  &  &  &  &  \\ \hline
Military Service: Did not serve & \cellcolor{black!25} &  & \cellcolor{black!25} &  &  \\ \hline
Military Service: Served &  & \cellcolor{black!25} &  & \cellcolor{black!25} & \cellcolor{black!25} \\ \hline
Religion: None & \cellcolor{black!25} & \cellcolor{black!25} &  & \cellcolor{black!25} & \cellcolor{black!25} \\ \hline
Religion: Jewish &  &  &  &  &  \\ \hline
Religion: Catholic &  &  & \cellcolor{black!25} &  &  \\ \hline
Religion: Mainline protestant &  &  &  &  &  \\ \hline
Religion: Evangelical protestant &  &  &  &  &  \\ \hline
Religion Mormon &  &  &  &  &  \\ \hline
College: No BA &  &  &  &  &  \\ \hline
College: Baptist college &  &  &  &  &  \\ \hline
College: Community college & \cellcolor{black!25} &  &  &  &  \\ \hline
College: State university &  &  &  &  &  \\ \hline
College: Small college &  & \cellcolor{black!25} &  &  &  \\ \hline
College: Ivy League university &  &  & \cellcolor{black!25} & \cellcolor{black!25} & \cellcolor{black!25} \\ \hline
Income: 32K & \cellcolor{black!25} &  &  &  &  \\ \hline
Income: 54K &  &  &  &  &  \\ \hline
Income: 65K  &  & \cellcolor{black!25} &  &  &  \\ \hline
Income: 92K &  &  &  &  &  \\ \hline
Income: 210K &  &  & \cellcolor{black!25} & \cellcolor{black!25} & \cellcolor{black!25} \\ \hline
Income 5.1M &  &  &  &  &  \\ \hline
Profession: Business owner & \cellcolor{black!25} &  &  &  &  \\ \hline
Profession: Lawyer &  &  & \cellcolor{black!25} & \cellcolor{black!25} & \cellcolor{black!25} \\ \hline
Profession: Doctor &  &  &  &  &  \\ \hline
Profession: High school teacher &  & \cellcolor{black!25} &  &  &  \\ \hline
Profession: Farmer &  &  &  &  &  \\ \hline
Profession: Car dealer &  &  &  &  &  \\ \hline
Race/Ethnicity: White & \cellcolor{black!25} &  &  &  &  \\ \hline
Race/Ethnicity: Native American &  &  &  &  &  \\ \hline
Race/Ethnicity: Black &  &  &  & \cellcolor{black!25} &  \\ \hline
Race/Ethnicity: Hispanic &  & \cellcolor{black!25} &  &  &  \\ \hline
Race/Ethnicity: Caucasian &  &  &  &  &  \\ \hline
Race/Ethnicity: Asian American &  &  & \cellcolor{black!25} &  & \cellcolor{black!25} \\ \hline
Gender: Male & \cellcolor{black!25} & \cellcolor{black!25} &  & \cellcolor{black!25} &  \\ \hline
Gender: Female &  &  & \cellcolor{black!25} &  & \cellcolor{black!25} \\ \hline

\end{tabular}

\caption{Attributes of outside options, optimal solution and heuristic solution for \candidate LC-MNL model with $K = 20$ segments. \label{table:candidate_solutions} }

\end{table}

\section{Conclusions}
\label{sec:conclusions}

In this paper, we have studied the logit-based share-of-choice product design problem. While we have showed that this problem is theoretically intractable even in the simplest case of two customer types and is moreover NP-Hard to approximate to within a reasonable factor, we nevertheless showed how it is possible to transform this problem into a mixed-integer convex optimization problem and in particular, a mixed-integer (exponential) cone program, which allows us to leverage cutting edge solvers that can handle these types of problems such as Mosek. Our numerical experiments show how our approach can obtain high quality solutions to large instances, whether generated synthetically or from real conjoint data, within reasonable time limits. To the best of our knowledge, this is the first methodology for solving the logit-based share-of-choice product design problem to provable optimality. 

There are a number of interesting directions for future research. In this paper, we have focused on exact methods for solving the SOCPD problem to provable optimality. A different direction is to consider approximation algorithms for the SOCPD problem. In Section~\ref{sec:approximation_algorithm} of the ecompanion, we present an approximation algorithm for the problem that has an exponential dependence on the number of customer types $K$ and on the input size (through a parameter that measures the range of partworth values), but is otherwise polynomial in $n$. Whether it is possible to eliminate this dependence on input size and to obtain a fully polynomial time approximation scheme (FPTAS) when $K$ is assumed to be constant (mirroring known results for assortment optimization under the mixture of MNL model; see \citealt{desir2022capacitated}) is an open question. 

Another interesting direction is to generalize the methodology here to problems involving logit probabilities outside of the product design space. Problems that involve logit probabilities appear in many other contexts. One example is the paper of \cite{jagabathula2020conditional}, which proposes a method for learning mixture of logit models from transaction data based on the Frank Wolfe method, where the key step in the method is to solve a subproblem that involves a weighted combination of logit probabilities. Another example is the paper of \cite{liu2022planning}, which proposes a method for planning urban bike lanes that also involves maximizing a sum of logit probabilities. 

Lastly, another interesting direction is to consider uncertainty in the logit-based SOCPD problem. While we have considered two robust optimization formulations of the uncertain logit-based SOCPD problem in Section~\ref{sec:robust} of the ecompanion, both of these approaches are limited to uncertainty in the partworth vectors $\betab_1,\dots, \betab_K$ only. An interesting question here is whether there exist other approaches that can address uncertainty in the type probabilities $\lambda_1,\dots, \lambda_K$ and the partworth vectors $\betab_1,\dots, \betab_K$ simultaneously.

\section*{Acknowledgments}
The authors thank the authors of \cite{toubia2003fast}, \cite{hainmueller2014causal} and \cite{allenby1995using} for making their data sets available, which were used in the experiments in Section~\ref{subsec:numerical_experiments_real}. The authors also thank the department editor George Shanthikumar, the associate editor and three anonymous referees for their helpful and constructive feedback that has significantly improved the quality of the paper.

\bibliographystyle{plainnat}
\bibliography{logit_SOCPD_lit_v2}

\ECSwitch

\section{Omitted proofs}

\subsection{Proof of Theorem~\ref{theorem:NPHard_Keq2} (NP-Hardness when $K = 2$)}
\label{proof:theorem_NPHard_Keq2}

To show this result, we will show that the partition problem, a well-known NP-complete problem \citep{garey2002computers}, can be reduced to the decision form of the logit-based SOCPD problem. The partition problem can be stated as follows:
\vspace{1em}

\begin{center}
\fbox{
\parbox{0.85\textwidth}{
\textbf{Partition:} \\
\textbf{Inputs}: 
\begin{itemize}
\item Integer $n$; 
\item integers $c_1,\dots, c_n$.
\end{itemize}
\textbf{Question}: Does there exist a set $S \subseteq [n]$ such that $\sum_{i \in S} c_i = \sum_{i \notin S} c_i$? 
}
}
\end{center}
\vspace{1em}

The decision form of the logit-based SOCPD problem can be stated as follows:\\

\begin{center}
\fbox{
\parbox{0.85\textwidth}{
\textbf{Logit-based SOCPD problem with $K = 2$ (decision form):} \\
\textbf{Inputs}: 
\begin{itemize} 
\item Integer $n$; 
\item utility parameters $\beta_{1,0}, \dots, \beta_{1,n}, \beta_{2,0}, \dots, \beta_{2,n}$; 
\item customer type probabilities $\lambda_1, \lambda_2 \geq 0$ such that $\lambda_1 + \lambda_2 = 1$; 
\item target share-of-choice value $\theta$.
\end{itemize}
\textbf{Question}: Does there exist an $\ab \in \Acal = \{0,1\}^n$ such that 
$$\lambda_1 \sigma( u_1(\ab)) + \lambda_2 \sigma(u_2(\ab)) \geq \theta$$
is satisfied? 
}
}
\end{center}

\vspace{1em}

Given an instance of the partition problem, we construct an instance of the decision form of the logit-based SOCPD problem such that the answer to the partition problem is yes if and only if the answer to the decision form of the logit-based SOCPD problem is yes.

Let $c_1,\dots, c_n$ be the sizes of the $n$ items in the partition problem. Let $T = \sum_{i=1}^n c_i$ be the total of all of the sizes. Note that the equality $\sum_{i \in S} c_i = \sum_{i \notin S} c_i$ implies
\begin{align*}
\sum_{i \in S} c_i & = \sum_{i \notin S} c_i \\
\Rightarrow \sum_{i \in S} c_i + \sum_{i \in S} c_i & = \sum_{i \in S} c_i + \sum_{i \notin S} c_i \\
\Rightarrow 2 \sum_{i \in S} c_i & = T \\
\Rightarrow \sum_{i \in S} c_i & = T/2.
\end{align*}
Thus, a set $S$ answers the partition problem in the affirmative if and only if $\sum_{i \in S} c_i = T/2$ if and only if $\sum_{i \notin S} c_i = T/2$. 

Consider an instance of the decision form of the logit-based SOCPD problem defined as follows. Let $\lambda_1 = \lambda_2 = 0.5$. Let the number of attributes be the same as the number of items $n$, and let $\Acal = \{0,1\}^n$. Let $p_U = 0.9$ and $p_L = 0.1$, and define $q_U = \log \frac{p_U}{1-p_U}$ and $q_L = \log \frac{p_L}{1-p_L}$ as the logits corresponding to $p_U$ and $p_L$ respectively. Define the utility parameters as follows:
\begin{align*}
\beta_{1,0} & = q_L + (1 - T/2) \cdot (q_U - q_L), \\
\beta_{1,i} & = (q_U - q_L) \cdot c_i, \quad \forall \ i \in [n],\\
\beta_{2,0} & = q_L + (T/2 + 1) \cdot (q_U - q_L), \\
\beta_{2,i} & = -(q_U - q_L) \cdot c_i, \quad \forall \ i \in [n].
\end{align*}
The utility functions $u_1, u_2: \Acal \to \Rbb$ are then
\begin{align*}
u_1(\ab) & = \beta_{1,0} + \sum_{i=1}^n \beta_{1,i} a_i \\
& = q_L + (1 - T/2) \cdot (q_U - q_L) + \sum_{i=1}^n (q_U - q_L) \cdot c_i \cdot a_i \\
& = q_L + (q_U - q_L) \cdot \left[ \sum_{i=1}^n c_i a_i - T/2 + 1 \right], \\
u_2(\ab) & = \beta_{2,0} + \sum_{i=1}^n \beta_{2,i} a_i \\
& = q_L + (T/2 + 1) \cdot (q_U - q_L) + \sum_{i=1}^n -(q_U - q_L) \cdot c_i \cdot a_i \\
& = q_L + (q_U - q_L) \cdot \left[ \sum_{i=1}^n - c_i a_i + T/2 + 1 \right].
\end{align*}
Finally, let $\theta = p_U = 0.9$. 

We now show that the answer to the partition problem is yes if and only if the answer to the logit-based SOCPD problem with $K = 2$ is yes. \\

\noindent \emph{Partition is yes $\Rightarrow$ Logit-based SOCPD is yes}. To prove this direction of the equivalence, let $S$ be the set for which $\sum_{i \in S} c_i = \sum_{i \notin S} c_i$. As discussed earlier, this implies that $\sum_{i \in S} c_i = T/2$ and $\sum_{i \notin S} c_i = T/2$. Let the product vector $\ab = (a_1,\dots, a_n)$ be defined as
\begin{align*}
a_i = \Ibb\{ i \in S\}.
\end{align*}
Observe now that:
\begin{align*}
u_1(\ab) & = q_L + (q_U - q_L) \cdot \left[ \sum_{i=1}^n c_i a_i - T/2 + 1 \right] \\
& = q_L + (q_U - q_L) \cdot \left[ \sum_{i \in S} c_i - T/2 + 1 \right] \\
& = q_L + (q_U - q_L) \cdot \left[ T/2 - T/2 + 1 \right] \\
& = q_L + (q_U - q_L) \cdot 1 \\
& = q_U, \\
u_2(\ab) & = q_L + (q_U - q_L) \cdot \left[ \sum_{i=1}^n - c_i a_i + T/2 + 1 \right] \\
& = q_L + (q_U - q_L) \cdot \left[ \sum_{i \in S} - c_i  + T/2 + 1 \right] \\
& = q_L + (q_U - q_L) \cdot \left[ -T/2  + T/2 + 1 \right] \\
& = q_L + (q_U - q_L) \\
& = q_U. 
\end{align*}
This implies that the objective value of $\ab$ is 
\begin{align*}
& \lambda_1 \sigma( u_1(\ab)) + \lambda_2 \sigma(u_2(\ab)) \\
& = 0.5 \cdot \sigma(q_U) + 0.5 \cdot \sigma(q_U) \\
& = (0.5)(0.9) + (0.5)(0.9) \\
& = 0.9,
\end{align*}
which implies that the answer to the decision form of the logit-based SOCPD problem is yes, as required.\\

\noindent \emph{Partition is no $\Rightarrow$ Logit-based SOCPD is no}. To prove the other direction of the equivalence, let $\ab$ be any product attribute vector. We need to show that the objective value of $\ab$ in the logit-based SOCPD problem is strictly less than 0.9. To see this, observe that if we define $S = \{i \in [n] \mid a_i = 1\}$, we obtain a subset of $[n]$. Since the answer to the partition problem is no, we know that $\sum_{i \in S} c_i \neq \sum_{i \notin S} c_i$. This is equivalent to $\sum_{i \in S} c_i \neq T/2$. 

There are now two possible cases to consider for where $\sum_{i \in S} c_i$ is in relation to $T/2$. If $\sum_{i \in S} c_i > T/2$, then because the $c_i$'s are integers, this means that $\sum_{i \in S} c_i \geq T/2 + 1/2$. This implies that the utility of segment 2 for product vector $\ab$ can be upper bounded as follows:
\begin{align*}
u_2(\ab) & = q_L + (q_U - q_L) \cdot \left[ \sum_{i=1}^n - c_i a_i + T/2 + 1 \right] \\
& =  q_L + (q_U - q_L) \cdot \left[ \sum_{i \in S} - c_i + T/2 + 1 \right] \\
& \leq q_L + (q_U - q_L) \cdot \left[ -T/2 - 1/2 + T/2 + 1 \right] \\
& = q_L + (q_U - q_L) \cdot (1/2) \\
& = (q_L + q_U)/2 \\
& = 0
\end{align*}
which implies that the objective value of $\ab$ is bounded from above as 
\begin{align*}
& \lambda_1 \sigma( u_1(\ab)) + \lambda_2 \sigma(u_2(\ab)) \\
& \leq 0.5 \cdot 1 + 0.5 \cdot \sigma(0) \\
& = 0.5 + (0.5)(0.5) \\
& = 0.75 \\
& < 0.9.
\end{align*}

Alternatively, if $\sum_{i \in S} c_i < T/2$, then we know that $\sum_{i \in S} c_i \leq T/2 - 1/2$. This implies that the utility of segment 1 for $\ab$ can be upper bounded as follows:
\begin{align*}
u_1(\ab) & = q_L + (q_U - q_L) \cdot \left[ \sum_{i=1}^n c_i a_i - T/2 + 1 \right] \\
& = q_L + (q_U - q_L) \cdot \left[ \sum_{i \in S} c_i - T/2 + 1 \right] \\
& \leq q_L + (q_U - q_L) \cdot \left[ T/2 - 1/2 - T/2 + 1 \right] \\
& = (q_L + q_U)/2 \\
& = 0, 
\end{align*}
which again implies that the objective value of $\ab$ is bounded from above as 
\begin{align*}
& \lambda_1 \sigma( u_1(\ab)) + \lambda_2 \sigma(u_2(\ab)) \\ 
& \leq \lambda_1 \sigma( 0 ) + \lambda_2 \cdot 1 \\
& = (0.5)(0.5) + (0.5)(1) \\
& = 0.75 \\
& < 0.9.
\end{align*}

This shows that if the answer to the partition problem is no, then the answer to our instance of the decision form of the logit-based SOCPD problem is also no. \\ 

Since our instance of the logit-based SOCPD problem can be constructed in polynomial time from the instance of the partition problem, it follows that the logit-based SOCPD problem is NP-Hard even when the number of segments $K$ is equal to 2. \Halmos

\subsection{Proof of Theorem~\ref{theorem:APXHard}}
\label{proof:theorem_APXHard}

To prove this result, we will leverage a known inapproximability result for the maximum independent set (MAX-IS) problem. In the MAX-IS problem, we are given an undirected graph $G = (V,E)$, where $V$ is the set of vertices and $E$ is the set of edges. An independent set $U \subseteq V$ is a set of vertices such that for any pair of vertices $v, v' \in U$, $v \neq v'$, there does not exist an edge between them, that is, $(v,v') \notin E$. The goal in the MAX-IS problem is to find an independent set whose size is maximal. The MAX-IS problem is known to be NP-Hard to approximate to within a factor $O(n^{1-\epsilon})$ for any $\epsilon > 0$ \citep{hastad1996clique}. 

In what follows we will construct an approximation-preserving reduction that maps an instance of the MAX-IS problem to an instance of the unconstrained logit-based SOCPD problem. Given a graph $G = (V,E)$, let the number of attributes $n = |V|$, the number of segments $K = n$, and let $V = \{v_1,\dots, v_n\}$ be an enumeration of the vertices in $V$. Define the parameters $p_L$ and $p_U$ as
\begin{align*}
p_L & = \frac{1}{100n}, \\
p_U & = 1 - \frac{1}{100}.
\end{align*}
Observe that both $p_L$ and $p_U$ can be regarded as probabilities. Using $p_L, p_U$, define the parameters $q_L$ and $q_U$ as the logits corresponding to these probabilities:
\begin{align*}
q_L & = \log \frac{p_L}{1 - p_L}, \\
q_U & = \log \frac{p_U}{1 - p_U}.
\end{align*}
Let the utility parameters $\beta_{i,j}$ for $i \in [n]$, $j \in \{ 0 \} \cup [n]$ be defined as follows:
\begin{align*}
\beta_{i,j} & = \left\{ \begin{array}{ll} 
q_L & \text{if} \ j = 0, \\
q_U - q_L & \text{if}\ j = i, \\
q_L - q_U & \text{if}\ j < i \ \text{and}\ (v_i, v_j) \in E, \\
0 & \text{otherwise}, \end{array} \right.
\end{align*}
Note that by construction, the highest possible value that $\sigma(u_i(\ab))$ can attain is $p_U$, which occurs if $a_{i'} = 0$ for $i ' < i$ with $(v_{i'}, v_i) \in E$, and $a_i = 1$. Otherwise, for any other $\ab$, $u_i(\ab)$ satisfies $u_i(\ab) \leq q_L$, and so $\sigma(u_i(\ab)) \leq p_L = 1/(100n)$. 

Let the weight $\lambda_k$ of each segment $k$ be set to $1 / n$. Finally, let $F: \Acal \to [0,1]$ be defined as the share of choice objective function:
\begin{align*}
F(\ab) & \equiv \frac{1}{n} \sum_{i=1}^n \sigma( u_i(\ab)).
\end{align*}

To establish the result we need to verify two claims. 
\begin{enumerate}
	\item {\bf Claim 1}. For any independent set $U \subseteq V$, there exists a product $\ab$ such that $ F(\ab) \geq \frac{99}{100n} | U |$.
	\item {\bf Claim 2}. For any product $\ab$ with share-of-choice given by $F(\ab)$, there exists an independent set $U \subseteq V$ such that $|U| \geq \lfloor \frac{100n}{99} F(\ab) \rfloor$. 
\end{enumerate}

\noindent \emph{Proof of Claim 1}. Let $U \subseteq V$ be an independent set. Consider the product vector $\ab = (a_1, \dots, a_n)$ where $a_i = \Ibb\{ v_i \in U \}$. For each $i$ such that $v_i \in U$, we have:
\begin{align*}
u_i(\ab) & = \beta_{i,0} + \sum_{j=1}^n \beta_{i,j} a_j \\
& = q_L + \sum_{\substack{j = 1:\\ (v_i, v_j) \in E}}^{i-1} (q_L - q_U) a_j + (q_U - q_L) a_i \\
& = q_L + 0 + (q_U - q_L) \cdot 1 \\
& = q_U
\end{align*}
where the second equality follows by how the attribute utilities $\beta_{i,j}$ are defined; the third equality follows because $U$ is an independent set, so $a_j = 0$ for all attributes $j$ such that there exists an edge between $v_i$ and  $v_j$; the fourth follows by algebra. Thus, we have:
\begin{align*}
F(\ab) & = \frac{1}{n} \sum_{i=1}^n \frac{ \exp( u_i(\ab)) }{1 + \exp( u_i(\ab))} \\
& \geq \frac{1}{n} \cdot \sum_{i: v_i \in U} \frac{ \exp( u_i(\ab)) }{1 + \exp( u_i(\ab))} \\ 
& = \frac{1}{n} \cdot \sum_{i: v_i \in U} \frac{ \exp( q_U) }{1 + \exp(q_U)} \\
& = \frac{1}{n} \cdot |U| \cdot p_U \\
& = \frac{99}{100n} \cdot |U|.
\end{align*}

\noindent \emph{Proof of Claim 2}. Let $\ab$ be an attribute vector. Let us define the set $U$ as follows:
\begin{equation*}
U = \{ v_i \in V \mid \sigma(u_i(\ab)) \geq p_U \}.
\end{equation*}
In other words, we retrieve those vertices for which the corresponding segment's purchase probability is at least $p_U$. 

We argue that this set $U$ is an independent set. To see this, let us suppose for the sake of a contradiction that it is not. Then there exist two distinct vertices $v_i, v_{i'} \in U$ such that $(v_i, v_{i'}) \in E$. Without loss of generality, let us assume that $i < i'$. Observe that if we calculate the logit of segment $i'$, we have
\begin{align*}
u_{i'}(\ab) & = \beta_{i',0} + \sum_{j=1}^n \beta_{i', j} a_j \\
& = q_L + \sum_{ \substack{j=1:\\ (v_j, v_{i'}) \in E}}^{i'-1} (q_L - q_U) a_j + (q_U - q_L) a_{i'} \\
& \leq q_L + (q_L - q_U) a_i + (q_U - q_L) a_{i'} \\
& = q_L + (q_L - q_U) \cdot 1 + (q_U - q_L) \cdot 1 \\
& = q_L,
\end{align*}
where the inequality follows because $q_L - q_U < 0$. This implies that 
\begin{align*}
\frac{ \exp( u_{i'}(\ab)) }{1 + \exp(u_{i'}(\ab))} \leq p_L < p_U.
\end{align*}
This, however, leads to a contradiction, because $v_{i'}$ was assumed to be in $U$, which would imply that the corresponding purchase probability of segment $i'$ was higher than $p_U$. Therefore, it must be the case that $U$ is an independent set.

Now, we derive the desired bound on $|U|$. We have:
\begin{align*}
\left \lfloor \frac{100n}{99} F(\ab) \right \rfloor & = \left \lfloor \frac{100n}{99} \cdot \frac{1}{n} \sum_{i=1}^n \sigma( u_i(\ab)) \right \rfloor \\
& = \left \lfloor \frac{100}{99} \cdot \sum_{i: v_i \in U} \sigma( u_i(\ab)) + \frac{100}{99} \cdot \sum_{i: v_i \notin U} \sigma( u_i(\ab))  \right \rfloor \\ 
& \leq \left \lfloor \frac{100}{99} \cdot |U| \cdot \frac{99}{100} + \frac{100}{99} \cdot (n - |U|) \frac{1}{100n}  \right \rfloor \\ 
& \leq \left \lfloor |U| + \frac{100}{99} \cdot n \cdot \frac{1}{100n} \right \rfloor \\ 
& = \left \lfloor |U| + \frac{1}{99} \right  \rfloor \\ 
& = |U| 
\end{align*}
where the first step follows by definition of $F$; the second step follows by algebra; the third step follows because the floor function is monotonic, and because by definition of the utility parameters $\{\beta_{i,j}\}$, $\sigma(u_i(\ab)) \leq p_U = 99/100$ for all $i$, while for $i$ such that $v_i \notin U$, it is the case that $\sigma(u_i(\ab)) \leq p_L = 1/(100n)$; the fourth step again follows by monotonicity of the floor function and the fact that $(n - |U|) \leq n$; the fifth step follows by algebra; and the last step follows by the fact that $|U|$ is an integer while $1/99$ is strictly less than 1.

We now show how Claims 1 and 2 imply the required result. Suppose that we have access to a $\gamma$-approximation algorithm for the logit-based SOCPD problem. 

Consider any instance of the MAX-IS problem corresponding to a graph $(V,E)$. Consider the logit-based SOCPD instance that corresponds to this instance of MAX-IS, according to the reduction above. Let $\ab^*$ be the optimal attribute vector. We run the $\gamma$-approximation algorithm on this logit-based SOCPD instance, and obtain an attribute vector $\ab$ that satisfies
\begin{equation}
\gamma F(\ab^*) \leq F(\ab). \label{eq:MAXIS_gamma_approx_F}
\end{equation}
Let $U^*$ be any optimal independent set for MAX-IS. By Claim 1, there exists an attribute vector $\tilde{\ab}$ that satisfies
\begin{equation}
F(\tilde{\ab}) \geq \frac{99}{100n} | U^* |. \label{eq:MAXIS_tildeab_claim1}
\end{equation}
If we combine \eqref{eq:MAXIS_tildeab_claim1} and \eqref{eq:MAXIS_gamma_approx_F}, we get
\begin{equation}
\gamma \frac{99}{100n} |U^*| \leq \gamma F(\tilde{\ab}) \leq \gamma F(\ab^*) \leq F(\ab),
\end{equation}
where the second inequality follows because $\tilde{\ab}$ is a feasible solution of the logit-based SOCPD problem. This inequality simplifies to
\begin{equation}
\gamma \frac{99}{100n} |U^*| \leq F(\ab)
\end{equation}
or equivalently
\begin{equation}
\gamma |U^*| \leq \frac{100n}{99} F(\ab). \label{eq:MAXIS_gamma_Ustar_ab}
\end{equation}
We now use Claim 2 to assert that for $\ab$, there exists an independent set $\tilde{U}$ such that
\begin{equation}
|\tilde{U}| \geq \left \lfloor \frac{100n}{99} F(\ab) \right \rfloor. \label{eq:MAXIS_claim2}
\end{equation}
Taking floors of the left and right hand side of \eqref{eq:MAXIS_gamma_Ustar_ab}, and combining it with \eqref{eq:MAXIS_claim2}, we get
\begin{equation}
\lfloor \gamma | U^* | \rfloor \leq \left \lfloor \frac{100n}{99} F(\ab) \right \rfloor \leq |\tilde{U}|
\end{equation}
or simply $\lfloor \gamma | U^* | \rfloor \leq |\tilde{U}|$.

Therefore, if there exists a polynomial-time algorithm that approximates the SOCPD problem with approximation ratio $C /n^{1-\epsilon}$ for all instances, where $C$ is a constant, then there exists a polynomial-time algorithm $\Pcal$ that find an independent set $\tilde{U}$ with cardinality
\begin{align*}
|\tilde{U}| \geq \Bigg\lfloor \frac{C \cdot |U^*| }{n^{1-\epsilon}}  \Bigg\rfloor.
\end{align*}
Now, consider another polynomial-time algorithm $\Pcal'$ for the maximum independent set problem as follows. The algorithm $\Pcal'$ first runs $\Pcal$ to obtain an independent set $\tilde{U}$. If $|\tilde{U}| > 0$, $\Pcal'$ returns $\tilde{U}$. Otherwise, if $|\tilde{U}| = 0$, then $\Pcal'$ discards this $\tilde{U}$, picks an arbitrary item $i$ from the vertex set, and returns $\tilde{U} = \{i \}$, which is also an independent set and is of cardinality 1. (Note that this case can occur if the solution $\ab$ returned for the SOCPD problem is a zero vector $\zerob$.) Consequently, the solution $\tilde{U}$ returned by $\Pcal'$ is guaranteed to have size
\begin{align*}
|\tilde{U}| \geq \max \Bigg\lbrace \Bigg\lfloor \frac{C \cdot |U^*| }{n^{1-\epsilon}}  \Bigg\rfloor, 1 \Bigg\rbrace,
\end{align*} 
for all instances.

Next, we show that $\tilde{U}$ approximates the maximum independent set $U^*$ with a factor of $O(1/n^{1-\epsilon})$ for all problem instances. For problem instances where $C |U^*| \geq n^{1-\epsilon}$, or equivalently $C |U^*| / n^{1-\epsilon} \geq 1$, we have
\begin{align*}
\tilde{U} \geq  \max \Bigg\lbrace \Bigg\lfloor \frac{C \cdot |U^*| }{n^{1-\epsilon}}  \Bigg\rfloor, 1 \Bigg\rbrace  \geq  \Bigg \lfloor \frac{C \cdot |U^*| }{n^{1-\epsilon}} \Bigg \rfloor \geq \frac{ (C/2)\cdot |U^*| }{n^{1-\epsilon}} ,
\end{align*}
where the third inequality follows by the fact that $\lfloor x \rfloor \geq x/2$ whenever $x \geq 1$. This implies that $\Pcal'$ can approximate $U^*$ with a factor of $(C/2) / n^{1-\epsilon}$. For problem instances where $C |U^*| < n^{1-\epsilon}$, or equivalently $C |U^*| / n^{1-\epsilon} < 1$, we have
\begin{align*}
\tilde{U} \geq  \max \Bigg\lbrace \Bigg \lfloor \frac{C \cdot |U^*| }{n^{1-\epsilon}}  \Bigg\rfloor, 1 \Bigg\rbrace = 1 > \frac{C |U^*|}{n^{1-\epsilon}} \geq \frac{(C/2) |U^*|}{n^{1-\epsilon}} ,
\end{align*}
also implying that $\Pcal'$ approximates $U^*$ within a factor of $(C/2) / n^{1-\epsilon}$. Therefore, for all instances, $\Pcal'$ is a polynomial-time $O(1/n^{1-\epsilon})$-approximation algorithm to the maximum independent set problem; thus, by the result of \cite{hastad1996clique}, it follows that the logit-based SOCPD problem is also NP-Hard to approximate to a factor better than $O(1/n^{1-\epsilon})$ for any $\epsilon > 0$. \Halmos

\subsection{Proof of Theorem~\ref{theorem:APXHard_MAX3SAT}}
\label{proof:theorem_APXHard_MAX3SAT}

To prove this result, we will develop an approximation preserving reduction between the MAX-3SAT problem and the logit-based SOCPD problem. The MAX-3SAT problem is the problem of setting a collection of binary variables so as to maximize the number of clauses, which are disjunctions of up to three literals, that are satisfied. More precisely, we have $n$ binary variables, $x_1, \dots, x_n$, and a Boolean formula $c_1 \wedge c_2 \wedge \dots  \wedge c_K$, where the symbol $\wedge$ denotes the ``and'' operator. Each $c_k$ is a disjunction involving three literals where a literal is one of the binary variables or the negation of one of the binary variables. For example, a clause could be $x_1 \vee x_4 \vee \neg x_9$ where $\vee$ denotes the ``or'' operator and $\neg$ denotes negation; in this example, the clause evaluates to 1 if $x_1 = 1$, or $x_4 = 1$, or $x_9 = 0$, and evaluates to zero if $x_1 = 0$, $x_4 = 0$ and $x_9 = 1$. The MAX-3SAT problem is to determine how $x_1, \dots, x_n$ should be set so that the number of the clauses $c_1,\dots, c_K$ that are true is maximized. 

Given an instance of the MAX-3SAT problem, we show how the instance can be transformed into an instance of the logit-based SOCPD problem. %

In the instance of the logit-based SOCPD problem that we will construct, we let the number of attributes $n$ be equal to the number of binary variables in the MAX-3SAT instance, and we let the set of permissible attribute vectors $\Acal$ simply be equal to $\{0,1\}^n$. Each attribute of our product will correspond to one of the binary variables. We let each customer type $k$ correspond to one of the $K$ clauses, and we set $\lambda_k = 1 / K$. To aid in defining the partworths of each customer type, we will define the parameters $p_L$ and $p_U$ as 
\begin{align}
p_L & = \frac{1}{MK}, \\
p_U & =1 - \frac{1}{MK},
\end{align}
and we define the utilities $Q_L$, $Q_U$ as the inverse of the logistic response function of each of these:
\begin{align}
Q_L & = \log \left[ \frac{ 1 / (MK) }{1 - 1/(MK) } \right], \\
Q_U & = \log \left[ \frac{ 1 - 1 / (MK) }{1/(MK) } \right],
\end{align}
where $M \geq 3$ is an arbitrary integer. 

Now, for each customer type $k$, let $J_k \in \{0,1,2,3\}$ denote the number of negative literals in the corresponding clause $k$ of the MAX-3SAT instance (i.e., how many literals of the form $\neg x_i$ appear in $c_k$). We define the partworths $\beta_{k,1}, \dots, \beta_{k,n}$ of customer type $k$ as follows:
\begin{align}
\beta_{k,i} & = \begin{cases}
0 & \text{if variable $x_i$ does not appear in any literal of clause $k$}, \\
Q_U - Q_L & \text{if the literal $x_i$ appears in clause $k$}, \\
Q_L - Q_U & \text{if the literal $\neg x_i$ appears in clause $k$}, \\
\end{cases}
\end{align}
for each $i \in \{1,\dots, n\}$, and we define the constant part of the utility $\beta_{k,0}$ as 
\begin{align}
\beta_{k,0} & = Q_L + J_k \cdot (Q_U - Q_L). 
\end{align}
The rationale for this choice is that the utility of an attribute vector $\ab$ will be equal to $Q_L$ if the attributes are set in a way such that none of the literals of clause $k$ are satisfied, and will be equal to $Q_U$ or higher if the attributes are set so that the clause is satisfied (i.e., at least one literal is true). For example, if clause $k$ is $c_k = x_1 \vee x_4 \vee \neg x_9$, then the corresponding utility function of customer type $k$ has the form:
\begin{align*}
u_k(\ab) & = Q_L + 1 \cdot (Q_U - Q_L) + (Q_U - Q_L) a_1 + (Q_U - Q_L) a_4 + (Q_L - Q_U) a_9 \\
& = Q_U + (Q_U - Q_L) a_1 + (Q_U - Q_L) a_4 + (Q_L - Q_U) a_9.
\end{align*}
If $a_1 = 1$, $a_4 = 0$ and $a_9 =1$, the clause evaluates to 1 ($= 1 \vee 0 \vee \neg 1$); the utility is 
\begin{align*}
u_k(\ab) & = Q_U + (Q_U - Q_L) \cdot 1 + (Q_U - Q_L) \cdot 0 + (Q_L - Q_U) \cdot 1 \\
& = Q_U. 
\end{align*}
If $a_1 = 0$, $a_4 = 0$, $a_9 = 1$, the clause evaluates to 0 ($ = 0 \vee 0 \vee \neg 1$), and the utility is 
\begin{align*}
u_k(\ab) & = Q_U + (Q_U - Q_L) \cdot 0 + (Q_U - Q_L) \cdot 0 + (Q_L - Q_U) \cdot 1 \\
& = Q_L,
\end{align*}
as expected.

Lastly, before we verify that this reduction is valid, it is helpful to introduce some additional notation to model the MAX-3SAT. Given a binary vector $\xb \in \{0,1\}^n$, we let $g_k(\xb) = 1$ if clause $k$ is satisfied and 0 if clause $k$ is not satisfied. The MAX-3SAT problem can then be written simply as 
\begin{equation*}
\max_{\xb \in \{0,1\}^n} \sum_{k=1}^K g_k(\xb).
\end{equation*}

Let $G^*$ denote the optimal objective value of the MAX-3SAT problem, and let $\xb^*$ denote an optimal solution of this problem. Similarly, let $F(\ab)$ denote the share-of-choice of a product vector $\ab \in \{0,1\}^n$, and let $F^*$ denote the optimal objective value of the logit-based SOCPD problem.

We now establish two claims.

\begin{enumerate}
\item {\bf Claim 1}. For any MAX-3SAT solution $\xb$, the SOCPD solution $\ab = \xb$ is such that $F(\ab) \leq \frac{1}{K} G(\xb) + p_L$. 
\item {\bf Claim 2}. For any product $\ab$, the MAX-3SAT solution $\xb = \ab$ is such that $G(\xb) \leq \frac{K}{p_U} F(\ab)$. 
\end{enumerate}

We first verify Claim 1. First, observe that
\begin{align*}
K \cdot F(\ab) & = K \cdot \frac{1}{K} \sum_{k=1}^K \sigma( u_k(\ab) ) \\
& = \sum_{k=1}^K \sigma( u_k(\ab) ) \\
& \leq \sum_{k=1}^K (g_k(\xb) + p_L) \\
& = G(\xb) + K p_L,
\end{align*}
where the inequality follows because for each $k$, $\sigma( u_k(\ab) ) = p_L$ if $g_k(\xb) = 0$, and otherwise, if $g_k(\xb) = 1$, then we simply use the fact that $\sigma(u_k(\ab)) < 1$. rearranging this inequality gives
\begin{equation}
F(\ab) \leq \frac{1}{K} G(\xb) + p_L,
\end{equation}
and we thus establish Claim 1. 

We now verify Claim 2. Observe that 
\begin{align*}
p_U G(\xb) & = \sum_{k=1}^K p_U g_k(\xb) \\
& \leq \sum_{k=1}^K \sigma( u_k(\ab) ) \\
& = K \cdot F(\ab),
\end{align*}
where the inequality follows because when $g_k(\xb) = 1$, $\sigma(u_k(\ab)) \geq p_U$, and otherwise, when $g_k(\xb) = 0$, we can use the fact that $\sigma(u_k(\ab)) > 0$. 

With these two properties in hand, let us now proceed as follows. Suppose that we obtain a $\gamma$ approximation of the logit-based SOCPD problem, which means that we identify an $\ab \in \{0,1\}^n$ such that
\begin{equation}
F(\ab) \geq \gamma F^*.
\end{equation}
Let $\xb = \ab$ be a candidate solution of MAX-3SAT. We will show now that $\xb$ satisfies the approximation guarantee
\begin{equation}
G(\xb) \geq \left( \gamma - \frac{(\gamma+1)}{M} \right) G^*.
\end{equation}

With regard to $F^*$, let $\ab'$ be any arbitrary vector in $\{0,1\}^n$, and let $\xb' = \ab'$. Then by Claim 2 we have
\begin{align}
F^* & \geq F(\ab') \\
& \geq (p_U / K) G(\xb'),
\end{align}
which implies that for all $\xb' \in \{0,1\}^n$, $KF^* / p_U$ is an upper bound on $G(\xb')$; hence, we have $KF^* / p_U \geq G^*$, or equivalently,
\begin{equation}
F^* \geq \frac{p_U}{K} G^*. \label{eq:Fstar_lower_bound}
\end{equation}
Combining \eqref{eq:Fstar_lower_bound} with Claim 1, we get
\begin{align}
\frac{1}{K} G(\xb) + p_L \geq \gamma \cdot \frac{p_U}{K} G^*.
\end{align}
Multiplying left and right by $K$ and re-arranging terms, we get
\begin{align}
G(\xb) \geq \gamma \cdot p_U G^* - K p_L 
\end{align}
To simplify the right hand side, we use the definitions of $p_U$ and $p_L$ to get
\begin{align}
G(\xb) & \geq \gamma \cdot p_U G^* - K p_L  \\
& = \gamma \cdot (1 - \frac{1}{MK}) G^* - K \cdot \frac{1}{MK} \\
& = \gamma \cdot (1 - \frac{1}{MK}) G^* - \frac{1}{M} \\
& \geq \gamma \cdot (1 - \frac{1}{M}) G^* - \frac{1}{M} G^* \\
& = \left(\gamma - \frac{(\gamma + 1)}{M}  \right) G^*,
\end{align}
where in the inequality we use the fact that $K \geq 1$ and that $G^* \geq 1$ (it is always possible to satisfy at least one clause). 

Since our SOCPD instance is constructed in polynomial time from the MAX-3SAT instance, since the candidate solution $\xb$ for MAX-3SAT is obtained trivially from the approximate solution $\ab$ (and therefore in polynomial time), and since MAX-3SAT is APX-complete \citep{papadimitriou1991optimization}, it follows that the logit-based SOCPD problem is APX-Hard.

For the last part of the theorem, observe that by \cite{haastad2001some}, the MAX-E3SAT problem is NP-hard to approximate to a factor of $7/8 + \epsilon$ for any positive $\epsilon$. (Recall that the MAX-E3SAT problem is a specific case of the MAX-3SAT problem where each clause consists of exactly three literals; in contrast to MAX-E3SAT, a MAX-3SAT instance could have some clauses consisting of fewer than three literals.) In the case of MAX-E3SAT, the exact same reduction described above goes through. For the logit-based SOCPD problem, any $\gamma > 7/8$ and $M$ sufficiently large would imply an approximation algorithm with factor greater than $7/8$ for MAX-E3SAT. Hence, when each customer type has exactly three non-zero partworths, it also follows that the logit-based SOCPD problem is NP-hard to approximate to a factor greater than $7/8 + \epsilon$ for any $\epsilon > 0$. \Halmos

\subsection{Proof of Proposition~\ref{proposition:P_leq_RA}}
\label{proof:proposition_P_leq_RA}

Let $(\bar{\ab}, \bar{\wb}, \bar{\xb}, \bar{\yb})$ be an optimal solution of the continuous relaxation of formulation~\modelP, and let $\bar{\ub} \in \Rbb^K$ be the vector of utilities corresponding to $\ab$.  To establish the proposition, we will first prove that 
\begin{align*}
\bar{x}_{k,1} \geq \frac{1}{1+e^{-\bar{u}_k}}, \\
\bar{x}_{k,0} \leq \frac{1}{1+e^{\bar{u}_k}},
\end{align*}
for each $k \in [K]$. To see why this is the case, consider the following optimization problem, which involves finding the maximum value of $x_{k,1}$ given the fixed value of $\bar{\ab}$, subject to the constraints of \modelP, and with the additional restriction that $y_{k,i}$ is exactly equal to the product of $\bar{a}_i$ and $x_{k,1}$:
\begin{subequations}
\begin{alignat}{2}
& \underset{ w_k, \xb_k, \yb_k}{ \text{maximize}} & \quad & x_{k,1} \\
& \text{subject to} & & y_{k,i} = \bar{a}_i \cdot x_{k,1}, \quad \forall \ i \in [n], \label{prob:P_nonlinear_bilinear}\\
& & & w_k = \beta_{k,0} x_{k,1} + \sum_{i=1}^n \beta_{k,i} y_{k,i}, \label{prob:P_nonlinear_w}\\
& & & x_{k,1} + x_{k,0} = 1, \\
& & & x_{k,1} + x_{k,1} e^{ \frac{- w_k}{x_{k,1}}} \leq 1, \label{prob:P_nonlinear_perspective}\\
& & & y_{k,i} \leq \bar{a}_i, \quad \forall \ i \in [n], \\
& & & y_{k,i} \leq x_{k,1}, \quad \forall \ i \in [n], \\
& & & y_{k,i} \geq \bar{a}_i - 1 + x_{k,1}, \quad \forall \ i \in [n], \\
& & & y_{k,i} \geq 0.
\end{alignat}
\label{prob:P_nonlinear}
\end{subequations}
Observe that the optimal solution of this problem is 
\begin{align}
x^*_{k,1} & = \frac{1}{1+e^{-\bar{u}_k}}, \label{eq:P_opt_soln_1}\\
x^*_{k,0} & = \frac{1}{1+e^{\bar{u}_k}}, \\
y^*_{k,i} & = \bar{a}_i \cdot x_{k,1}, \quad \forall \ i \in [n], \\
w^*_k & = x_{k,1} \cdot \bar{u}_k. \label{eq:P_opt_soln_5}
\end{align}
To see this, observe that the above solution is feasible for the problem~\eqref{prob:P_nonlinear}. In addition, observe that constraints~\eqref{prob:P_nonlinear_bilinear} and \eqref{prob:P_nonlinear_w} imply that $w_k$ must be equal to  $u^*_k \cdot x_{k,1}$. As a result, \eqref{prob:P_nonlinear_perspective} implies that $x_{k,1}$, which is the objective, is upper bounded in the following way:
\begin{align}
x_{k,1} + x_{k,1} e^{ \frac{-w_k}{x_{k,1}}} & \leq 1 \\
\Rightarrow x_{k,1} (1 + e^{\frac{-x_{k,1} \cdot \bar{u}_k}{x_{k,1}}}) & \leq 1 \\
\Rightarrow x_{k,1} (1 + e^{-\bar{u}_k}) & \leq 1 \\
\Rightarrow x_{k,1} & \leq \frac{1}{1 + e^{-\bar{u}_k}} 
\end{align}
Since the proposed solution in \eqref{eq:P_opt_soln_1} - \eqref{eq:P_opt_soln_5} attains this upper bound, it must also be optimal.

Next, consider what happens if we relax the constraint~\eqref{prob:P_nonlinear_bilinear}. In doing so we obtain the following program:
\begin{subequations}
\begin{alignat}{2}
& \underset{ w_k, \xb_k, \yb_k}{ \text{maximize}} & \quad & x_{k,1} \\
& \text{subject to} & & w_k = \beta_{k,0} x_{k,1} + \sum_{i=1}^n \beta_{k,i} y_{k,i}, \label{prob:P_relaxed_w}\\
& & & x_{k,1} + x_{k,0} = 1, \\
& & & x_{k,1} + x_{k,1} e^{ \frac{- w_k}{x_{k,1}}} \leq 1, \label{prob:P_relaxed_perspective}\\
& & & y_{k,i} \leq \bar{a}_i, \quad \forall \ i \in [n], \\
& & & y_{k,i} \leq x_{k,1}, \quad \forall \ i \in [n], \\
& & & y_{k,i} \geq \bar{a}_i - 1 + x_{k,1}, \quad \forall \ i \in [n], \\
& & & y_{k,i} \geq 0.
\end{alignat}
\label{prob:P_relaxed}
\end{subequations}
Since this problem is a relaxation, the optimal solution $(w'_k, \xb'_k, \yb'_k)$ must do at least as well as $(w^*_k, \xb^*_k, \yb^*_k)$. This means that $x'_{k,1} \geq x^*_{k,1} = 1/(1+e^{-\bar{u}_k})$, and similarly that $x'_{k,0} \leq 1/(1+e^{\bar{u}_k})$. 

Coming back to the solution $(\bar{\ab}, \bar{\wb}, \bar{\xb}, \bar{\yb})$ of the relaxation of formulation~\modelP, observe that utilizing the argument above and the fact that the objective is a nonnegative weighted combination of the $x_{k,1}$ variables, we will have that $\bar{x}_{k,1} \geq 1 / (1 + e^{-\bar{u}_k})$ and that $\bar{x}_{k,0} \leq 1 / (1 + e^{ \bar{u}_k})$ for each $k$ for which $\lambda_k > 0$. Without loss of generality, we can also assume that these inequalities hold for all $k \in [K]$, since there is no contribution to the objective function of formulation~\modelP from the term $\lambda_k x_{k,1}$ for any $k$ with $\lambda_k = 0$. 

With this property of $(\bar{\ab}, \bar{\wb}, \bar{\xb}, \bar{\yb})$ established, we now claim that $(\bar{\ab}, \bar{\ub}, \bar{\wb}, \bar{\xb}, \bar{\yb})$ is a feasible solution of the relaxation of formulation~\modelRA. Note that this amounts to verifying that $(\bar{\ab}, \bar{\ub}, \bar{\wb}, \bar{\xb}, \bar{\yb})$ satisfies the representative agent constraint
\begin{equation}
w_k - x_{k,1} \log x_{k,1} - x_{k,0} \log x_{k,0} \geq \log(1 + e^{u_k})
\end{equation}
for every $k$, since the other constraints in \modelRA are already present in \modelP. To see why this constraint is satisfied by our solution $(\bar{\ab}, \bar{\ub}, \bar{\wb}, \bar{\xb}, \bar{\yb})$, observe that:
\begin{align}
\bar{x}_{k,1} + \bar{x}_{k,1} e^{ \frac{- \bar{w}_{k}}{\bar{x}_{k,1}}} & \leq 1 \\
\Rightarrow \quad \bar{x}_{k,1} e^{ \frac{- \bar{w}_{k}}{\bar{x}_{k,1}}} & \leq \bar{x}_{k,0} \\
\Rightarrow \quad \log \bar{x}_{k,1} + \frac{- \bar{w}_{k}}{\bar{x}_{k,1}} & \leq \log \bar{x}_{k,0} \\
\Rightarrow \quad \bar{x}_{k,1} \log \bar{x}_{k,1} - \bar{w}_{k} & \leq (1 - \bar{x}_{k,0}) \log \bar{x}_{k,0} \\
\Rightarrow \quad \bar{w}_{k} - \bar{x}_{k,0} \log \bar{x}_{k,0} - \bar{x}_{k,1} \log \bar{x}_{k,1} &  \geq - \log \bar{x}_{k,0} \label{eq:almost_RA_constraint}
\end{align}
Now, recall that $\bar{x}_{k,0} \leq 1 / (1 + e^{ \bar{u}_k})$, or equivalently (after taking logs and multiplying by -1):
\begin{align}
- \log \bar{x}_{k,0} & \geq \log(1 + e^ {\bar{u}_k}).  \label{eq:minuslogx0_geq_softplus}
\end{align}
Inequality~\eqref{eq:almost_RA_constraint} and \eqref{eq:minuslogx0_geq_softplus} together imply that
\begin{equation}
\bar{w}_k - \bar{x}_{k,1} \log \bar{x}_{k,1} - \bar{x}_{k,0} \log \bar{x}_{k,0} \geq \log(1 + e^{\bar{u}_k}),
\end{equation}
which is exactly the representative agent constraint of formulation~\modelRA. As a result, we have established that $(\bar{\ab}, \bar{\ub}, \bar{\wb}, \bar{\xb}, \bar{\yb})$ is a feasible solution of the relaxation of formulation~\modelRA. Since the two formulations share the same objective functions, it thus follows that $Z^*_{\modelP} \leq Z^*_{\modelRA}$, as required. \Halmos

\clearpage
\pagebreak

\section{Extra modeling details}

In this section, we provide some additional discussion of the modeling capability of our mixed-integer convex programming models discussed in Section~\ref{sec:micp}. Section~\ref{subsec:extra_Acal} provides some examples of what can be modeled using the linear constraint $\Cb \ab \leq \db$ that defines $\Acal$, while Section~\ref{subsec:extra_profit} discusses how the three formulations (\modelRA, \modelP and \modelPRPT) can be modified for the purpose of expected profit maximization.

\subsection{Set of feasible product designs}
\label{subsec:extra_Acal}

The constraint $\Cb \ab \leq \db$ which defines the set $\Acal$ can be used to encode a variety of requirements on the attribute vectors $\ab$ as linear constraints. For example, if the product has two attributes, where the first attribute has three levels and the second attribute has four levels, then one can model the product through the vector $\ab = (a_1, a_2, a_3, a_4, a_5)$, where $a_1$ and $a_2$ are dummy variables to represent two out of the three levels of the first attribute and $a_3, a_4, a_5$ are dummy variables to represent three out of the four levels of the second attribute. One would then need to enforce the constraints
\begin{align}
& a_1 + a_2 \leq 1, \label{eq:attribute1_constraint}\\
& a_3 + a_4 + a_5 \leq 1 \label{eq:attribute2_constraint}
\end{align}
to ensure that at most one out of the variables $a_1, a_2$ is set to 1 and at most one variable out of $a_3, a_4, a_5$ is set to 1. This can be achieved by specifying $\Cb$ and $\db$ as
\begin{align*}
\Cb = \left[ \begin{array}{ccccc} 1 & 1 & 0 & 0 & 0 \\ 0 & 0 & 1 & 1 & 1 \end{array} \right], \ \db = \left[ \begin{array}{c} 1 \\ 1 \end{array} \right] \label{eq:multi_level_attribute_Cb_db}
\end{align*}   

Beside the ability to represent multi-level attributes, one can use the constraint $\Cb \ab \leq \db$ to represent design requirements such as weight and cost; for example, one may be interested in imposing the constraint 
\begin{equation*}
b_0 + \sum_{i=1}^n b_i a_i \leq B,
\end{equation*}
where $b_0$ is the base weight of the product, $b_i$ is the incremental weight added to the product from attribute $i$ and $B$ is a limit on the overall weight of the product. This constraint can be modeled by specifying $\Cb$ and $\db$ as 
\begin{align*}
\Cb = \left[ \begin{array}{cccc} b_1 & b_2 & \cdots & b_n \end{array} \right], \ \db = \left[ B - b_0 \right]. 
\end{align*}

\subsection{Extension to expected profit maximization}
\label{subsec:extra_profit}

While all three of our formulations \modelRA, \modelP and \modelPRPT corresponds to the share-of-choice objective, it turns out that it is straightforward to generalize these models so as to optimize a profit-based objective. In particular, suppose that the marginal profit of a design $\ab$ is given by a function $R(\ab)$ defined as 
\begin{equation}
R(\ab) = r_0 + \sum_{i=1}^n r_i a_i. \label{eq:profit_linear}
\end{equation}
In other words, the profit $R(\ab)$ is a linear function of the binary attributes. One can model various types of profit structures with this assumption. For example, if all of the attributes correspond to non-price features that affect the cost of the product, then one can set $r_0$ to be the price of the product (a positive quantity), and each $r_i$ to be the marginal incremental cost of attribute $i$ (a negative quantity). 

With this assumption, the logit-based expected profit product design problem can be written as 
\begin{equation}
\underset{ \ab \in \Acal}{\text{maximize}} \quad R(\ab) \cdot \left[ \sum_{k=1}^K \lambda_k \cdot \frac{ \exp( u_k(\ab)) }{1 + \exp(u_k(\ab))} \right]. \label{prob:EPPD_abstract}
\end{equation}
In terms of the $x_{k,1}$ decision variables that appear in formulations~\modelRA, \modelP and \modelPRPT, the objective function can be re-written as
\begin{align*}
R(\ab) \cdot \left[ \sum_{k=1}^K \lambda_k \cdot x_{k,1} \right] & = \left( r_0 + \sum_{i=1}^n r_i a_i \right) \cdot \left[ \sum_{k=1}^K \lambda_k \cdot x_{k,1} \right] \\
& = \sum_{k=1}^K \lambda_k \cdot \left[ r_0 x_{k,1} + \sum_{i=1}^n r_i \cdot a_i x_{k,1} \right].
\end{align*}
Notice that this last expression includes terms of the form $a_i x_{k,1}$, which we can already represent through the variables $y_{k,i}$ that appear in all three formulations. We can therefore re-write the objective function of problem~\eqref{prob:EPPD_abstract} as
\begin{equation*}
\sum_{k=1}^K \lambda_k \cdot \left[ r_0 x_{k,1} + \sum_{i=1}^n r_i \cdot y_{k,i} \right] 
\end{equation*}
Thus, the expected profit product design problem can be handled by modifying the objective function of formulation \modelRA / \modelP / \modelPRPT.

\subsection{Extension to product line design}
\label{subsec:extra_PLD}

In this section, we discuss how two of our formulations can be adapted for the setting where the firm must design more than one product, i.e., a product line. 

In this setting, we let $J$ denote the number of products to be designed (the width of the product line), and let $j$ denote the index of a product in the product line, which ranges from 1 to $J$. We use the binary decision variable $a_{j,i}$ to denote whether product $j$ in the product line has attribute $i$ or not. We let $\ab_j$ denote the vector of attributes for the $j$th product. As in the single product design case, we assume $\Acal$ is the set of feasible attribute vectors, and we again assume that $\Acal = \{ \ab \in \{0,1\}^n \mid \Cb \ab \leq \db \}$, for some choice of a matrix $\Cb$ and vector $\db$. 

We assume that each customer chooses between one of the $J$ products and the no-purchase option. For each customer type $k$, we assume that the utility of the $j$th product is
\begin{equation}
u_k(\ab) = \beta_{k,0} + \sum_{i=1}^n \beta_{k,i} a_{j,i}.
\end{equation}
The choice probability of the $j$th product for customer type $k$ is then given by 
\begin{equation}
\frac{ \exp( u_k(\ab_j)) }{ \sum_{j'=1}^J \exp( u_k(\ab_{j'})) + 1},
\end{equation}
where we assume without loss of generality that the no-purchase option utility is normalized to zero. As before, we assume that $\lambda_k$ is the probability of a customer belonging to type $k$. 

As in Section~\ref{subsec:extra_profit}, we assume that the marginal profit of a product is a linear function of its attributes, given by equation~\eqref{eq:profit_linear}. The product line design problem can then be written abstractly as follow:
\begin{equation}
\max_{(\ab_1,\dots, \ab_J) \in \Acal^J} \sum_{k=1}^K \lambda_k \cdot \frac{ \sum_{j=1}^J R(\ab_j) e^{u_k(\ab_j)}}{ \sum_{j=1}^J e^{u_k(\ab_j)} + 1}.
\end{equation}

We will now show how one can formulate this problem using both the representative agent approach, leading to a model akin to formulation~\modelRA, and using the perspectification approach, leading to a model analogous to formulation~\modelP. 

Let us define the decision variable $x_{k,j}$ for $j = 1,\dots, J$, which is the choice probability of product $j$ for customer type $k$; additionally, let $x_{k,0}$ be the choice probability of the no-purchase option. Let the decision variable $y_{k,j,i}$ represent the linearization of $x_{k,j} \cdot a_{j,i}$, and let $w_{k,j}$ represent the linearization of $x_{k,j} \cdot u_{k,j}$. Then the formulation can be written as follows.
\begin{subequations}
\begin{alignat}{2}
& \underset{\ab, \ub, \wb, \xb, \yb  }{\text{maximize}} & \quad & \sum_{k=1}^K \lambda_k \cdot \sum_{j=1}^J \left( r_0 x_{k,j} + \sum_{i=1}^n r_i y_{k,j,i} \right) \\
& \text{subject to} & & \sum_{j=0}^J x_{k,j} = 1, \quad \forall k \in [K], \\
& & & \sum_{j=1}^J w_{k,j} - \sum_{j=0}^J x_{k,j} \log x_{k,j} \geq \log(1 + \sum_{j=1}^J e^{u_{k,j}}), \quad \forall k \in [K], \label{prob:PLD_RA_bilevel_constraint} \\
& & & u_{k,j} = \beta_{k,0} + \sum_{i=1}^n \beta_{k,i} a_{j,i}, \quad \forall k \in [K], j \in [J], \\
& & & w_{k,j} = \beta_{k,0} x_{k,j} + \sum_{i=1}^n \beta_{k,i} y_{k,j,i}, \quad \forall k \in [K], j \in [J], \label{prob:PLD_RA_w_def}\\ 
& & & y_{k,j,i} \leq x_{k,j}, \quad \forall k \in [K], j \in [J], i \in [n], \label{prob:PLD_RA_y_1} \\
& & & y_{k,j,i} \leq a_{j,i}, \quad \forall k \in [K], j \in [J], i \in [n], \\
& & & y_{k,j,i} \geq a_{j,i} - 1 + x_{k,j,i}, \quad \forall k \in [K], j \in [J], i \in [n], \\
& & & y_{k,j,i} \geq 0, \quad \forall k \in [K], j \in [J], i \in [n], \label{prob:PLD_RA_y_4} \\
& & & \Cb \ab_j \leq \db, \quad \forall j \in [J], \\
& & & a_{j,i} \in \{0,1\}, \quad \forall j \in [J], i \in [n], \\
& & & x_{k,j} \geq 0, \quad \forall k \in [K], j \in [J].
\end{alignat}
\label{prob:PLD_RA}%
\end{subequations}
This formulation behaves similarly to formulation~\modelRA, and is the generalization of that formulation to the product line setting. In particular, when $\ab_1,\dots, \ab_J$ are binary, constraints~\eqref{prob:PLD_RA_y_1} - \eqref{prob:PLD_RA_y_4} ensure that $y_{k,j,i} = x_{k,j} \cdot a_{j,i}$, and \eqref{prob:PLD_RA_w_def} ensures that $w_{k,j} = u_{k,j} \cdot x_{k,j}$. Consequently, constraint~\eqref{prob:PLD_RA_bilevel_constraint} ensures that the objective value of the representative agent problem with $J+1$ alternatives, with utilities $u_{k,1},\dots, u_{k,J}, 0$, is at least as high as the optimal objective value, which is $\log(1 + \sum_{j=1}^J e^{u_{k,j}})$. This can only occur if each $x_{k,j} = e^{u_{k,j}} / (\sum_{j'=1}^J e^{u_{k,j'}} + 1)$, which are exactly the choice probabilities from the multinomial logit choice model.

Now, let us see how to generalize formulation~\modelP to the product line setting. We retain the same definitions of the $x_{k,j}$ and $a_{j,i}$ decision variables from the previous formulation. We now define $v_{k,j,j'}$ as the linearization of $x_{k,j} u_{k,j'}$, for each $j \in \{0,\dots,J\}$, $j' \in [J]$ and each $k \in [K]$. We also define $z_{k,j,j',i}$ as the linearization of $x_{k,j} \cdot a_{j',i}$, for each $k \in [K]$, $j \in \{0,\dots,J\}$, $j' \in [J]$ and $i \in [n]$.  

To see how perspectification applies in the product line setting, let us start from the ideal constraint we would want to impose, which is the following non-convex inequality:
\begin{equation}
x_{k,j} \leq \frac{ e^{u_{k,j}}}{\sum_{j'=1}^J e^{u_{k,j'}} + 1}.
\end{equation}
By dividing the numerator and the denominator on the right by $e^{u_{k,j}}$, we get
\begin{equation}
x_{k,j} \leq \frac{ 1 }{\sum_{j'=1}^J e^{u_{k,j'} - u_{k,j}} + e^{-u_{k,j}}}.
\end{equation}
We now re-arrange the constraint to get
\begin{equation}
\sum_{j'=1}^J x_{k,j} e^{u_{k,j'} - u_{k,j}} + x_{k,j} e^{-u_{k,j}} \leq 1.
\end{equation}
This constraint is non-convex, but by applying perspectification and using the definition of the variables $v_{k,j,j'}$, we obtain 
\begin{equation}
\sum_{j'=1}^J x_{k,j} e^{ \frac{v_{k,j,j'} - v_{k,j,j}}{ x_{k,j}}} + x_{k,j} e^{ - v_{k,j,j} / x_{k,j} } \leq 1,
\end{equation}
which is convex. 

Using this, together with appropriate constraints for linearization, we arrive at the following formulation, which generalizes \modelP to the product line setting.

\begin{subequations}
\begin{alignat}{2}
& \underset{\ab, \vb, \xb, \zb}{ \text{maximize} } & & \sum_{k=1}^K \lambda_k \cdot \sum_{j=1}^J \left( r_0 x_{k,j} + \sum_{i=1}^n r_i z_{k,j,j,i} \right) \\
& \text{subject to} & \quad & \sum_{j'=1}^J x_{k,j} e^{ \frac{v_{k,j,j'} - v_{k,j,j}}{ x_{k,j}}} + x_{k,j} e^{ - v_{k,j,j} / x_{k,j} } \leq 1, \quad \forall k \in [K], j \in [J], \label{prob:PLD_P_perspective_j}\\
& & & \sum_{j'=1}^J x_{k,0} e^{ \frac{v_{k,0,j'}}{ x_{k,0}}} + x_{k,0} \leq 1, \quad \forall k \in [K], \label{prob:PLD_P_perspective_0} \\
& & & \sum_{j=0}^J z_{k,j,j',i} = a_{j',i}, \quad \forall k \in [K], j' \in [J], i \in [n], \label{prob:PLD_P_z_1} \\
& & & z_{k,j,j',i} \leq x_{k,j}, \quad \forall k \in [K], j \in \{0,\dots,J\}, j' \in [J], i \in [n] \label{prob:PLD_P_z_2} \\
& & & z_{k,j,j',i} \geq 0, \quad \forall k \in [K], j \in \{0,\dots,J\}, j' \in [J], i \in [n] \label{prob:PLD_P_z_3} \\
& & & v_{k,j,j'} = \beta_{k,0} x_{k,j} + \sum_{i=1}^n \beta_{k,i} z_{k,j,j',i}, \quad \forall k \in [K], j \in \{0,\dots,J\}, j' \in [J], \label{prob:PLD_P_v} \\
& & & \sum_{j=0}^J x_{k,j} = 1, \quad \forall k \in [K], \\
& & & \Cb \ab_j \leq \db, \quad \forall j \in [J], \\
& & & a_{j,i} \in \{0,1\}, \quad \forall j \in [J], i \in [n], \\
& & & x_{k,j} \geq 0, \quad \forall k \in [K], j \in \{0,\dots,J\}. 
\end{alignat}
\label{prob:PLD_P}%
\end{subequations}
In the above formulation, constraint~\eqref{prob:PLD_P_perspective_j} is exactly the perspectified constraint derived above, which constraints $x_{k,j}$ for each $j \in [J]$; constraint~\eqref{prob:PLD_P_perspective_j} is the same constraint but for the no-purchase probability $x_{k,0}$. Constraints~\eqref{prob:PLD_P_z_1} - \eqref{prob:PLD_P_z_3} ensure that when the $a_{j,i}$ variables are binary, that $z_{k,j,j',i}$ is indeed forced to be equal to $x_{k,j} \cdot a_{j',i}$, while constraint~\eqref{prob:PLD_P_v} subsequently ensures that $v_{k,j,j'}$ is equal to $x_{k,j} \cdot u_{k,j'}$. Note that unlike problem~\eqref{prob:PLD_RA}, and similarly to formulation~\modelP, it is not necessary to include $u_{k,j}$ as an explicit decision variable. 

Lastly, we note that it is possible to apply RPT to formulation~\eqref{prob:PLD_P} to obtain a stronger formulation, similarly to how \modelPRPT improves on \modelP in the single product case. However, the resulting formulation becomes quite cumbersome to write due to the extremely large number of additional decision variables and constraints; we consequently do not pursue this formulation further.

\pagebreak

\section{Mixed-integer second order cone formulation}
\label{sec:SOC}

In this section, we describe an alternate mixed-integer conic formulation of the logit-based SOCPD problem. This formulation is distinct from formulations \modelRA, \modelP and \modelPRPT in that it does not involve the exponential cone; in fact, this formulation is built entirely using the second-order cone. We consequently denote this formulation by \modelSOC. Before continuing, we note that the structure of this formulation is numerically challenging for second-order cone solvers such as Gurobi and Mosek, for reasons that we will describe below. Nevertheless, this formulation is theoretically interesting and serves to illustrate the richness of the logit-based SOCPD problem. Section~\ref{subsec:SOC_model} presents this formulation, while Section~\ref{subsec:SOC_numerics} presents the results of a small set of numerical experiments to compare \modelSOC against our exponential cone-based models (\modelRA, \modelP and \modelPRPT).

\subsection{Formulation \modelSOC}
\label{subsec:SOC_model}

To motivate this formulation, let us start with the observation that one could model the logit-based SOCPD problem using the $x_{k,0}$, $x_{k,1}$ and $u_k$ variables exactly, if one could enforce the following two constraints, where we use $r_k$ to denote an additional hypograph variable:
\begin{align}
& x_{k,0} (1 + r_k) \geq 1, \quad \forall \ k \in [K], \label{eq:xk0_RSOC} \\
& r_k \leq e^{u_k}, \quad \forall \ k \in [K]. \label{eq:rk_hypograph}
\end{align}
Together, the two constraints enforce the constraint $x_{k,0} \geq 1/(1 + e^{u_k})$; since $x_{k,1} = 1 - x_{k,0}$, this ensures that $x_{k,1}$ is exactly the choice probability of customer type $k$. The constraint~\eqref{eq:xk0_RSOC} is a convex constraint (as we shall discuss in some more detail below). However, constraint~\eqref{eq:rk_hypograph} is not, because $e^{u_k}$ is convex in $u_k$. Surprisingly, it turns out that it is possible to model it with a finite number of convex inequality constraints. To illustrate the idea of how this constraint can be modeled, suppose that $n = 5$, in which case the constraint can be equivalently written as 
\begin{equation}
r_k \leq e^{\beta_{k,0} + \sum_{i=1}^5 \beta_{k,i} a_i}.
\end{equation}
Observe that this constraint is equivalent to 
\begin{equation}
r_k \leq e^{\beta_{k,0}} \cdot \prod_{i=1}^5 e^{\beta_{k,i} a_i}. \label{eq:rk_leq_prod}
\end{equation}
Next, observe that each term $e^{\beta_{k,i} a_i}$ can be rewritten as
\begin{equation}
e^{\beta_{k,i} a_i} = e^{\beta_{k,i}} a_i + 1 \cdot (1 - a_i).
\end{equation}
This equivalence holds because the $a_i$ variables are binary, and the two expressions on either side of the above equality are equal whenever $a_i$ is exactly 0 or 1. Thus, constraint~\eqref{eq:rk_leq_prod} is equivalent to
\begin{equation}
r_k \leq e^{\beta_{k,0}} \cdot \prod_{i=1}^5 \left[ e^{\beta_{k,i}} a_i + 1 \cdot (1 - a_i) \right]. \label{eq:SOC_rk_leq_e_uk}
\end{equation}
It now turns out that this constraint can be equivalently re-written using second-order cones. Specifically, recall that the constraint $xy \geq z^2$, for $x, y \geq 0$, can be represented using second-order cones (this constraint is known as a rotated second-order cone constraint); this constraint is equivalent to $x^{1/2} y^{1/2} \geq |z|$. Let us now introduce the auxiliary variables $\tilde{r}_{k,1}, \dots, \tilde{r}_{k,9}$, and consider the following constraints:
\begin{align}
\tilde{r}_{k,5} & = e^{ \beta_{k,1} \cdot 2^3} a_1 + 1 \cdot (1 - a_1), \label{eq:SOC_example_a1}\\
\tilde{r}_{k,6} & = e^{ \beta_{k,2} \cdot 2^3} a_2 + 1 \cdot (1 - a_2), \\
\tilde{r}_{k,7} & = e^{ \beta_{k,3} \cdot 2^3} a_3 + 1 \cdot (1 - a_3), \\
\tilde{r}_{k,8} & = e^{ \beta_{k,4} \cdot 2^3} a_4 + 1 \cdot (1 - a_4), \\
\tilde{r}_{k,9} & = e^{ \beta_{k,5} \cdot 2^1} a_5 + 1 \cdot (1 - a_5), \label{eq:SOC_example_a5} \\
\tilde{r}^{1/2}_{k,5} \cdot \tilde{r}^{1/2}_{k,6} & \geq \tilde{r}_{k,3}, \\
\tilde{r}^{1/2}_{k,7} \cdot \tilde{r}^{1/2}_{k,8} & \geq \tilde{r}_{k,4}, \\
\tilde{r}^{1/2}_{k,3} \cdot \tilde{r}^{1/2}_{k,4} & \geq \tilde{r}_{k,2}, \\
\tilde{r}^{1/2}_{k,9} \cdot \tilde{r}^{1/2}_{k,2} & \geq \tilde{r}_{k,1}, \label{eq:SOC_example_last} \\
r_k & = e^{\beta_{k,0}} \cdot \tilde{r}_{k,1}, \label{eq:SOC_example_addu0}\\
x_{k,0} \cdot (1 + r_k) & \geq 1, \label{eq:SOC_example_xk0_r_rotatedcone}
\end{align}
We draw the reader's attention to two important observations. First, all of the constraints in this set of constraints are either equality constraints, or rotated second order cone constraints; this includes the last constraint~\eqref{eq:SOC_example_xk0_r_rotatedcone}, which enforces a lower bound on $x_{k,0}$. Second, observe that if one combines constraints~\eqref{eq:SOC_example_a1} - \eqref{eq:SOC_example_addu0} and sequentially substitutes in/projects out the $\tilde{r}_{k,j}$ variables, one exactly obtains the constraint~\eqref{eq:SOC_rk_leq_e_uk}. To see this, observe that after performing such substitutions, one obtains that 
\begin{align}
\tilde{r}_{k,1} & \leq (e^{ \beta_{k,1} \cdot 2^3} a_1 + (1 - a_1))^{1/2^3}  \cdot (e^{ \beta_{k,2} \cdot 2^3} a_2 + (1 - a_2))^{1/2^3}  \cdot (e^{ \beta_{k,3} \cdot 2^3} a_3 + (1 - a_3))^{1/2^3} \nonumber \\
& \phantom{\leq} \cdot (e^{ \beta_{k,4} \cdot 2^3} a_4 + (1 - a_4))^{1/2^3} \cdot (e^{ \beta_{k,5} \cdot 2^1} a_5 + (1 - a_5))^{1/2^1} \label{eq:SOC_example_projectedout}
\end{align}
and here, we again observe that each term of the form $(e^{ \beta_{k,i} \cdot 2^d} a_i + (1 - a_i))^{1/2^d}$ is exactly equal to $e^{\beta_{k,i} \cdot 2^d \cdot 1/2^d \cdot a_i} = e^{\beta_{k,i} a_i}$ whenever $a_i$ is exactly 0 or 1. This last expression also explains why the $2^3$ and $2^1$ factors in constraints~\eqref{eq:SOC_example_a1} - \eqref{eq:SOC_example_a5} are needed: in each constraint of the form $\tilde{r}^{1/2}_{k,j} \tilde{r}^{1/2}_{k,j'} \geq \tilde{r}_{k,j''}$, the variables on the left hand side are raised to the power of $1/2$, which can be viewed as shrinking the $\beta_{k,i'}$ coefficients that enter either $\tilde{r}_{k,j}$ or $\tilde{r}_{k,j'}$ by a factor of 1/2. Thus, to counteract this reduction, the $\beta_{k,i}$'s that appear in \eqref{eq:SOC_example_a1} - \eqref{eq:SOC_example_a5} have to be scaled up by either $2^3$ (for $i = 1, \dots, 4$) or $2^1$ (for $i = 5$). 

\begin{figure}
\centering
\includegraphics[width=0.3\textwidth]{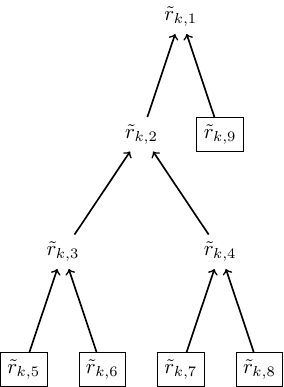}
\caption{Visualization of rotated second order cone constraint for $n = 5$ example in constraints~\eqref{eq:SOC_example_a1} - \eqref{eq:SOC_example_xk0_r_rotatedcone}. \label{figure:SOC_constraint_tree}}
\end{figure}

Figure~\ref{figure:SOC_constraint_tree} visualizes the above constraint set as a binary tree, where each node corresponds to one of the constraints, and the root node corresponds to $\tilde{r}_{k,1}$. Each leaf node corresponds to one of the equality constraints~\eqref{eq:SOC_example_a1} - \eqref{eq:SOC_example_a5}, while each internal node corresponds to one of the rotated second order cone constraints. This visualization also makes clear where the exponent in the $2^d$ expressions in constraints~\eqref{eq:SOC_example_a1} - \eqref{eq:SOC_example_a5} come from: the $d$ in $2^d$ corresponds to the depth of each of the leaf nodes, i.e., the length of the path from the root node to the leaf node. 

With this insight in mind, we now define our general formulation. For each $k \in [K]$, let $\Tcal_k$ denote a binary tree with exactly $n$ leaves. Let $\Lcal_k$ denote the set of leaves of $\Tcal$, and $\Scal_k$ denote the set of internal (non-leaf) nodes of $\Tcal_k$. Each internal node $j \in \Scal_k$ has exactly two child nodes, where $\leftnode(j)$ denotes the left child node, and $\rightnode(j)$ denotes the right child node. Let $d(j)$ denote the depth of each leaf $j \in \Lcal_k$, and let $p(j)$ denote the attribute assigned to leaf $j$. With these definitions, our formulation \modelSOC can be written as follows:
\begin{subequations}
\begin{alignat}{2}
\modelSOC: \quad & \text{maximize} & & \sum_{k=1}^K \lambda_k x_{k,1} \\
& \text{subject to} & \quad & \tilde{r}_{k,j} = e^{\beta_{k,p(j)} \cdot 2^{d(j)}} a_{p(j)} + (1 - a_{p(j)}), \quad \forall \ j \in \Lcal_k, \ k \in [K], \label{prob:SOC_leaves} \\
& & & \tilde{r}^{1/2}_{k,\leftnode(j)} \cdot \tilde{r}^{1/2}_{k,\rightnode(j)} \geq \tilde{r}_{k,j}, \quad \forall \ j \in \Scal_k, \ k \in [K], \label{prob:SOC_rotatedcone} \\
& & & r_{k} = e^{\beta_{k,0}} \cdot \tilde{r}_{k,1}, \quad \forall\ k \in [K], \\
& & & x_{k,0} (1 + r_{k}) \geq 1, \quad \forall \ k \in [K], \\
& & & x_{k,0} + x_{k,1} = 1, \quad \forall \ k \in [K], \\
& & & \Cb \ab \leq \db, \\
& & & a_i \in \{0,1\}, \quad \forall \ i \in [n].
\end{alignat} 
\end{subequations}
We make several important remarks about formulation~\modelSOC. First, as noted above, formulation~\modelSOC is special because it is a mixed-integer second order cone program; by using rotated second order cones and properties of binary variables, one can represent $e^{\beta_{k,0} + \sum_{i=1}^n \beta_{k,i}}$. Stated differently, this formulation illustrates that one can model the logit-based SOCPD problem exactly without the use of exponential cones. 

Second, note that the tree $\Tcal_k$ will generally contain $n$ leaf nodes and $n - 1$ internal nodes. In addition, one can always construct the tree $\Tcal_k$ such that the maximum depth of any leaf node is at most $\lceil \log_2 n \rceil$. (Simply construct a complete tree of depth $\lceil \log_2 n \rceil$, and prune it until one obtains exactly $n$ leaves.) Thus, the factor of $2^{d(j)}$ that appears in constraint~\eqref{prob:SOC_leaves} is at most $2^{\lceil \log_2 n \rceil } \leq 2n$.

Third, as we alluded to earlier, while this formulation is theoretically interesting, it does come with a major limitation. The limitation here comes from the magnitude of the coefficients. As noted in our second remark above, we expect the coefficient $e^{\beta_{k,p(j)} \cdot 2^{d(j)}}$ will be generally equal to $e^{\beta_{k,p(j)} \cdot 2n}$. In the case when $n$ is large and $\beta_{k,p(j)}$ is positive, this coefficient will be very large, and when $n$ is large and $\beta_{k,p(j)}$ is negative, this coefficient will be very small. This leads to numeric issues for mixed-integer second order cone solvers such as Gurobi and Mosek. With Gurobi in particular, large magnitudes of the $\beta_{k,i}$'s and even moderate values of $n$ cause numerical issues (for example, when solving the continuous relaxation, Gurobi erroneously reports the problem to be infeasible or unbounded). Thus, while this formulation is of theoretical interest, we expect that it will in general not be possible to solve. In Section~\ref{subsec:SOC_numerics} of the ecompanion, we present numerical results for a small set of instances with small partworth magnitudes, for which we solve the relaxation of \modelSOC and compare it to \modelRA, \modelP and \modelPRPT. Interestingly, \modelSOC does not lead to a tighter bound than \modelPRPT, but it sometimes leads to a tighter bound than \modelP. 

As an aside to the above two remarks, having a tree $\Tcal_k$ with $2n - 1$ nodes is necessary only when customer type $k$ has exactly $n$ non-zero partworths $\beta_{k,i}$. This is typically the case for choice models obtained from conjoint datasets, where the estimation methods of choice (latent-class MNL and hierarchical Bayesian mixture MNL models) generally do not enforce any kind of sparsity. However, in cases where the partworth matrix $[ \beta_{k,i} ]$ is sparse, then one can modify formulation~\modelSOC so that each customer type is represented with a different tree $\Tcal_k$, and each such tree has $n^*_k$ leaves and $n^*_k - 1$ internal nodes, where $n^*_k < n$ is the number of attributes with non-zero partworths for that customer type. In this case, each constraint tree $\Tcal_k$ can be constructed so that the depth of each leaf is at most $\lceil \log_2 n^*_k \rceil$, and thus the factor of $2^{d(j)}$ will be at most $2n^*_k$. In the case of sparsity, this modification can potentially allow for the aforementioned numerical issue to be partially mitigated. 

Lastly, we note that one can re-write the family of constraints~\eqref{prob:SOC_rotatedcone} using a single power cone constraint. The $m+1$ dimensional power cone is the set 
\begin{equation}
\{ (x_1,\dots, x_m, y) \in \Rbb^{m+1} \mid y \leq x^{\alpha_1}_1 \cdot \dots \cdot x^{\alpha_m}_m; x_1,\dots, x_m, y, \geq 0 \}
\end{equation}
where $\alpha_1,\dots, \alpha_m \geq 0$ and $\sum_{j=1}^m \alpha_j = 1$ \citep{aps2021mosekcookbook}. Thus, one could write 
\begin{equation}
\tilde{r}_{k,1} \leq \prod_{i=1}^n \left[ e^{\beta_{k,i} \cdot n} a_i + (1 - a_i) \right]^{1/n},
\end{equation}
which is power cone-representable. Mosek provides support for conic programs involving power cone constraints and integer variables. However, we encountered errors when attempting to employ Mosek in this manner, and subsequently did not pursue this alternate formulation approach further.

\subsection{Numerical experiments with formulation~\modelSOC} %
\label{subsec:SOC_numerics}

In this section, we present two small numerical experiments to compare formulation \modelSOC with \modelRA, \modelP and \modelPRPT. Recall that \modelSOC differs from \modelRA, \modelP and \modelPRPT in that it models the logit probabilities of each customer type entirely using second-order cone constraints, whereas \modelRA, \modelP and \modelPRPT are based on exponential cone constraints.

In our first experiment, we compare the relaxation bounds of the four formulations. We consider the same overall set of problem instances as in Section~\ref{subsec:numerical_experiments_synthetic}, but with the following restrictions. We restrict to $n = 30$ and vary $K \in \{10,20,30\}$. For the scale parameter $c$, we vary $c \in \{0.1, 0.2, 0.4\}$. The reason for restricting to $n = 30$ and smaller values of $c$ is due to the numerical issue we described in Section~\ref{subsec:SOC_model}; in particular, with $c \geq 0.5$, the relaxation of \modelSOC begins to pose numerical problems for Mosek. 

Additionally, we vary the intercept term $\beta_{k,0}$ in the utility function of each customer type. In particular, we let $\beta_{k,0} = \bar{\beta}$, where we vary $\bar{\beta} \in \{-4, -3, -2, -1, 0\}$. Recall that in our experiments in Section~\ref{subsec:numerical_experiments_synthetic}, $\beta_{k,0}$ is always fixed to -3 across all customer types, and as noted there, one can interpret this choice of constant term as assigning a choice probability of roughly 5\% ($ = \sigma(-3)$) to a product with no attributes ($\ab = \zerob$). By now varying $\bar{\beta}$, we effectively make this no-attribute product more attractive, and the overall market share that one can garner from an optimal product becomes higher. 

By varying the tuple of parameters $(K, c, \bar{\beta})$ in the manner described above, we obtain $3 \times 3 \times 5 = 45$ sets of 20 instances, for a total of 900 instances. For each instance, we solve the continuous relaxations of \modelRA, \modelP, \modelPRPT and \modelSOC. We solve all four formulations using Mosek. 

For \modelSOC, an important aspect of the model is the choice of tree $\Tcal_k$ for each customer type $k$. Here, to ensure numerical stability, we proceed as follows. We construct a complete tree of depth $\lceil \log_2 n \rceil$. We prune this tree back to a tree with exactly $n$ leaves. We sort the partworths in decreasing order of magnitude, i.e., we obtain indices $i_1,\dots, i_n$ such that $|\beta_{k,i_1}| \geq | \beta_{k,i_2}| \geq \dots \geq | \beta_{k,i_n}|$. We then iterate from $t = 1,\dots,n$ and set $p(j) = t$, where $j$ is the leaf with the smallest depth (i.e., closest to the the root) for which $p(j)$ has not yet been set, with ties broken arbitrarily. By setting up the tree in this way, no leaf is at a depth greater than $\lceil \log_2 n \rceil$, which means that the scale factor $2^{d(j)}$ that appears in formulation~\modelSOC is never more than $2n$; additionally, larger attributes are assigned to shallower leaves, which prevents $e^{ \beta_{k,p(j)} \cdot 2^{d(j)}}$ from being too large or too small. We note that this is only one way of setting the trees $\Tcal_1,\dots, \Tcal_K$, and other ways are possible; our goal in setting the tree in this way is to avoid numerical issues. In our preliminary experimentation with \modelSOC, we have found that the choice of trees $\Tcal_1,\dots, \Tcal_K$ can affect the tightness of the relaxation, but the effect appears to be minor. Understanding how the choice of trees $\Tcal_1,\dots, \Tcal_K$ can affect the formulation is an interesting question for future research. 

For each instance, we compute the metric
\begin{equation}
\tilde{I}_m = 100\% \times \frac{Z_{m,rlx} - Z_{best,rlx}}{Z_{best,rlx}},
\end{equation}
where $Z_{m,rlx}$ is the relaxation bound from model $m$, and $Z_{best,rlx} = \min_m Z_{m,rlx}$ is the lowest (i.e., tightest) relaxation bound. In words, this metric measures the gap of the relaxation bound of a given formulation relative to the lowest (i.e., tightest) relaxation bound obtained for each instance.

Table~\ref{table:R2_SOC_relaxations} shows the average value of $\tilde{I}_m$ for all four formulations, and across all values of $(K, c, \bar{\beta})$. We were able to solve the relaxation of \modelSOC to optimality in all 900 instances. Note that in this table, the average values of $\tilde{I}_m$ are generally large because the tightest bound $Z_{best,rlx}$ is usually very small, leading to large gaps. 

From this table, what we find is that \modelSOC often has a lower average gap than \modelP when $c$ is small and $\bar{\beta}$ is small. As either $c$ increases or $\bar{\beta}$ increases, the average gap of \modelP improves and in a small number of parameter settings \modelP has a lower average gap (for example, $K = 10$, $c = 0.4$, $\bar{\beta} = 0$). With regard to \modelRA, we find that in all but two instances out of the 900, \modelSOC is tighter than \modelRA. With regard to \modelPRPT, we find that \modelPRPT gives the best relaxation bound in all 900 instances, and is always tighter than \modelSOC. This experiment thus further underscores the value of \modelPRPT in producing effective bounds for the logit-based SOCPD problem. 

\begin{table}
\centering 
\begin{tabular}{llrrrrr} \toprule
$K$ & $c$ & $\bar{\beta}$ & $\tilde{I}_{\modelRA}$ & $\tilde{I}_{\modelP}$ & $\tilde{I}_{\modelSOC}$ & $\tilde{I}_{\modelPRPT}$ \\ \midrule
10 & 0.1 & -4 & 256.1 & 57.8 & 9.1 & 0.0 \\ 
  10 & 0.1 & -3 & 238.4 & 55.5 & 9.6 & 0.0 \\ 
  10 & 0.1 & -2 & 199.3 & 49.8 & 10.3 & 0.0 \\ 
  10 & 0.1 & -1 & 101.8 & 37.9 & 9.6 & 0.0 \\ 
  10 & 0.1 &  0 & 45.1 & 13.0 & 6.5 & 0.0 \\ 
  10 & 0.2 & -4 & 359.4 & 87.7 & 17.9 & 0.0 \\ 
  10 & 0.2 & -3 & 312.3 & 85.0 & 21.9 & 0.0 \\ 
  10 & 0.2 & -2 & 176.9 & 77.0 & 27.5 & 0.0 \\ 
  10 & 0.2 & -1 & 92.7 & 43.7 & 29.2 & 0.0 \\ 
  10 & 0.2 &  0 & 43.5 & 15.3 & 19.6 & 0.0 \\ 
  10 & 0.4 & -4 & 244.7 & 97.9 & 33.4 & 0.0 \\ 
  10 & 0.4 & -3 & 147.6 & 82.2 & 49.9 & 0.0 \\ 
  10 & 0.4 & -2 & 98.3 & 54.0 & 60.3 & 0.0 \\ 
  10 & 0.4 & -1 & 62.2 & 29.5 & 50.0 & 0.0 \\ 
  10 & 0.4 &  0 & 32.8 & 12.7 & 29.2 & 0.0 \\ \midrule
  20 & 0.1 & -4 & 259.2 & 61.0 & 8.4 & 0.0 \\ 
  20 & 0.1 & -3 & 242.9 & 58.9 & 9.2 & 0.0 \\ 
  20 & 0.1 & -2 & 205.6 & 53.7 & 10.6 & 0.0 \\ 
  20 & 0.1 & -1 & 106.5 & 42.1 & 11.2 & 0.0 \\ 
  20 & 0.1 &  0 & 47.7 & 15.4 & 8.3 & 0.0 \\ 
  20 & 0.2 & -4 & 359.3 & 89.1 & 15.7 & 0.0 \\ 
  20 & 0.2 & -3 & 310.6 & 86.9 & 20.3 & 0.0 \\ 
  20 & 0.2 & -2 & 179.5 & 79.5 & 27.3 & 0.0 \\ 
  20 & 0.2 & -1 & 96.4 & 48.2 & 31.9 & 0.0 \\ 
  20 & 0.2 &  0 & 47.7 & 18.7 & 23.8 & 0.0 \\ 
  20 & 0.4 & -4 & 243.5 & 96.7 & 32.4 & 0.0 \\ 
  20 & 0.4 & -3 & 147.9 & 83.9 & 49.6 & 0.0 \\ 
  20 & 0.4 & -2 & 100.2 & 58.0 & 63.3 & 0.0 \\ 
  20 & 0.4 & -1 & 66.5 & 33.8 & 56.4 & 0.0 \\ 
  20 & 0.4 &  0 & 38.1 & 16.3 & 36.0 & 0.0 \\ \midrule
  30 & 0.1 & -4 & 258.1 & 61.6 & 7.8 & 0.0 \\ 
  30 & 0.1 & -3 & 243.2 & 59.7 & 8.8 & 0.0 \\ 
  30 & 0.1 & -2 & 208.5 & 54.8 & 10.5 & 0.0 \\ 
  30 & 0.1 & -1 & 108.5 & 43.8 & 11.7 & 0.0 \\ 
  30 & 0.1 &  0 & 49.1 & 16.7 & 9.2 & 0.0 \\ 
  30 & 0.2 & -4 & 359.1 & 89.3 & 14.3 & 0.0 \\ 
  30 & 0.2 & -3 & 314.5 & 87.3 & 19.1 & 0.0 \\ 
  30 & 0.2 & -2 & 181.7 & 80.8 & 26.8 & 0.0 \\ 
  30 & 0.2 & -1 & 98.0 & 50.0 & 32.6 & 0.0 \\ 
  30 & 0.2 &  0 & 49.8 & 20.3 & 25.7 & 0.0 \\ 
  30 & 0.4 & -4 & 245.7 & 97.3 & 29.9 & 0.0 \\ 
  30 & 0.4 & -3 & 148.6 & 84.7 & 48.2 & 0.0 \\ 
  30 & 0.4 & -2 & 101.3 & 59.6 & 64.2 & 0.0 \\ 
  30 & 0.4 & -1 & 68.6 & 35.7 & 59.3 & 0.0 \\ 
  30 & 0.4 &  0 & 40.6 & 18.0 & 39.2 & 0.0 \\ \bottomrule 
\end{tabular}
\caption{Comparison of relaxation bounds, measured as gaps relative to best relaxation bound for each instance, for \modelRA, \modelP, \modelPRPT and \modelSOC. \label{table:R2_SOC_relaxations} }
\end{table}

In our second experiment, we compare the optimality gap of solving model~\modelSOC as a mixed-integer conic program, and compare its optimality gap with that of models \modelRA, \modelP and \modelPRPT. We restrict our attention to $n = 30$ and $\bar{\beta} = -3$, and vary $K \in \{10,20,30\}$ and $c \in \{0.1, 0.2, 0.4\}$. We solve both \modelSOC and \modelP using Mosek, with a time limit of one hour. For each instance, we calculate the optimality gap $O_{m}$ of formulation $m$ in the same way as in Section~\ref{subsec:numerical_experiments_synthetic}, and additionally record the computation time $T_m$.

Table~\ref{table:R2_SOC_MIPs} below displays the average optimality gap as a percentage for both \modelSOC and \modelP and the average computation time for each formulation, as $K$ and $c$ vary. From this table, we can see that in general, \modelSOC is able to attain a lower optimality gap than \modelP and \modelRA within a one hour time limit, and for the cases where it is possible to solve the problem to full optimality, the time required for \modelSOC is lower than for \modelP and \modelRA. With regard to \modelPRPT, when the instances can be solved within the one hour time limit, it appears that \modelPRPT is faster (e.g., $K = 10$, $c = 0.4$), while for larger instances, the two formulations seem to exhibit similar performance.

\begin{table}
\begin{center}
\begin{tabular}{llrrrrrrrr} \toprule
$K$ & $c$ & $O_{\modelRA}$ & $O_{\modelP}$ & $O_{\modelPRPT}$ & $O_{\modelSOC}$ & $T_{\modelRA}$ & $T_{\modelP}$ & $T_{\modelPRPT}$ & $T_{\modelSOC}$ \\ 
  \hline
 10 & 0.10 & 48.21 & 3.46 & 0.00 & 0.00 & 3617.73 & 2897.79 & 70.50 & 52.17 \\ 
   10 & 0.20 & 56.77 & 6.99 & 0.00 & 0.00 & 3622.70 & 2902.37 & 350.10 & 543.53 \\ 
   10 & 0.40 & 56.85 & 0.50 & 0.00 & 0.00 & 3625.64 & 1844.50 & 718.98 & 1178.93 \\[0.25em]
   20 & 0.10 & 53.47 & 13.12 & 0.00 & 0.00 & 3615.60 & 3615.05 & 470.06 & 812.32 \\ 
   20 & 0.20 & 63.11 & 22.95 & 4.19 & 3.29 & 3616.22 & 3614.74 & 3225.59 & 2932.75 \\ 
   20 & 0.40 & 64.92 & 27.21 & 22.95 & 23.12 & 3619.70 & 3614.63 & 3601.47 & 3601.11 \\[0.25em]
   30 & 0.10 & 56.03 & 17.27 & 0.00 & 1.61 & 3614.84 & 3615.24 & 1572.22 & 2953.90 \\ 
   30 & 0.20 & 67.80 & 31.04 & 15.81 & 14.75 & 3614.80 & 3615.55 & 3601.66 & 3601.38 \\ 
   30 & 0.40 & 71.53 & 38.62 & 33.13 & 34.97 & 3615.95 & 3616.01 & 3601.64 & 3600.95 \\ \bottomrule
\end{tabular}
\end{center}
\caption{Comparison of optimality gaps and computation times for integer versions of \modelRA, \modelP, \modelPRPT and \modelSOC. \label{table:R2_SOC_MIPs} }
\end{table}

As noted before, we expect \modelP to be more useful than \modelSOC in most cases. However, in certain cases where the partworths are sufficiently small or are sparse, \modelSOC could potentially be useful. Whether or not mixed-integer second order cone program solvers can be designed to robustly solve \modelSOC, and whether formulation~\modelSOC can be further modified to sidestep the numerical magnitude issue, are both interesting questions for future study.

\clearpage
\pagebreak

\section{Additional details for numerical experiments in Section~\ref{sec:numerical_experiments}}

\subsection{Computation times and optimality gaps for \modelRA and \modelPRPT on synthetic instances}
\label{subsec:R2_PRPT_gap_time}

In this section, we provide additional results for the optimality gap and computation time of formulations~\modelPRPT and \modelRA when they are used to solve the integer versions of the synthetic instances used in Section~\ref{subsec:numerical_experiments_synthetic}. 

Table~\ref{table:R2_RA_gap_time} shows the optimality gap and computation time for formulation~\modelRA. While formulation~\modelRA is generally able to achieve a low average gap across all the instance sets and in many cases is able to solve instances to full optimality within a two hour time limit, it generally has larger computation times and gaps than \modelP. Additionally, in some instances Mosek encountered fatal errors ({\tt MOSEK fatal error stoptask} and {\tt MOSEK fatal error stopenv} in {\tt src/prslv/prlog.c}). This is likely due to the ill-posedness issue highlighted in Section~\ref{subsec:micp_P} (namely that the principal constraint in \modelRA must hold at equality for integer solutions). The $(c,n,K)$ combinations for which this occurred are indicated in the table by asterisks. In total, there were 6 instances out of the 900 total instances where this occurred: 1 instance with $(c,n,K) = (20,50,30)$; 1 instance with $(c,n,K) = (20,40,30)$; 1 instance with $(c,n,K) = (20,40,20)$; 1 instance with $(c,n,K) = (10,40,30)$; and 2 instances with $(c,n,K) = (10,30,30)$.

Table~\ref{table:R2_PRPT_gap_time} shows the optimality gap and computation time for formulation~\modelPRPT. Overall, although formulation~\modelPRPT provides a stronger relaxation bound compared to \modelP, it is generally more difficult to solve it as an integer program than \modelP. In particular, the average optimality gap and computation time for \modelPRPT are generally larger than for \modelP across all of the $(c,n,K)$ combinations. Of the three formulations - \modelRA, \modelP and \modelPRPT - it appears that \modelPRPT is generally the slowest. Furthermore, in some instances, Mosek encountered a fatal error (code 5010) when solving the integer version of \modelPRPT; the $(c,n,K)$ combinations for which this occurred are indicated in Table~\ref{table:R2_PRPT_gap_time} by asterisks. In total, there were 21 instances out of the 900 total instances where this occurred: 1 instance with $(c,n,K) = (5,50,30)$; 8 instances with $(c,n,K) = (10,50,30)$; and 12 instances with $(c,n,K) = (20,50,30)$.

\begin{table}
\begin{tabular}{lllrr} \toprule
$c$ & $n$ & $K$ & $O_{\modelRA}$ & $T_{\modelRA}$ \\ \midrule
 5 &  30 &  10 & 0.00 & 147.87 \\ 
    5 &  30 &  20 & 0.28 & 2633.99 \\ 
    5 &  30 &  30 & 9.82 & 6685.95 \\ 
    5 &  40 &  10 & 0.01 & 2660.31 \\ 
    5 &  40 &  20 & 4.83 & 7068.59 \\ 
    5 &  40 &  30 & 13.91 & 7213.60 \\ 
    5 &  50 &  10 & 0.04 & 1853.00 \\ 
    5 &  50 &  20 & 5.01 & 7211.76 \\ 
    5 &  50 &  30 & 19.48 & 7212.95 \\ 
    5 &  60 &  10 & 0.03 & 749.20 \\ 
    5 &  60 &  20 & 2.57 & 5877.82 \\ 
    5 &  60 &  30 & 16.19 & 7214.92 \\ 
    5 &  70 &  10 & 0.01 & 372.67 \\ 
    5 &  70 &  20 & 0.92 & 6132.60 \\ 
    5 &  70 &  30 & 8.82 & 7212.24 \\ \bottomrule
    \end{tabular} \qquad
\begin{tabular}{lllrr} \toprule
$c$ & $n$ & $K$ & $O_{\modelRA}$ & $T_{\modelRA}$ \\ \midrule
   10 &  30 &  10 & 0.00 & 37.87 \\ 
   10 &  30 &  20 & 0.00 & 489.93 \\ 
   10 &  30 &  30 & 4.45$^*$ & 6058.49$^*$ \\ 
   10 &  40 &  10 & 0.00 & 71.34 \\ 
   10 &  40 &  20 & 0.47 & 3475.42 \\ 
   10 &  40 &  30 & 6.30$^*$ & 7106.71$^*$ \\ 
   10 &  50 &  10 & 0.00 & 22.96 \\ 
   10 &  50 &  20 & 1.55 & 5016.24 \\ 
   10 &  50 &  30 & 11.77 & 7212.94 \\ 
   10 &  60 &  10 & 0.00 & 38.98 \\ 
   10 &  60 &  20 & 0.27 & 2346.55 \\ 
   10 &  60 &  30 & 6.80 & 7211.44 \\ 
   10 &  70 &  10 & 0.00 & 26.36 \\ 
   10 &  70 &  20 & 0.02 & 1733.51 \\ 
   10 &  70 &  30 & 2.50 & 5268.35 \\ \bottomrule
   \end{tabular} \qquad
\begin{tabular}{lllrr} \toprule
$c$ & $n$ & $K$ & $O_{\modelRA}$ & $T_{\modelRA}$ \\ \midrule
   20 &  30 &  10 & 0.00 & 16.81 \\ 
   20 &  30 &  20 & 0.00 & 385.26 \\ 
   20 &  30 &  30 & 1.75 & 4610.50 \\ 
   20 &  40 &  10 & 0.00 & 28.27 \\ 
   20 &  40 &  20 & 0.29$^*$ & 2193.60$^*$ \\ 
   20 &  40 &  30 & 4.35$^*$ & 6468.67$^*$ \\ 
   20 &  50 &  10 & 0.00 & 28.29 \\ 
   20 &  50 &  20 & 0.53 & 1747.12 \\ 
   20 &  50 &  30 & 8.47$^*$ & 6663.90$^*$ \\ 
   20 &  60 &  10 & 0.00 & 33.18 \\ 
   20 &  60 &  20 & 0.00 & 835.86 \\ 
   20 &  60 &  30 & 3.85 & 6679.17 \\ 
   20 &  70 &  10 & 0.00 & 38.67 \\ 
   20 &  70 &  20 & 0.00 & 812.00 \\ 
   20 &  70 &  30 & 1.01 & 5212.16 \\ \bottomrule
\end{tabular}
\caption{Optimality gap and computation time of formulation \modelRA as $c$, $n$ and $K$ vary. Asterisks denote $(c,n,K)$ combinations for which Mosek encountered errors (see text for details). \label{table:R2_RA_gap_time} }
\end{table}

\begin{table}
\begin{tabular}{lllrr} \toprule
$c$ & $n$ & $K$ & $O_{\modelPRPT}$ & $T_{\modelPRPT}$ \\ \midrule
5 &  30 &  10 & 0.00 & 125.44 \\ 
    5 &  30 &  20 & 0.51 & 3773.07 \\ 
    5 &  30 &  30 & 7.10 & 6816.63 \\ 
    5 &  40 &  10 & 0.00 & 641.40 \\ 
    5 &  40 &  20 & 5.47 & 6452.49 \\ 
    5 &  40 &  30 & 19.15 & 7213.38 \\ 
    5 &  50 &  10 & 0.00 & 366.49 \\ 
    5 &  50 &  20 & 3.80 & 6605.27 \\ 
    5 &  50 &  30 & 20.54$^*$ & 7214.95$^*$ \\ 
    5 &  60 &  10 & 0.00 & 649.10 \\ 
    5 &  60 &  20 & 5.80 & 6578.93 \\ 
    5 &  60 &  30 & 34.89 & 7217.10 \\ 
    5 &  70 &  10 & 0.00 & 978.12 \\ 
    5 &  70 &  20 & 3.21 & 5706.16 \\ 
    5 &  70 &  30 & 24.20 & 7202.02 \\ \bottomrule
\end{tabular} \qquad
\begin{tabular}{lllrr} \toprule
$c$ & $n$ & $K$ & $O_{\modelPRPT}$ & $T_{\modelPRPT}$ \\ \midrule
   10 &  30 &  10 & 0.00 & 122.70 \\ 
   10 &  30 &  20 & 0.17 & 2713.76 \\ 
   10 &  30 &  30 & 9.81 & 7088.07 \\ 
   10 &  40 &  10 & 0.00 & 500.30 \\ 
   10 &  40 &  20 & 6.95 & 6438.01 \\ 
   10 &  40 &  30 & 17.41 & 7065.37 \\ 
   10 &  50 &  10 & 0.00 & 687.43 \\ 
   10 &  50 &  20 & 5.60 & 6654.84 \\ 
   10 &  50 &  30 & 19.82$^*$ & 7215.15$^*$ \\ 
   10 &  60 &  10 & 0.00 & 1195.47 \\ 
   10 &  60 &  20 & 8.92 & 7092.99 \\ 
   10 &  60 &  30 & 59.05 & 7218.98 \\ 
   10 &  70 &  10 & 0.00 & 1787.23 \\ 
   10 &  70 &  20 & 7.38 & 6174.05 \\ 
   10 &  70 &  30 & 81.86 & 7222.57 \\ \bottomrule
\end{tabular} \qquad
\begin{tabular}{lllrr} \toprule
$c$ & $n$ & $K$ & $O_{\modelPRPT}$ & $T_{\modelPRPT}$ \\ \midrule
   20 &  30 &  10 & 0.00 & 165.68 \\ 
   20 &  30 &  20 & 0.68 & 3147.87 \\ 
   20 &  30 &  30 & 13.52 & 7042.75 \\ 
   20 &  40 &  10 & 0.00 & 454.27 \\ 
   20 &  40 &  20 & 8.34 & 6712.30 \\ 
   20 &  40 &  30 & 18.46 & 7216.36 \\ 
   20 &  50 &  10 & 0.00 & 879.30 \\ 
   20 &  50 &  20 & 9.37 & 6994.03 \\ 
   20 &  50 &  30 & 19.72$^*$ & 7218.72$^*$ \\ 
   20 &  60 &  10 & 0.03 & 1745.49 \\ 
   20 &  60 &  20 & 14.76 & 7179.80 \\ 
   20 &  60 &  30 & 83.74 & 7223.19 \\ 
   20 &  70 &  10 & 1.61 & 3232.79 \\ 
   20 &  70 &  20 & 10.20 & 6869.17 \\ 
   20 &  70 &  30 & 81.09 & 7116.10 \\ \bottomrule
\end{tabular}
\caption{Optimality gap and computation time of formulation \modelPRPT as $c$, $n$ and $K$ vary. Asterisks denote $(c,n,K)$ combinations for which Mosek encountered errors (see text for details). \label{table:R2_PRPT_gap_time} }
\end{table}

\subsection{Comparison of computation times for heuristic approaches and formulation \modelP on synthetic instances}
\label{subsec:R1_synthetic_time}

Table~\ref{table:R1_synthetic_time} below compares the average computation time (in seconds), where the average is taken over the 20 instances for each $(c,n,K)$ triple, for each of the heuristics (\Greedy, \LS, \KKDP, \GM) and formulation~\modelP. 

\begin{table}
\centering
\begin{tabular}{rrrrrrrr} \toprule
$c$ & $n$ & $K$ & $T_{\Greedy}$ & $T_{\LS}$ & $T_{\KKDP}$ & $T_{\GM}$ & $T_{\modelP}$ \\ \midrule
5 &  30 &  10 & 0.00 & 0.00 & 0.00 & 0.07 & 15.54 \\ 
    5 &  30 &  20 & 0.00 & 0.00 & 0.00 & 0.11 & 202.62 \\ 
    5 &  30 &  30 & 0.00 & 0.00 & 0.00 & 0.19 & 2333.12 \\ 
    5 &  40 &  10 & 0.00 & 0.01 & 0.00 & 0.10 & 44.99 \\ 
    5 &  40 &  20 & 0.00 & 0.01 & 0.00 & 18.56 & 2393.76 \\ 
    5 &  40 &  30 & 0.00 & 0.00 & 0.00 & 0.30 & 7026.27 \\ 
    5 &  50 &  10 & 0.00 & 0.01 & 0.00 & 0.23 & 21.62 \\ 
    5 &  50 &  20 & 0.00 & 0.01 & 0.00 & 0.43 & 2809.89 \\ 
    5 &  50 &  30 & 0.00 & 0.01 & 0.00 & 0.82 & 7213.79 \\ 
    5 &  60 &  10 & 0.00 & 0.02 & 0.00 & 0.42 & 47.54 \\ 
    5 &  60 &  20 & 0.00 & 0.02 & 0.00 & 3.25 & 2317.13 \\ 
    5 &  60 &  30 & 0.00 & 0.02 & 0.00 & 2.88 & 7213.76 \\ 
    5 &  70 &  10 & 0.00 & 0.03 & 0.00 & 0.67 & 18.33 \\ 
    5 &  70 &  20 & 0.00 & 0.03 & 0.00 & 3.70 & 1267.83 \\ 
    5 &  70 &  30 & 0.00 & 0.03 & 0.00 & 7.60 & 6410.63 \\  \midrule
   10 &  30 &  10 & 0.00 & 0.00 & 0.00 & 0.20 & 14.96 \\ 
   10 &  30 &  20 & 0.00 & 0.00 & 0.00 & 0.21 & 80.05 \\ 
   10 &  30 &  30 & 0.00 & 0.00 & 0.00 & 0.27 & 1107.93 \\ 
   10 &  40 &  10 & 0.00 & 0.00 & 0.00 & 0.32 & 17.78 \\ 
   10 &  40 &  20 & 0.00 & 0.00 & 0.00 & 0.76 & 399.22 \\ 
   10 &  40 &  30 & 0.00 & 0.00 & 0.00 & 0.63 & 5316.60 \\ 
   10 &  50 &  10 & 0.00 & 0.01 & 0.00 & 0.38 & 17.92 \\ 
   10 &  50 &  20 & 0.00 & 0.01 & 0.00 & 2.56 & 1351.77 \\ 
   10 &  50 &  30 & 0.00 & 0.01 & 0.00 & 4.41 & 7161.57 \\ 
   10 &  60 &  10 & 0.00 & 0.02 & 0.00 & 1.25 & 22.60 \\ 
   10 &  60 &  20 & 0.00 & 0.02 & 0.00 & 15.24 & 709.22 \\ 
   10 &  60 &  30 & 0.00 & 0.01 & 0.00 & 21.05 & 6692.13 \\ 
   10 &  70 &  10 & 0.00 & 0.02 & 0.00 & 2.84 & 16.69 \\ 
   10 &  70 &  20 & 0.00 & 0.02 & 0.00 & 14.89 & 87.16 \\ 
   10 &  70 &  30 & 0.00 & 0.01 & 0.00 & 143.68 & 3310.76 \\ \midrule
   20 &  30 &  10 & 0.00 & 0.00 & 0.00 & 1.57 & 12.23 \\ 
   20 &  30 &  20 & 0.00 & 0.00 & 0.00 & 0.28 & 66.61 \\ 
   20 &  30 &  30 & 0.00 & 0.00 & 0.00 & 0.28 & 614.31 \\ 
   20 &  40 &  10 & 0.00 & 0.00 & 0.00 & 0.36 & 18.42 \\ 
   20 &  40 &  20 & 0.00 & 0.00 & 0.00 & 26.03 & 496.22 \\ 
   20 &  40 &  30 & 0.00 & 0.00 & 0.00 & 1.39 & 4346.31 \\ 
   20 &  50 &  10 & 0.00 & 0.01 & 0.00 & 0.51 & 24.78 \\ 
   20 &  50 &  20 & 0.00 & 0.01 & 0.00 & 4.84 & 817.61 \\ 
   20 &  50 &  30 & 0.00 & 0.02 & 0.00 & 23.71 & 5550.03 \\ 
   20 &  60 &  10 & 0.00 & 0.01 & 0.00 & 1.22 & 38.61 \\ 
   20 &  60 &  20 & 0.00 & 0.01 & 0.00 & 6.87 & 114.48 \\ 
   20 &  60 &  30 & 0.00 & 0.01 & 0.00 & 82.74 & 5115.71 \\ 
   20 &  70 &  10 & 0.00 & 0.01 & 0.00 & 0.93 & 21.44 \\ 
   20 &  70 &  20 & 0.00 & 0.02 & 0.00 & 12.15 & 202.43 \\ 
   20 &  70 &  30 & 0.00 & 0.01 & 0.00 & 113.62 & 2610.99 \\  \bottomrule
\end{tabular}
\caption{Comparison of computation times for heuristic approaches and formulation~\modelP on synthetic instances. \label{table:R1_synthetic_time} }  
\end{table}

\subsection{Additional details on real data sets in Section~\ref{subsec:numerical_experiments_real}}
\label{subsec:R1_real_data_details}

In this section, we provide some additional details on the four data sets used in Section~\ref{subsec:numerical_experiments_real}. As noted in Section~\ref{subsec:numerical_experiments_real}, these four data sets are conjoint analysis data sets, and specifically choice-based conjoint data sets. In choice-based conjoint analysis, a respondent is shown two or more hypothetical products formulated in terms of the attributes that are being studied, and is asked to choose between them. The choice between these hypothetical products (also known as \emph{profiles}) is called a \emph{task}. Based on the responses given by each customer to each of their tasks, one can estimate a discrete choice model, such as a latent-class logit model, that predicts how the customer will choose and provides a measure of the utility for each attribute. As an alternative to choice-based conjoint analysis, there also exists what is called ratings-based or metric conjoint analysis, where a customer is shown a single profile and asked to provide a numeric rating. Based on the responses to such rating tasks, one can use ordinary least squares to determine the utility of each attribute.

In all four data sets (\bank, \candidate, \immigrant and \toubia), each task consists of choosing between two profiles, which is also known as a paired comparison task in the conjoint analysis literature. The number of paired comparison tasks varies for each data set. For \bank, each respondent performed between 14 and 17 paired comparison tasks; for \candidate, between 3 and 6 tasks; for \immigrant, exactly 5 tasks; and for \toubia, exactly 16 tasks. (We note that for \toubia, the paired comparison also included a metric/rating component, where respondents were asked to specify the degree to which one profile was preferred to the other. In this experiment, respondents were allowed to specify being indifferent between the two profiles; this happened in 297 out of 5280 total responses. Since the estimation of latent-class and mixture MNL models requires a choice and since this indifference happened in a relatively small number of responses, we removed these responses from the \toubia data set when conducting our estimation procedures.)

While conjoint studies sometimes involve tasks where respondents can select a no-purchase option (for example, a respondent is shown two or more profiles and a ``none of the above'' option), none of the four data sets we used include explicit information on the no-purchase option, and none of them included any task where the respondent was asked to choose between a product profile and the no-purchase option. Thus, we instead assume the existence of several competitive products, and assume that the customer is allowed to choose the product we have designed or one of the competitive products. The utility of each of the competitive offerings is calculated using the same partworths that are used to calculate the utility of the product we are designing. To provide an example, suppose that $\ab'$, $\ab''$, $\ab'''$ are the attribute vectors of three competitive products. If $\ab$ is the attribute vector of our product, then the no-purchase probability of customer type $k$ would be given by 
\begin{equation}
\frac{ e^{ \sum_{i=1}^n \beta_{k,i} a'_i } + e^{ \sum_{i=1}^n \beta_{k,i} a''_i } + e^{ \sum_{i=1}^n \beta_{k,i} a'''_i} }{ e^{\sum_{i=1}^n \beta_{k,i} a_i} + e^{ \sum_{i=1}^n \beta_{k,i} a'_i } + e^{ \sum_{i=1}^n \beta_{k,i} a''_i } + e^{ \sum_{i=1}^n \beta_{k,i} a'''_i} }. \label{eq:no_purchase_competitive_products}
\end{equation}
Recall that the no-purchase probability under the model described in Section~\ref{subsec:model_problem_definition} can be expressed as
\begin{align}
& \frac{ 1 }{1 + e^{ \beta_{k,0} + \sum_{i=1}^n \beta_{k,i} a_i}} \nonumber \\
& = \frac{ e^{-\beta_{k,0} } }{e^{- \beta_{k,0}} + e^{ \sum_{i=1}^n \beta_{k,i} a_i}}. \label{eq:no_purchase_betak0}
\end{align}
Thus, to calibrate $\beta_{k,0}$ so that \eqref{eq:no_purchase_competitive_products} and \eqref{eq:no_purchase_betak0} are equal, we simply set $\beta_{k,0}$ as 
\begin{equation}
\beta_{k,0} = - \log \left( e^{ \sum_{i=1}^n \beta_{k,i} a'_i } + e^{ \sum_{i=1}^n \beta_{k,i} a''_i } + e^{ \sum_{i=1}^n \beta_{k,i} a'''_i}  \right) .
\end{equation}

\subsection{Attributes for real data instances in Section~\ref{subsec:numerical_experiments_real}}

Tables~\ref{table:bank_attributes}, \ref{table:candidate_attributes}, \ref{table:immigrant_attributes} and \ref{table:toubia_attributes} display the attributes and attribute levels for the \bank, \candidate, \immigrant and \toubia datasets, respectively.

\begin{table}[!ht]
\centering
\begin{tabular}{ll} \toprule
Attribute & Levels \\ \midrule
Interest Rate & High Fixed Rate, Medium Fixed Rate, \\
& Low Fixed Rate, Medium Variable Rate \\
Rewards & 1, 2, 3, 4 \\
Annual Fee & High, Medium, Low \\
Bank & Bank A, Bank B, Out of State Bank \\
Rebate & Low, Medium, High \\
Credit Line & Low, High \\
Grace Period & Short, Long \\ \bottomrule
\end{tabular}
\caption{Attributes for \bank dataset. \label{table:bank_attributes}}
\end{table}

\begin{table}[!ht]
\centering
\begin{tabular}{ll} \toprule
Attribute & Levels \\ \midrule
Age & 36, 45, 52, 60, 68, 75 \\
Military Service & Did Not serve, Served \\
Religion & None, Jewish, Catholic, Mainline Protestant, \\
& Evangelical Protestant, Mormon  \\
College & No BA, Baptist College, Community College, \\
&  State University, Small College, \\
& Ivy League University \\
Income & 32K, 54K, 65K, 92K, 210K, 5.1M \\
Profession & Business Owner, Lawyer, Doctor, \\
& High School Teacher, Farmer, Car Dealer \\
Race/Ethnicity & White, Native American, Black, \\
& Hispanic, Caucasian, Asian American \\
Gender & Male, Female  \\ \bottomrule
\end{tabular}
\caption{Attributes for \candidate dataset. \label{table:candidate_attributes} }
\end{table}

\begin{table}[!ht]
\centering
\begin{tabular}{ll} \toprule
Attribute & Levels \\ \midrule
Education & No Formal, 4th Grade, 8th Grade, High School, \\
& Two-Year College, College Degree, Graduate Degree \\
Gender & Female, Male \\
Origin & Germany, France, Mexico, Philippines, Poland, \\
& India, China, Sudan, Somalia, Iraq \\
Application Reason & Reunite With Family, Seek Better Job, Escape Persecution \\
Profession & Janitor, Waiter, Child Care Provider, Gardener, Financial Analyst, \\
& Construction Worker, Teacher, Computer Programmer, \\
& Nurse, Research Scientist, Doctor \\
Job Experience & None, 1-2 Years, 3-5 Years, 5+ Years \\
Job Plans & Contract With Employer, Interviews With Employer, \\
& Will Look For Work, No Plans To Look For Work \\
Prior Trips to US & Never, Once As Tourist, Many Times As Tourist, \\
& Six Months With Family, Once Without Authorization \\
Language & Fluent English, Broken English, \\
&  Tried English But Unable, Used Interpreter \\ \bottomrule
\end{tabular}
\caption{Attributes for \immigrant dataset. \label{table:immigrant_attributes}}
\end{table}

\begin{table}[!ht]
\centering
\begin{tabular}{ll} \toprule
Attribute & Levels \\ \midrule
Price & \$70, \$75, \$80, \$85, \$90, \$95, \$100 \\
Size & Normal, Large \\
Color & Black, Red \\
Logo & No, Yes \\
Handle & No, Yes \\
PDA Holder & No, Yes \\
Cellphone Holder & No, Yes \\
Velcro Flap & No, Yes \\
Protective Boot & No, Yes \\ \bottomrule
\end{tabular}
\caption{Attributes for \toubia dataset. \label{table:toubia_attributes} }
\end{table}

\clearpage
\pagebreak

\subsection{Hierarchical Bayesian model specification}
\label{appendix:numerics_HB_specification}

For our hierarchical Bayesian model, we assume that each respondent's partworth vector $\betab = (\beta_1,\dots, \beta_n)$ is drawn as
\begin{equation}
\betab \sim N( \bar{\betab}, \mathbf{V}_{\betab}), \label{eq:HB_first_stage_betab}
\end{equation}
where $N(\mub, \Sigmab)$ denotes a multivariate normal distribution with mean $\mu$ and covariance matrix $\Sigmab$. The distributions of the mean $\bar{\betab}$ and covariance matrix $\mathbf{V}_{\betab}$ are then specified as
\begin{align}
\bar{\betab} \sim N( \zerob, \alpha \mathbf{V}_{\betab}), \label{eq:HB_second_stage_barbetab}\\
\mathbf{V}_{\betab} \sim IW(\nu, \mathbf{V} ), \label{eq:HB_second_stage_V}
\end{align}
where $IW(\nu, \mathbf{W})$ denotes an inverse Wishart distribution with degrees of freedom $\nu$ and scale matrix $\mathbf{W}$. This model specification is implemented in {\tt bayesm}, using the {\tt rhierBinLogit} function. We use {\tt bayesm}'s defaults for $\nu$, $\mathbf{V}$ and $\alpha$.

\subsection{Additional constraints for \immigrant dataset}
\label{appendix:numerics_immigrant_constraints}

As discussed  in Section~\ref{subsec:numerical_experiments_real}, we define $\Acal$ with some additional constraints, which we describe here:
\begin{itemize}
\item If the immigrant's profession attribute is set to ``doctor'', ``research scientist'', ``computer programmer'' or ``financial analyst'', then the immigrant's education attribute is set to ``college degree'' or ``graduate degree''.
\item If the immigrant's profession attribute is set to ``teacher'' or ``nurse'', then the immigrant's education attribute is set to ``high school'', ``two-year college'', ``college degree'' or ``graduate degree''. 
\item If the immigrant's application reason attribute is set to ``escape persecution'', then the country of origin attribute is set to Sudan, Somalia or Iraq. 
\item Either the immigrant's application reason attribute is set to ``seek better job'' or the immigrant's job plan attribute is set to ``no plans to work'', but they cannot both be set in this way. 
\end{itemize}

\subsection{Competitive offerings for Section~\ref{subsec:numerical_experiments_real}}
\label{appendix:numerics_competitive_offerings}

Tables~\ref{table:bank_nopurchase}, \ref{table:candidate_nopurchase}, \ref{table:immigrant_nopurchase} and \ref{table:toubia_nopurchase} display the attributes of the competitive offerings for the \bank, \candidate, \immigrant and \toubia datasets, respectively. We note that for \toubia, we follow the same competitive offerings used in other optimization work that has used this dataset \citep{belloni2008optimizing,bertsimas2017robust,bertsimas2019exact}.

\begin{table}[!ht]
\centering
\scriptsize
\begin{tabular}{|l|c|c|c|} \hline
Attribute & Outside  & Outside  &  Outside  \\
& Option 1 & Option 2 & Option 3 \\ \hline
Interest Rate: High fixed rate &  &  &  \\ \hline
Interest Rate: Medium fixed rate & \cellcolor{black!25} &  &  \\ \hline
Interest Rate: Low fixed rate &  & \cellcolor{black!25} &  \\ \hline
Interest Rate: Medium variable rate &  &  & \cellcolor{black!25} \\ \hline
Rewards: 1 &  &  &  \\ \hline
Rewards: 2 &  & \cellcolor{black!25} &  \\ \hline
Rewards: 3 & \cellcolor{black!25} &  &  \\ \hline
Rewards: 4 &  &  & \cellcolor{black!25} \\ \hline
Annual Fee: High &  &  & \cellcolor{black!25} \\ \hline
Annual Fee: Medium & \cellcolor{black!25} &  &  \\ \hline
Annual Fee: Low &  & \cellcolor{black!25} &  \\ \hline
Bank: Bank A & \cellcolor{black!25} & \cellcolor{black!25} & \cellcolor{black!25} \\ \hline
Bank: Bank B &  &  &  \\ \hline
Bank: Out of state bank &  &  &  \\ \hline
Rebate: Low &  & \cellcolor{black!25} &  \\ \hline
Rebate: Medium & \cellcolor{black!25} &  &  \\ \hline
Rebate: High &  &  & \cellcolor{black!25} \\ \hline
Credit Line: Low &  & \cellcolor{black!25} &  \\ \hline
Credit Line: High & \cellcolor{black!25} &  & \cellcolor{black!25} \\ \hline
Grace Period: Short & \cellcolor{black!25} &  & \cellcolor{black!25} \\ \hline
Grace Period: Long &  & \cellcolor{black!25} &  \\ \hline
\end{tabular}
\caption{Outside options for \bank dataset problem instances. \label{table:bank_nopurchase} }
\end{table}

\begin{table}[!ht]
\centering
\scriptsize
\begin{tabular}{|l|c|c|c|} \hline
Attribute & Outside  & Outside  &  Outside  \\
& Option 1 & Option 2 & Option 3 \\ \hline
Age: 36 & \cellcolor{black!25} &  &  \\ \hline
Age: 45 &  &  &  \\ \hline
Age: 52 &  & \cellcolor{black!25} &  \\ \hline
Age: 60 &  &  &  \\ \hline
Age: 68 &  &  & \cellcolor{black!25} \\ \hline
Age: 75 &  &  &  \\ \hline
Military Service: Did not serve & \cellcolor{black!25} &  & \cellcolor{black!25} \\ \hline
Military Service: Served &  & \cellcolor{black!25} &  \\ \hline
Religion: None & \cellcolor{black!25} & \cellcolor{black!25} &  \\ \hline
Religion: Jewish &  &  &  \\ \hline
Religion: Catholic &  &  & \cellcolor{black!25} \\ \hline
Religion: Mainline protestant &  &  &  \\ \hline
Religion: Evangelical protestant &  &  &  \\ \hline
Religion Mormon &  &  &  \\ \hline
College: No BA &  &  &  \\ \hline
College: Baptist college &  &  &  \\ \hline
College: Community college & \cellcolor{black!25} &  &  \\ \hline
College: State university &  &  &  \\ \hline
College: Small college &  & \cellcolor{black!25} &  \\ \hline
College: Ivy League university &  &  & \cellcolor{black!25} \\ \hline
Income: 32K & \cellcolor{black!25} &  &  \\ \hline
Income: 54K &  &  &  \\ \hline
Income: 65K  &  & \cellcolor{black!25} &  \\ \hline
Income: 92K &  &  &  \\ \hline
Income: 210K &  &  & \cellcolor{black!25} \\ \hline
Income 5.1M &  &  &  \\ \hline
Profession: Business owner & \cellcolor{black!25} &  &  \\ \hline
Profession: Lawyer &  &  & \cellcolor{black!25} \\ \hline
Profession: Doctor &  &  &  \\ \hline
Profession: High school teacher &  & \cellcolor{black!25} &  \\ \hline
Profession: Farmer &  &  &  \\ \hline
Profession: Car dealer &  &  &  \\ \hline
Race/Ethnicity: White & \cellcolor{black!25} &  &  \\ \hline
Race/Ethnicity: Native American &  &  &  \\ \hline
Race/Ethnicity: Black &  &  &  \\ \hline
Race/Ethnicity: Hispanic &  & \cellcolor{black!25} &  \\ \hline
Race/Ethnicity: Caucasian &  &  &  \\ \hline
Race/Ethnicity: Asian American &  &  & \cellcolor{black!25} \\ \hline
Gender: Male & \cellcolor{black!25} & \cellcolor{black!25} &  \\ \hline
Gender: Female &  &  & \cellcolor{black!25} \\ \hline
\end{tabular}
\caption{Outside options for \candidate dataset problem instances. \label{table:candidate_nopurchase}}
\end{table}

\begin{table}[!ht]
\centering
\scriptsize
\begin{tabular}{|l|c|c|c|} \hline
Attribute & Outside  & Outside  &  Outside  \\
& Option 1 & Option 2 & Option 3 \\ \hline
Education: No formal &  &  &  \\ \hline
Education: 4th grade &  &  &  \\ \hline
Education: 8th grade &  &  &  \\ \hline
Education: High school &  & \cellcolor{black!25} &  \\ \hline
Education: Two-year college &  &  &  \\ \hline
Education: College degree & \cellcolor{black!25} &  &  \\ \hline
Education: Graduate degree &  &  & \cellcolor{black!25} \\ \hline
Gender: Female & \cellcolor{black!25} & \cellcolor{black!25} &  \\ \hline
Gender: Male &  &  & \cellcolor{black!25} \\ \hline
Origin: Germany & \cellcolor{black!25} &  &  \\ \hline
Origin: France &  &  &  \\ \hline
Origin: Mexico &  &  &  \\ \hline
Origin: Philippines &  &  &  \\ \hline
Origin: Poland &  & \cellcolor{black!25} &  \\ \hline
Origin: India &  &  &  \\ \hline
Origin: China &  &  & \cellcolor{black!25} \\ \hline
Origin: Sudan &  &  &  \\ \hline
Origin: Somalia &  &  &  \\ \hline
Origin: Iraq &  &  &  \\ \hline
Application Reason: Reunite with family &  &  &  \\ \hline
Application Reason: Seek better job & \cellcolor{black!25} & \cellcolor{black!25} & \cellcolor{black!25} \\ \hline
Application Reason: Escape persecution &  &  &  \\ \hline
Profession: Janitor &  &  &  \\ \hline
Profession: Waiter &  &  &  \\ \hline
Profession: Child care provider &  & \cellcolor{black!25} &  \\ \hline
Profession: Gardener &  &  &  \\ \hline
Profession: Financial analyst &  &  &  \\ \hline
Profession: Construction worker &  &  &  \\ \hline
Profession: Teacher &  &  &  \\ \hline
Profession: Computer programmer &  &  &  \\ \hline
Profession: Nurse & \cellcolor{black!25} &  &  \\ \hline
Profession: Research scientist &  &  &  \\ \hline
Profession: Doctor &  &  & \cellcolor{black!25} \\ \hline
Job Experience: None &  &  & \cellcolor{black!25} \\ \hline
Job Experience: 1-2 years &  &  &  \\ \hline
Job Experience: 3-5 years & \cellcolor{black!25} & \cellcolor{black!25} &  \\ \hline
Job Experience: 5+ years &  &  &  \\ \hline
Job Plans: Contract with employer & \cellcolor{black!25} & \cellcolor{black!25} &  \\ \hline
Job Plans: Interviews with employer &  &  &  \\ \hline
Job Plans: Will look for work &  &  & \cellcolor{black!25} \\ \hline
Job Plans: No plans to look for work &  &  &  \\ \hline
Prior Trips to U.S.: Never &  &  & \cellcolor{black!25} \\ \hline
Prior Trips to U.S.: Once as tourist & \cellcolor{black!25} & \cellcolor{black!25} &  \\ \hline
Prior Trips to U.S.: Many times as tourist &  &  &  \\ \hline
Prior Trips to U.S.: Six months with family &  &  &  \\ \hline
Prior Trips to U.S.: Once without authorization &  &  &  \\ \hline
Language: Fluent English & \cellcolor{black!25} & \cellcolor{black!25} &  \\ \hline
Language: Broken English &  &  & \cellcolor{black!25} \\ \hline
Language: Tried English but unable &  &  &  \\ \hline
Language: Used interpreter &  &  &  \\ \hline
\end{tabular}
\caption{Outside options for \immigrant dataset problem instances. \label{table:immigrant_nopurchase} }
\end{table}

\begin{table}[!ht]
\centering
\scriptsize
\begin{tabular}{|l|c|c|c|} \hline
Attribute & Outside  & Outside  &  Outside  \\
& Option 1 & Option 2 & Option 3 \\ \hline
Price: \$70 & \cellcolor{black!25} &  &  \\ \hline
Price: \$75 &  &  &  \\ \hline
Price: \$80 &  &  &  \\ \hline
Price: \$85 &  & \cellcolor{black!25} &  \\ \hline
Price: \$90 &  &  &  \\ \hline
Price: \$95 &  &  &  \\ \hline
Price: \$100 &  &  & \cellcolor{black!25} \\ \hline
Size: Large &  &  & \cellcolor{black!25} \\ \hline
Color: Red &  &  & \cellcolor{black!25} \\ \hline
Logo: Yes &  &  & \cellcolor{black!25} \\ \hline
Handle: Yes &  & \cellcolor{black!25} & \cellcolor{black!25} \\ \hline
PDA Holder: Yes &  & \cellcolor{black!25} & \cellcolor{black!25} \\ \hline
Cellphone Holder: Yes &  & \cellcolor{black!25} & \cellcolor{black!25} \\ \hline
Mesh Pocket: Yes &  & \cellcolor{black!25} & \cellcolor{black!25} \\ \hline
Velcro Flap: Yes &  & \cellcolor{black!25} & \cellcolor{black!25} \\ \hline
Protective Boot: Yes &  & \cellcolor{black!25} & \cellcolor{black!25} \\ \hline
\end{tabular}
\caption{Outside options for \toubia dataset problem instances. (For ease of comparison, only one level of each binary attribute is shown.) \label{table:toubia_nopurchase} }
\end{table}

\clearpage

\subsection{Tuning of $K$ for LC-MNL models}
\label{appendix:numerics_tuning_K}

In this section, we seek to understand what is the best value for the number of customer classes $K$ in each of our four real data sets (\bank, \candidate, \immigrant and \toubia) in terms of model fit. As observed in our experiments using these data sets in Section~\ref{subsec:numerical_experiments_real}, as well as the results of our experiments using our synthetic problem instances in Section~\ref{subsec:numerical_experiments_synthetic}, a general characteristic of formulation~\modelP is that it requires more time to solve as $K$ increases. In this section, we will see that in all four data sets, the value of $K$ that fits the data best is generally small. 

There are many approaches that one can take to tuning the number of classes $K$. The approach we will consider is based on minimizing a model selection measure, specifically the Akaike information criterion (AIC), the Bayesian information criterion (BIC) and the consistent AIC (CAIC). The use of these model selection measures for tuning the number of segments $K$ is standard in the marketing literature (see for example \citealt{kamakura1989probabilistic}, \citealt{bucklin1992brand}, \citealt{andrews2002empirical}, \citealt{andrews2003comparison}; see also Section 14.4 of \citealt{train2009discrete}).

The AIC, BIC and CAIC are defined as 
\begin{align}
\text{AIC} & = -2 \cdot LL + 2 \cdot (K + nK), \\
\text{BIC} & = -2 \cdot LL + \log M \cdot (K + nK), \\
\text{CAIC} & = -2 \cdot LL + (\log M + 1) \cdot (K + nK).
\end{align}
In the above equations, $LL$ is the log likelihood of the LC-MNL model on the complete data set. The term $K + nK$ represents the number of parameters in the LC-MNL model; recall that there are $K$ values of $\lambda_1,\dots, \lambda_K$, and there are $nK$ partworth parameters $\{ \beta_{k,i}\}$. Lastly, the term $M$ is the number of choice tasks in the data set. The AIC, BIC and CAIC metrics can be viewed as the negative log likelihood penalized by a measure of the size of the model, where higher values correspond to models with a better size-adjusted fit to the data. Given a set of candidate values of $K$, the best value of $K$ is the one that minimizes the chosen metric (either AIC, BIC or CAIC). 

\begin{table}
\centering
\begin{tabular}{lrrrr} \toprule
$K$ & \bank & \candidate & \immigrant & \toubia \\ \midrule
 5 & 14649.4 & 2011.4 & 7626.0 & 5748.8 \\ 
  10 & 14366.0 & \bfseries 1992.1 & \bfseries 7365.4 & 5615.2 \\ 
  15 & 14288.8 & 2264.7 & 7383.3 & 5570.2 \\ 
  20 & \bfseries 14260.6 & 2586.7 & 7656.4 &  5545.8 \\ 
  30 & 14262.6 & 3222.6 & 8438.3 & \bfseries 5535.0 \\ 
  40 & 14384.5 & 3881.1 & 9228.4 & 5557.8 \\ 
  50 & 14529.5 & 4535.2 & 10056.0 & 5608.8 \\ \midrule
  $K^*$ & 20 & 10 & 10 & 30 \\ \bottomrule
  \end{tabular}
\caption{Results for tuning of $K$ using AIC. \label{table:real_fit_AIC} }
\end{table}

\begin{table}
\centering
\begin{tabular}{lrrrr} \toprule
$K$ & \bank & \candidate & \immigrant & \toubia \\ \midrule
 5 & \bfseries 15219.6 & \bfseries 2911.9 & \bfseries 9064.6 & \bfseries 6107.0 \\ 
  10 & 15506.3 & 3793.1 & 10242.8 & 6331.7 \\ 
  15 & 15999.4 & 4966.2 & 11699.3 & 6644.9 \\ 
  20 & 16541.3 & 6188.7 & 13411.1 & 6978.8 \\ 
  30 & 17683.6 & 8625.6 & 17070.3 & 7684.5 \\ 
  40 & 18945.9 & 11085.2 & 20737.8 & 8423.9 \\ 
  50 & 20231.2 & 13540.3 & 24442.7 & 9191.4 \\ \midrule
  $K^*$ & 5 & 5 & 5 & 5 \\ \bottomrule
 \end{tabular}
\caption{Results for tuning of $K$ using BIC. \label{table:real_fit_BIC} }
\end{table}

\begin{table}
\centering
\begin{tabular}{lrrrr} \toprule
$K$ & \bank & \candidate & \immigrant & \toubia \\ \midrule
 5 & \bfseries 15294.6 & \bfseries 3076.9 & \bfseries 9274.6 & \bfseries 6162.0 \\ 
  10 & 15656.3 & 4123.1 & 10662.8 & 6441.7 \\ 
  15 & 16224.4 & 5461.2 & 12329.3 & 6809.9 \\ 
  20 & 16841.3 & 6848.7 & 14251.1 & 7198.8 \\ 
  30 & 18133.6 & 9615.6 & 18330.3 & 8014.5 \\ 
  40 & 19545.9 & 12405.2 & 22417.8 & 8863.9 \\ 
  50 & 20981.2 & 15190.3 & 26542.7 & 9741.4 \\  \midrule
  $K^*$ & 5 & 5 & 5 & 5 \\ \bottomrule
 \end{tabular}
\caption{Results for tuning of $K$ using CAIC. \label{table:real_fit_CAIC} }
\end{table}

We compute the AIC, BIC and CAIC for the LC-MNL models estimated using each data set in its entirety; Tables~\ref{table:real_fit_AIC}, \ref{table:real_fit_BIC} and  \ref{table:real_fit_CAIC} present the results. In each table, the best value for a given data set is indicated in bold, and the bottom-most row summarizes the best value of $K$. For the AIC, we can see that the optimal value for each data set is either $K = 10$, $K = 20$ or $K = 30$, whereas for the BIC and CAIC, the optimal value is $K = 5$ for each data set. From a data fitting perspective, these metrics imply that the right number of customer classes should not be more than 20. Thus, while the runtime of formulation~\modelP does increase quickly as $K$ increases, in practice we expect that one will not have to solve it for extremely large values of $K$ when using LC-MNL models.

We additionally note that our findings here are consistent with empirical studies of latent class models in the marketing literature, where typically the number of segments used for latent-class logit modeling does not exceed $K = 20$. For example, \cite{kamakura1989probabilistic} use $K = 5$ segments and \cite{russell1994understanding} use $K = 8$ segments. \cite{andrews2002empirical} compared mixture logit models with continuous mixture distributions and latent class logit models that assume a discrete mixture distribution using a comprehensive simulation study, where purchase data is simulated according to a ground truth model and different mixture logit models are estimated using this data. They showed that even when the ground truth model is a continuous mixture logit model, the latent-class logit model with under ten segments performs better than a richer continuous mixture model. Lastly, Sawtooth Software, an industry leader in conjoint analysis software, provides a latent class MNL estimation tool that by default tests between $K = 1$ and $K = 5$ segments, and allows analysts to go up to $K = 30$ segments but not higher \citep{sawtooth2021latentclass}.  From this perspective, we do not expect that formulation~\modelP will need to be solved for LC-MNL models with $K > 30$.

\subsection{Distribution of partworth magnitudes for real data instances}
\label{appendix:partworth_distributions}

\begin{figure}
\begin{subfigure}{0.5\textwidth}
                \centering
                \includegraphics[width=\textwidth]{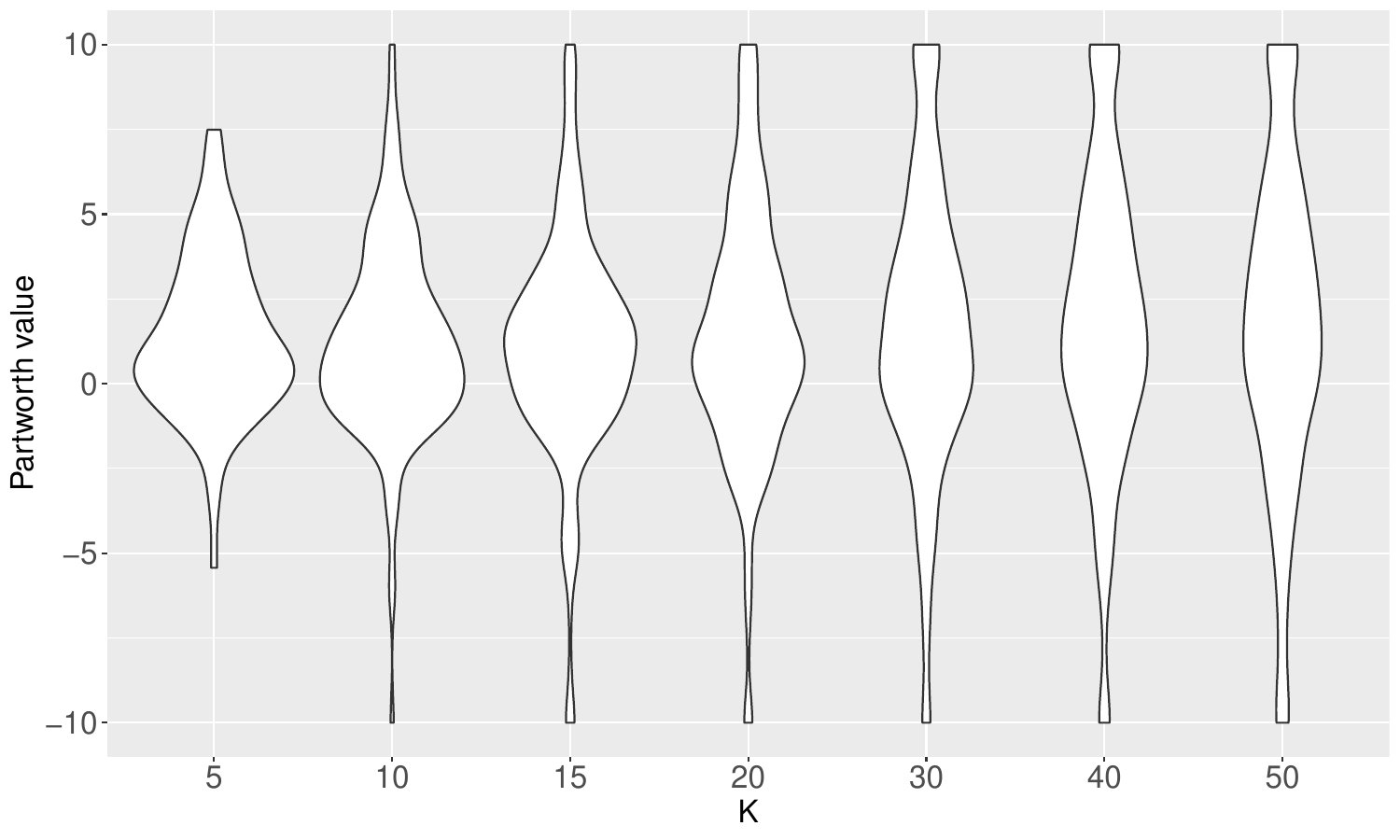}
                \caption{Partworth distributions for LC-MNL models for \bank data set. }
        \end{subfigure} \qquad  \qquad
\begin{subfigure}{0.5\textwidth}
                \centering
                \includegraphics[width=\textwidth]{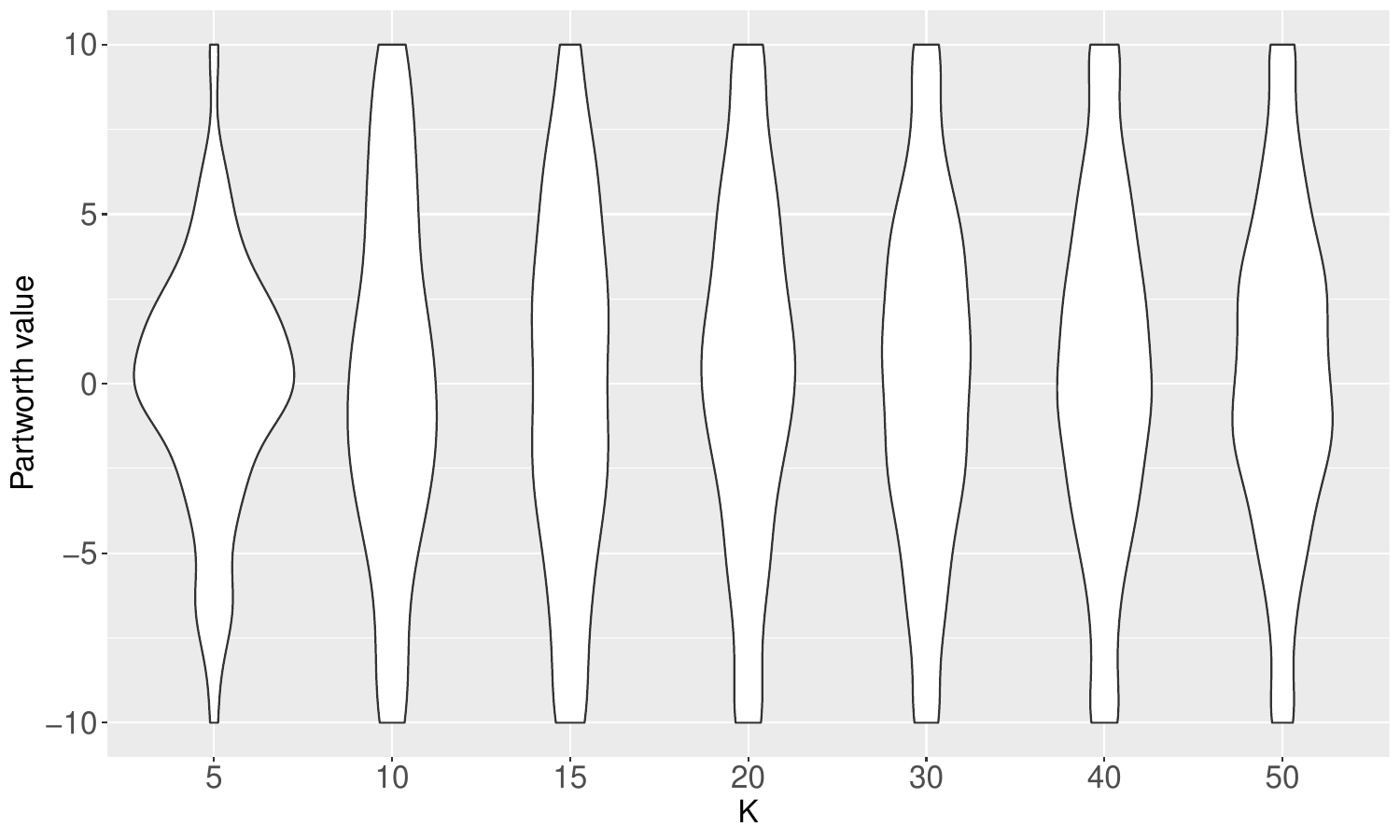}
                \caption{Partworth distributions for LC-MNL models for \candidate data set. }
        \end{subfigure}
\par\bigskip

\begin{subfigure}{0.5\textwidth}
                \centering
                \includegraphics[width=\textwidth]{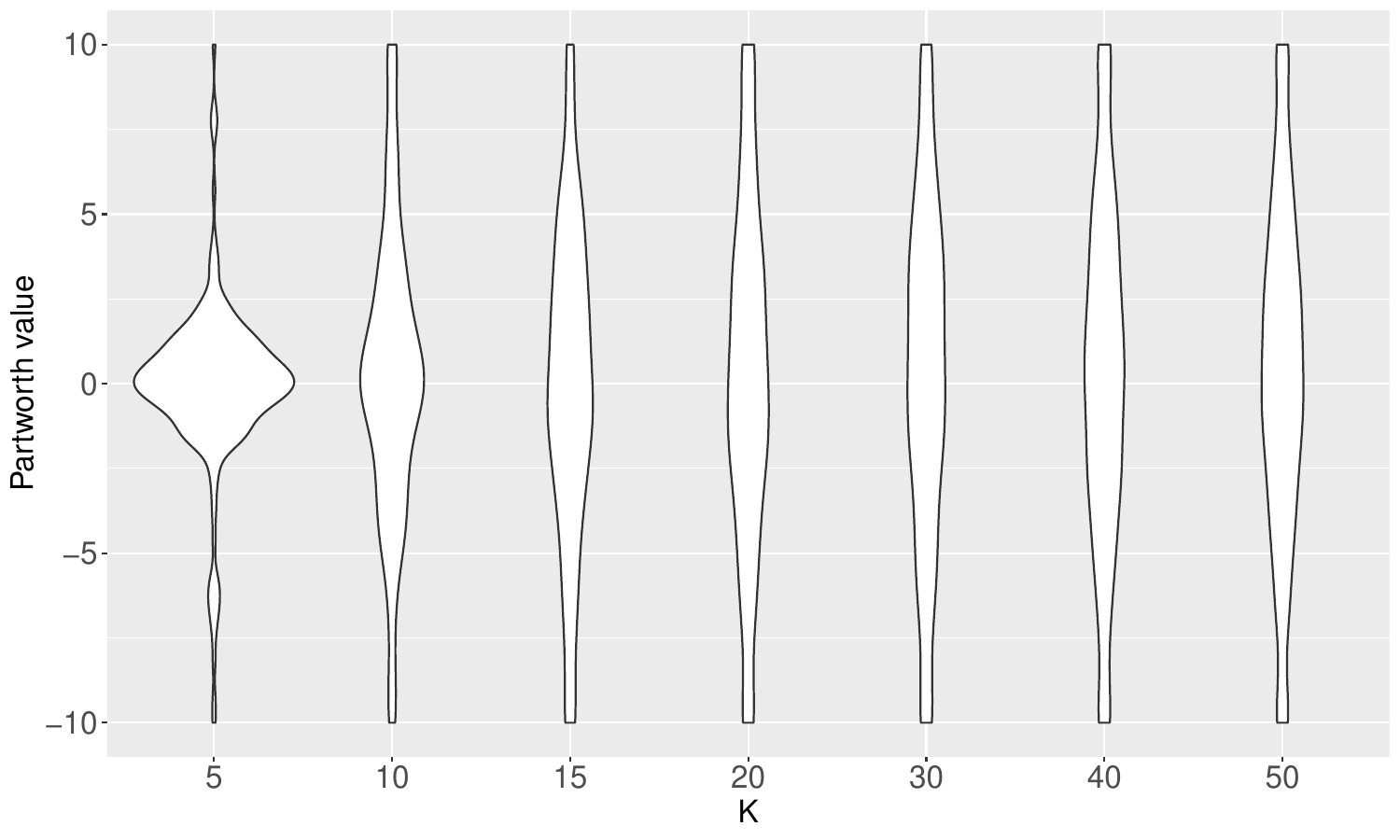}
                \caption{Partworth distributions for LC-MNL models for \immigrant data set. }
        \end{subfigure} \qquad \qquad 
\begin{subfigure}{0.5\textwidth}
                \centering
                \includegraphics[width=\textwidth]{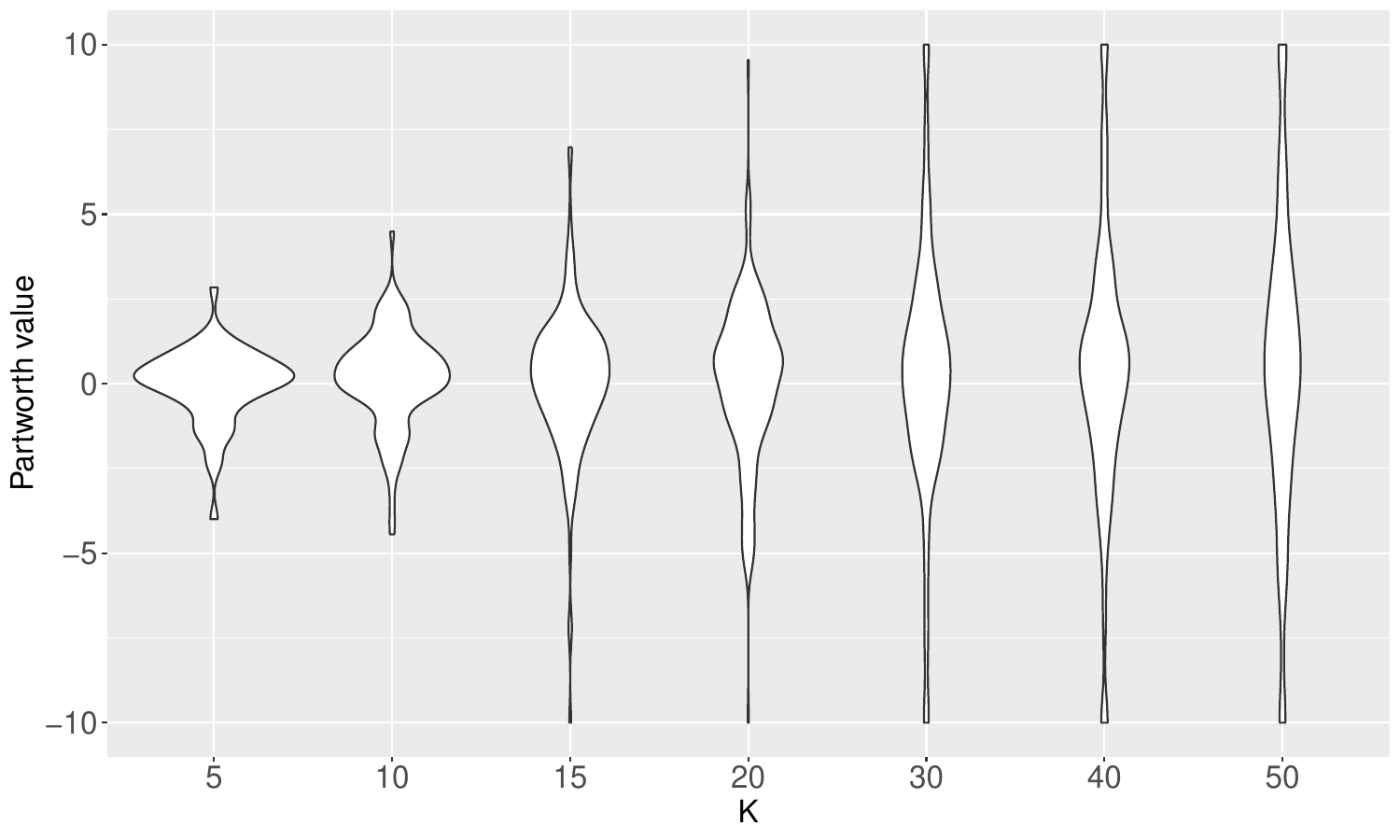}
                \caption{Partworth distributions for LC-MNL models for \toubia data set. }
        \end{subfigure}
\par \bigskip
\caption{Plots of partworth distributions for LC-MNL models estimated from real data sets in Section~\ref{subsec:numerical_experiments_real}. \label{figure:LC_partworth_distributions}}
\end{figure}

\begin{figure}
\begin{center}
\begin{subfigure}[t]{0.2\textwidth}
                \centering
                \includegraphics[width=\textwidth]{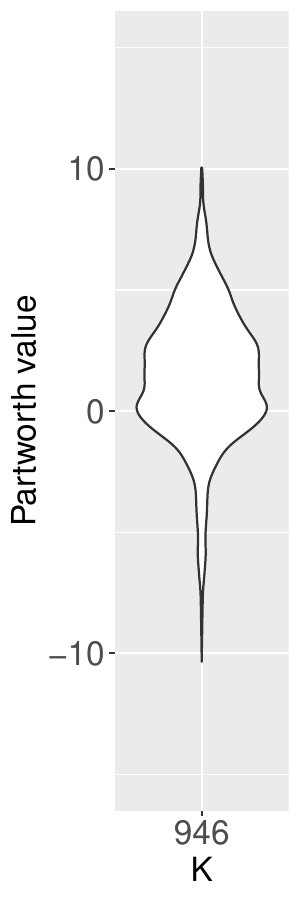}
                \caption{Partworth distributions for HB model for \bank data set. }
        \end{subfigure} \qquad
\begin{subfigure}[t]{0.2\textwidth}
                \centering
                \includegraphics[width=\textwidth]{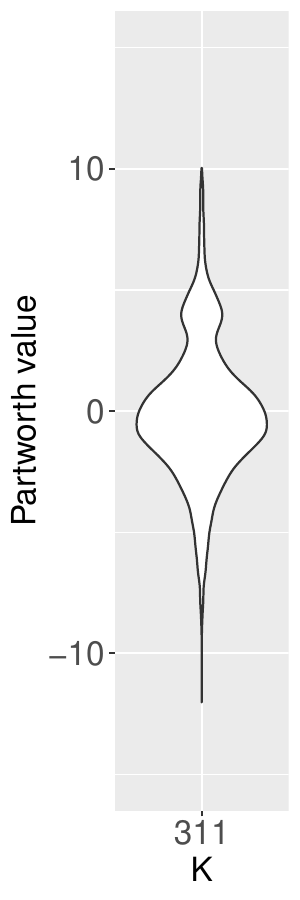}
                \caption{Partworth distributions for HB model for \candidate data set. }
        \end{subfigure} \qquad
\begin{subfigure}[t]{0.2\textwidth}
                \centering
                \includegraphics[width=\textwidth]{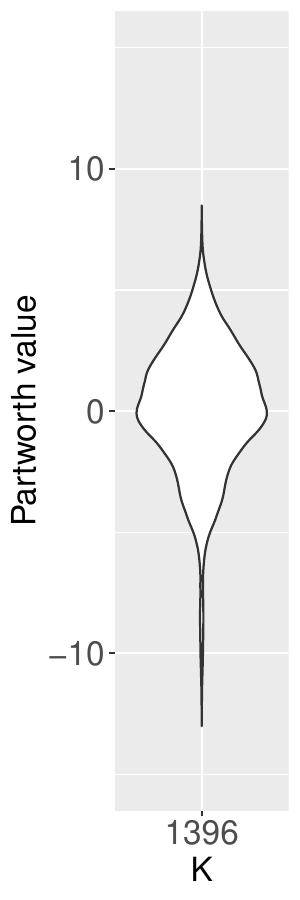}
                \caption{Partworth distributions for HB model for \immigrant data set. }
        \end{subfigure} \qquad
\begin{subfigure}[t]{0.2\textwidth}
                \centering
                \includegraphics[width=\textwidth]{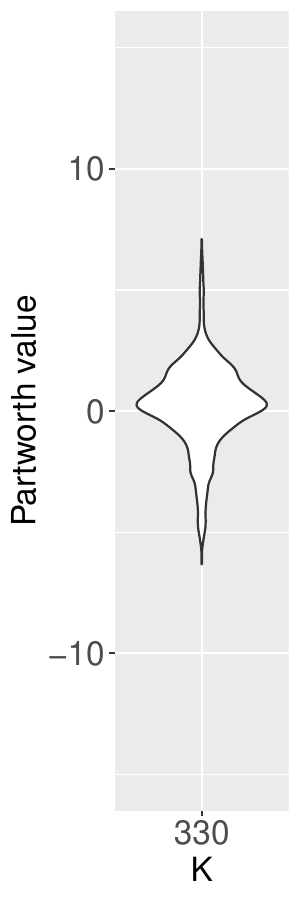}
                \caption{Partworth distributions for HB model for \toubia data set. }
        \end{subfigure}
\end{center}
\caption{Plots of partworth distributions for HB models estimated from real data sets in Section~\ref{subsec:numerical_experiments_real}. \label{figure:HB_partworth_distributions}}
\end{figure}

Figure~\ref{figure:LC_partworth_distributions} plots the distributions (processed via kernel density smoothing) of the $\beta_{k,i}$ values for the LC-MNL models estimated using the EM algorithm for the data sets in Section~\ref{subsec:numerical_experiments_real}. Figure~\ref{figure:HB_partworth_distributions} plots analogous distributions of the $\beta_{k,i}$ values for the HB models estimated using MCMC in Section~\ref{subsec:numerical_experiments_real}. Note that some of the distributions in Figure~\ref{figure:LC_partworth_distributions} are truncated at -10 and +10; as we note in Section~\ref{subsec:numerical_experiments_real}, in our implementation of EM we constrained each partworth parameter to be between -10 and +10 to avoid numerical issues.

\clearpage

\section{Heuristic based on geometric mean maximization}
\label{sec:extensions_geometric_mean}

In this section, we present a heuristic for approximately solving the logit-based SOCPD problem that is based on approximating the objective using a weighted geometric mean. Section~\ref{subsec:extensions_geometric_mean_main} presents the details of the method, while Sections~\ref{proof:theorem_GM_performance_guarantee} and \ref{proof:theorem_GM_NP_Hard} provide proofs of two theoretical results presented in Section~\ref{subsec:extensions_geometric_mean_main}. 

\subsection{Geometric mean heuristic}
\label{subsec:extensions_geometric_mean_main}

To motivate this approach, recall that the logit-based SOCPD problem~\eqref{prob:SOCPD_abstract} involves maximizing the fraction of customers who purchase a product. This objective function is formulated as the weighted sum of the logit probabilities of each customer purchasing the product. An alternate way of understanding this objective is that it is the weighted arithmetic mean of the logit probabilities of the customer types. 

Consequently, instead of formulating the objective of our product design problem as an arithmetic mean, we can instead consider formulating the problem using the geometric mean. This leads to the following optimization problem:
\begin{equation}
\underset{\ab \in \Acal}{\text{maximize}} \prod_{k=1}^K \left[ \frac{ \exp(u_k(\ab))}{1 + \exp(u_k(\ab))} \right]^{\lambda_k}. \label{prob:GM_abstract}
\end{equation}
In other words, rather than trying to optimize the weighted arithmetic mean of the purchase probabilities, this problem seeks to optimize the weighted geometric mean of the purchase probabilities, where the weights indicate the relative proportion of each customer type in the population. 

This formulation is interesting to consider because it provides a lower bound on the optimal value of problem~\eqref{prob:SOCPD_abstract}; the following simple result, which is based on the arithmetic-geometric mean inequality and is stated without proof, formalizes this.
\begin{proposition}
Let $Z^*_{AM}$ and $Z^*_{GM}$ be the optimal objective values of problems~\eqref{prob:SOCPD_abstract} and \eqref{prob:GM_abstract}, respectively. Then $Z^*_{AM} \geq Z^*_{GM}$. 
\end{proposition}
Thus, by solving problem~\eqref{prob:GM_abstract}, we obtain a lower bound on problem~\eqref{prob:SOCPD_abstract}; by evaluating the objective value of the optimal solution of \eqref{prob:GM_abstract} within problem~\eqref{prob:SOCPD_abstract}, we obtain an even stronger lower bound. The solution of the geometric mean problem~\eqref{prob:GM_abstract} can be used as an approximate solution of the arithmetic mean problem~\eqref{prob:SOCPD_abstract}. 

We can further analyze the approximation quality of the solution of problem~\eqref{prob:GM_abstract} with regard to the original problem. Let us use $\xb = (x_1,\dots, x_K)$ to denote the vector of purchase probabilities for the $K$ different customer types, and let us use $\xb(\ab)$ to denote the vector of purchase probabilities for a given product $\ab \in \Acal$:
\begin{equation*}
\xb(\ab) = (x_1(\ab), \dots, x_K(\ab)) = \left( \frac{ \exp(u_1(\ab))}{1 + \exp(u_1(\ab))}, \dots, \frac{ \exp(u_K(\ab))}{1 + \exp(u_K(\ab))} \right).
\end{equation*}
Let us also use $\Xcal$ be the set of achievable customer choice probabilities, given by
\begin{equation}
\Xcal = \left\{ \xb \in [0,1]^K \mid x_k = \frac{ \exp(u_k(\ab))}{1 + \exp(u_k(\ab))} \ \text{for some} \ \ab \in \Acal \right\}.
\end{equation}
Given a vector of choice probabilities $\xb$, we use the function $f: \Xcal \to \mathbb{R}$ to denote the weighted arithmetic mean of $\xb$, with the weights $\lambdab = (\lambda_1,\dots, \lambda_K)$:
\begin{equation}
f(\xb) = \sum_{k=1}^K \lambda_k x_k.
\end{equation}
Similarly, we use $g: \Xcal \to \mathbb{R}$ to denote the weighted geometric mean of $\xb$:
\begin{equation}
g(\xb) = \prod_{k=1}^K x_k^{\lambda_k}.
\end{equation}
Thus, in terms of these two functions, the original logit-based SOCPD problem can be written as $\max_{\ab \in \Acal} f( \xb(\ab) )$, while the geometric mean problem~\eqref{prob:GM_abstract} can be written as $\max_{\ab \in \Acal} g( \xb(\ab) )$. We then have the following guarantee on the performance of any solution of the geometric mean problem~\eqref{prob:GM_abstract} with respect to the objective of the original logit-based SOCPD problem~\eqref{prob:SOCPD_abstract}. 

\begin{theorem}
Let $L$ and $U$ be nonnegative numbers satisfying $L \leq x_k(\ab) \leq U$ for all $k \in \{1,\dots, K\}$ and $\ab \in \Acal$. Let $\ab^* \in \arg \max_{\ab \in \Acal} f(\xb(\ab))$ be a solution of the arithmetic mean problem, and $\hat{\ab} \in \arg \max_{\ab \in \Acal} g(\xb(\ab))$ be a solution of the geometric mean problem. Then the geometric mean solution $\hat{\ab}$ satisfies
\begin{equation*}
f( \xb( \hat{\ab})) \geq \frac{1}{ \sum_{k=1}^K \lambda_k \left( \frac{U}{L} \right)^{1 - \lambda_k} } \cdot f( \xb( \ab^* )).
\end{equation*}
\label{theorem:GM_performance_guarantee}
\end{theorem}

\begin{figure}
\centering
\includegraphics[width=0.7\textwidth]{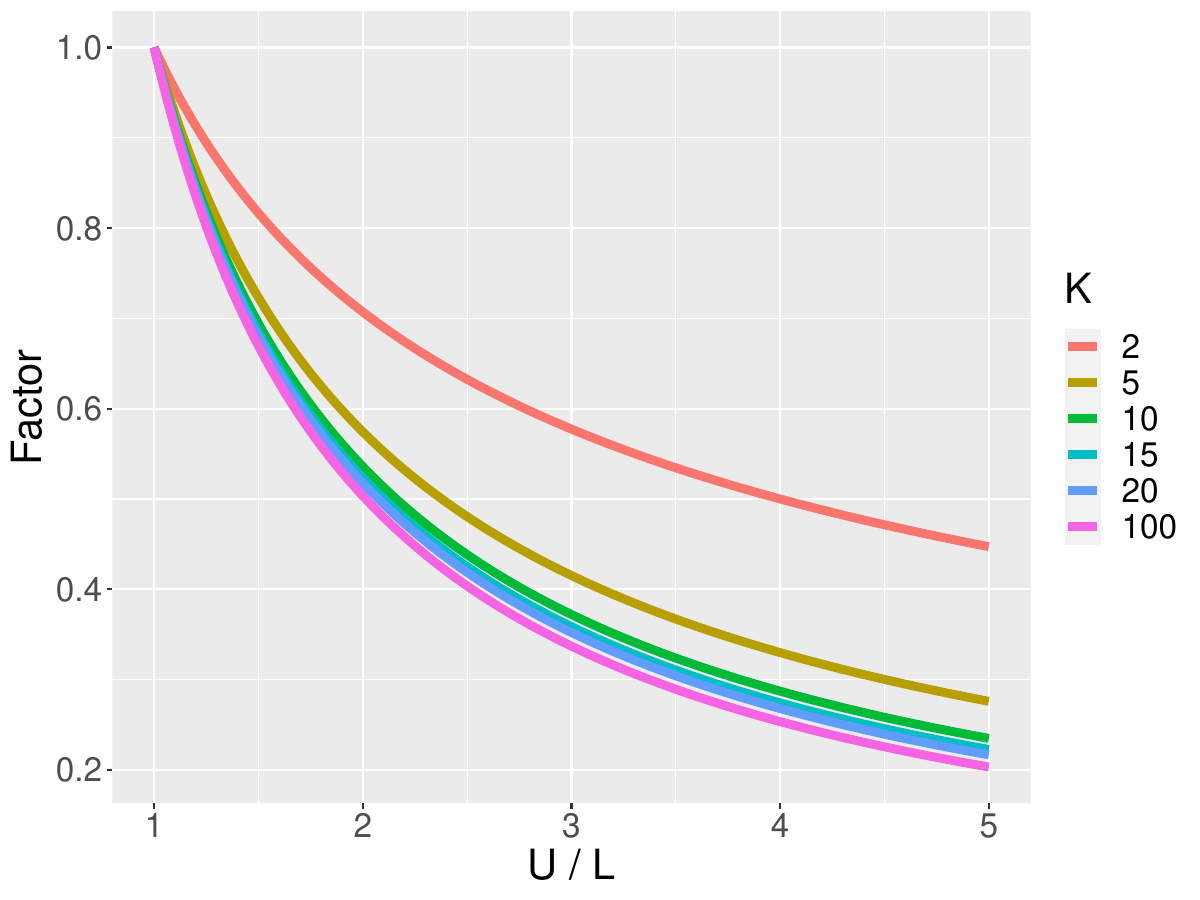}
\caption{Plot of the approximation factor $\Gamma$ as a function of the ratio $U/L$, for different values of $K$. Note that $\lambdab$ is assumed to be the uniform distribution, i.e., $\lambdab = (1/K, \dots, 1/K)$. \label{figure:GM_bound_factor} }
\end{figure}

The proof of Theorem~\ref{theorem:GM_performance_guarantee} (see Section~\ref{proof:theorem_GM_performance_guarantee} of the ecompanion) follows by finding constants $\underline{\alpha}$ and $\overline{\alpha}$ such that $\underline{\alpha} f(\xb) \leq g(\xb) \leq \overline{\alpha} g(\xb)$ for any vector of probabilities $\xb$, and then showing that a solution $\hat{\ab}$ that maximizes $g(\xb(\cdot))$ must be within a factor $\underline{\alpha} / \overline{\alpha}$ of the optimal objective of the arithmetic mean problem. Theorem~\ref{theorem:GM_performance_guarantee} is valuable because it provides some intuition for when a solution $\hat{\ab}$ obtained by solving the geometric mean problem~\eqref{prob:GM_abstract} will be close in performance to the optimal solution of the original (arithmetic mean) problem~\eqref{prob:SOCPD_abstract}. In particular, the factor $\Gamma$ defined as  
\begin{equation*}
\Gamma = \underline{\alpha} \, / \, \overline{\alpha} = \frac{1}{ \sum_{k=1}^K \lambda_k \left( \frac{U}{L} \right)^{1 - \lambda_k} }
\end{equation*}
is decreasing in the ratio $U/L$. Recall that $U$ is an upper bound on the highest purchase probability that can be achieved for any customer type, while $L$ is similarly a lower bound on the lowest purchase probability that can be achieved for any customer type. When the ratio $U/L$ is large, it implies that there is a large range of choice probabilities spanned by the set of product designs $\Acal$. On the other hand, when $U/L$ is small, then the range of choice probabilities is smaller. Thus, the smaller the range of choice probabilities spanned by the set $\Acal$ is small, the closer we should expect the geometric mean solution to be in performance to the optimal solution of the arithmetic mean problem. Figure~\ref{figure:GM_bound_factor} visualizes the dependence of the factor $\Gamma$ on $U / L$ when $\lambdab$ is assumed to be the discrete uniform distribution and $K$ is varied. 

In addition to the ratio $U/L$, the factor $\Gamma$ is also affected by $\lambdab$. It can be verified that the factor $\Gamma$ is minimized when the customer type distribution is uniform, i.e., $\lambdab = (1/K, \dots, 1/K)$. In addition, it can also be verified that when $\lambdab$ is such that $\lambda_k = 1$ for a single customer type (and $\lambda_{k'} = 0$ for all others), the factor $\Gamma$ becomes 1. Thus, the more ``unbalanced'' the customer type distribution $\lambdab$ is, the closer the geometric mean solution should be in performance to the optimal solution of the arithmetic mean problem.

Lastly, with regard to the bounds $U$ and $L$, we note that these can be found easily. In particular, for each customer type $k$, one can compute 
$u_{k,\max} = \max_{\ab \in \Acal} u_k(\ab)$ and $u_{k,\min} = \min_{\ab \in \Acal} u_k(\ab)$, which are the highest and lowest utilities that one can attain for customer type $k$; for many common choices of $\Acal$ this should be an easy problem. (For example, if $\Acal$ is simply $\{0,1\}^n$, we can find $u_{k,\max}$ by setting to 1 those attributes for which $\beta_{k,i} > 0$ and setting to 0 all other attributes; $u_{k,\min}$ can be found in a similar manner). One can then compute $L$ and $U$ as 
\begin{align*}
U & = \max_{k = 1,\dots, K} \frac{ \exp( u_{k,\max} )}{1 + \exp( u_{k,\max}) }, \\
L & = \min_{k = 1,\dots, K} \frac{ \exp( u_{k,\min} )}{1 + \exp( u_{k,\min}) }.
\end{align*}

We now turn our attention to how one can solve problem~\eqref{prob:GM_abstract}. While problem~\eqref{prob:GM_abstract} is still a challenging nonconvex problem, it is possible to transform it into a mixed-integer convex problem. To do so, we consider taking the logarithm of the objective function of \eqref{prob:GM_abstract}:
\begin{align*}
\log \prod_{k=1}^K \left[ \frac{ \exp(u_k(\ab))}{1 + \exp(u_k(\ab))} \right]^{\lambda_k} & = \sum_{k=1}^K \lambda_k \log \left( \frac{ \exp(u_k(\ab))}{1 + \exp(u_k(\ab))}  \right) \\
& = \sum_{k=1}^K \lambda_k \cdot \left( u_k(\ab) - \log(1 + \exp(u_k(\ab))) \right).
\end{align*}
This transformation is useful because the logarithm function is monotonic, so any solution that maximizes the logarithm of the objective function maximizes the objective function itself. This leads to the following mixed-integer convex program:
\begin{subequations}
\begin{alignat}{2}
& \underset{\ab, \ub}{ \text{maximize} } & \quad & \sum_{k=1}^K \lambda_k \cdot \left( u_k - \log(1 + \exp(u_k)) \right) \\
& \text{subject to} & & u_k = \beta_{k,0} + \sum_{i=1}^n \beta_{k,i} a_i, \quad \forall\ k \in \{1,\dots, K\}, \\
& & & \Cb \ab \leq \db, \\
& & & a_i \in \{0,1\}, \quad \forall i \in \{1,\dots,n\}. 
\end{alignat}%
\label{prob:GM_MICONVP}%
\end{subequations}
This formulation can be further reformulated as a mixed-integer exponential cone program, and solved using Mosek. In terms of complexity, we note that problem~\eqref{prob:GM_abstract} is still a hard problem, which is formalized in the proposition below.
\begin{theorem}
The geometric mean problem~\eqref{prob:GM_abstract} is NP-Hard. \label{theorem:GM_NP_Hard}
\end{theorem}
We refer readers to Section~\ref{proof:theorem_GM_NP_Hard} for the proof of this result. However, in spite of this result, our experience is that the conic reformulation of \eqref{prob:GM_MICONVP} can generally be solved quite quickly, and much faster than our exact formulations (\modelRA, \modelP and \modelPRPT); in our experiments in Section~\ref{subsec:numerical_experiments_synthetic}, we find that synthetic problem instances of up to $n = 70$ attributes and $K = 30$ customer types can be solved to full optimality in no more than two minutes on average. Despite this positive empirical result, we again note that the geometric mean approach only provides a heuristic and not an exact approach to solving the original problem~\eqref{prob:SOCPD_abstract}, and indeed, in Section~\ref{subsec:numerical_experiments_real}, we shall see that the \GM heuristic can be quite suboptimal in real data instances.

\subsection{Proof of Theorem~\ref{theorem:GM_performance_guarantee}}
\label{proof:theorem_GM_performance_guarantee}

Let $\xb^* = \xb( \ab^* )$ and $\hat{\xb} = \xb( \hat{\ab})$. To prove the result we proceed in three steps. \\

\textbf{Step 1:} The first step in our proof is to show that if there exist nonnegative constants $\overline{\alpha}$ and $\underline{\alpha}$ such that $g$ satisfies
\begin{equation}
\underline{\alpha} f(\xb) \leq g(\xb) \leq \overline{\alpha} f(\xb) \label{eq:uniform_constant_bound}
\end{equation}
for all $\xb \in \Xcal$, then $\hat{\xb}$ satisfies
\begin{equation}
f(\hat{\xb}) \geq (\underline{\alpha} / \overline{\alpha} ) \cdot f(\xb^*).
\end{equation}

To establish this, we will first bound the quantity $f(\xb^*) - f(\hat{\xb})$. We have
\begin{align*}
f(\xb^*) - f(\hat{\xb}) & = [f(\xb^*) - g(\xb^*)] + [ g(\xb^*) - g(\hat{\xb})] + [ g(\hat{\xb}) - f(\hat{\xb})] \\
& \leq f(\xb^*) - g(\xb^*) + g(\hat{\xb}) - f(\hat{\xb}) \\
& \leq f(\xb^*) - \underline{\alpha} f(\xb^*) + \overline{\alpha} f(\hat{\xb}) - f(\hat{\xb}) \\
& = (1 - \underline{\alpha}) f(\xb^*) - (1 - \overline{\alpha}) f(\hat{\xb}) \\
& = (1 - \overline{\alpha} + \overline{\alpha} - \underline{\alpha}) f(\xb^*) - (1 - \overline{\alpha}) f(\hat{\xb}) \\
& = (1 - \overline{\alpha})(f(\xb^*) - f(\hat{\xb})) + (\overline{\alpha} - \underline{\alpha}) f(\xb^*)
\end{align*}
where the first step follows by algebra; the second step follows since $g(\xb^*) \leq g(\hat{\xb})$, which is true by the definition of $\hat{\xb}$ as the vector of choice probabilities for an optimal product $\hat{\ab}$ for the function $g( \xb(\ab))$; the third step follows by \eqref{eq:uniform_constant_bound}; and the remaining steps by algebra.

Observe that by re-arranging the inequality
\begin{equation}
f(\xb^*) - f(\hat{\xb}) \leq (1 - \overline{\alpha})(f(\xb^*) - f(\hat{\xb})) + (\overline{\alpha} - \underline{\alpha}) f(\xb^*)
\end{equation}
we obtain that
\begin{equation}
\overline{\alpha} [ f(\xb^*) - f(\hat{\xb}) ] \leq (\overline{\alpha} - \underline{\alpha}) f(\xb^*).
\end{equation}
Since $\overline{\alpha}$ is nonnegative, dividing through by $\overline{\alpha}$ we obtain 
\begin{equation}
 f(\xb^*) - f(\hat{\xb}) \leq \frac{ (\overline{\alpha} - \underline{\alpha})}{\overline{\alpha}} f(\xb^*),
\end{equation}
and re-arranging, we obtain
\begin{align*}
f(\hat{\xb}) & \geq \left[ 1 - \frac{\overline{\alpha} - \underline{\alpha}}{\overline{\alpha}} \right] f(\xb^*) \\
& = (\underline{\alpha}  / \overline{\alpha}) \cdot f(\xb^*),
\end{align*}
which is the desired result. \\

\textbf{Step 2:} We now establish explicit values for the constants $\overline{\alpha}$ and $\underline{\alpha}$. Recall that by the arithmetic-geometric mean inequality, $g(\xb) \leq f(\xb)$ for all $\xb \in \Xcal$. Therefore, a valid choice of $\overline{\alpha}$ is 1. 

For $\underline{\alpha}$, we proceed as follows. Consider the ratio $f(\xb) / g(\xb)$. For any $\xb$, we have
{\allowdisplaybreaks
\begin{align*}
\frac{ f(\xb) }{ g(\xb)} & = \frac{ \sum_{k=1}^K \lambda_k x_k }{ \prod_{k=1}^K x_k^{\lambda_k}} \\
& = \sum_{k=1}^K \lambda_k \cdot x_k^{1 - \lambda_k} \cdot \prod_{k' \neq k} x_{k'}^{-\lambda_{k'}} \\ 
& \leq \sum_{k=1}^K \lambda_k \cdot U^{1 - \lambda_k} \cdot \prod_{k' \neq k} L^{-\lambda_{k'}} \\
& = \sum_{k=1}^K \lambda_k \cdot U^{1 - \lambda_k} \cdot L^{- \sum_{k' \neq k} \lambda_{k'}} \\
& = \sum_{k=1}^K \lambda_k \cdot U^{1 - \lambda_k} \cdot L^{\lambda_k - 1} \\
& = \sum_{k=1}^K \lambda_k  \left( \frac{U}{L} \right)^{1 - \lambda_k},
\end{align*}}
where the first step follows by the definitions of $f$ and $g$; the second by algebra; the third by the fact that the function $h(x) = x^{1 - \lambda_k}$ is increasing in $x$ (since $1 - \lambda_k \geq 0$), and that the function $\bar{h}(x) = x^{-\lambda_{k'}}$ is decreasing in $x$ (since $-\lambda_k \leq 0$); the fourth by algebra; the fifth by recognizing that $\sum_{k' = 1}^K \lambda_k = 1$, which implies that $\lambda_k - 1 = - \sum_{k' \neq k} \lambda_{k'}$; and the last by algebra. This implies that a valid choice of $\underline{\alpha}$ is 
\begin{equation}
\underline{\alpha} = \frac{1}{\sum_{k=1}^K \lambda_k  \left( \frac{U}{L} \right)^{1 - \lambda_k}}.
\end{equation}

\textbf{Step 3:} We conclude the proof by combining Steps 1 and 2. In particular, by using $\overline{\alpha} = 1$ and $\underline{\alpha} = [ \sum_{k=1}^K \lambda_k  \left( U / L \right)^{1 - \lambda_k} ]^{-1}$, we obtain that 
\begin{equation*}
f( \xb( \hat{\ab})) \geq \frac{1}{ \sum_{k=1}^K \lambda_k \left( \frac{U}{L} \right)^{1 - \lambda_k} } \cdot f( \xb( \ab^* )),
\end{equation*}
as required. \Halmos

\subsection{Proof of Theorem~\ref{theorem:GM_NP_Hard}}
\label{proof:theorem_GM_NP_Hard}

To prove this result, we will show that MAX-3SAT problem can be reduced to the geometric mean problem~\eqref{prob:GM_abstract}. 

Given an instance of the MAX-3SAT problem, we construct an instance of the geometric mean problem~\eqref{prob:GM_abstract} as follows. Let the number of attributes $n$ be equal to the number of binary variables in the MAX-3SAT instance, and we define $\Acal$ as $\{0,1\}^n$. Each attribute of our product will correspond to one of the binary variables. Let each customer type $k$ correspond to one of the $K$ clauses, and we set $\lambda_k = 1 / K$. We define the parameters $p_L$ and $p_U$ as
\begin{align}
p_L & = \frac{1}{100K}, \\
p_U & = {p_L}^{p_L} = \left( \frac{1}{100K}\right) ^{\frac{1}{100K}}
\end{align}
and we define the utilities $Q_L$ and $Q_U$ as
\begin{align}
Q_L &= \log \left(\frac{p_L}{1 - p_L} \right), \\
Q_U &= \log \left(\frac{p_U}{1 - p_U} \right).   
\end{align}
We define the partworth parameters as how we did in the proof of Theorem~\ref{proof:theorem_APXHard}. For each customer type $k$, let $J_k \in \{0,1,2,3\}$ denote the number of negative literals in the corresponding clause $k$ of the MAX-3SAT instance (i.e., how many literals of the form $\neg x_i$ appear in $c_k$). We define the partworths $\beta_{k,1}, \dots, \beta_{k,n}$ of customer type $k$ as follows:
\begin{align}
\beta_{k,i} & = \begin{cases}
0 & \text{if variable $x_i$ does not appear in any literal of clause $k$}, \\
Q_U - Q_L & \text{if the literal $x_i$ appears in clause $k$}, \\
Q_L - Q_U & \text{if the literal $\neg x_i$ appears in clause $k$}, \\
\end{cases}
\end{align}
for each $i \in \{1,\dots, n\}$, and we define the constant part of the utility $\beta_{k,0}$ as 
\begin{align}
\beta_{k,0} & = Q_L + J_k \cdot (Q_U - Q_L). 
\end{align}

Next, we need to show that, given an optimal solution $\ab$ to the geometric mean problem, the solution  $\xb$, which is obtained by setting $x_i = a_i$ for each $i \in \{1,\dots, n\}$, is an optimal solution of the MAX-3SAT problem. However, before we establish this, we make the following observation. Since $\log(p_L)$ is a constant and $\sum_{k = 1}^K \lambda_k = 1$, we can subtract $\sum_{k = 1}^K \lambda_k \log(p_L)$ from the objective function of the geometric mean problem and divide it by $-\log(p_L) > 0$ without changing the optimal solution of the problem. After this transformation, we obtain the following objective function:
\begin{equation}
\sum_{k = 1}^K \lambda_k \frac{-\log(p_L) + \log\left(\frac{\exp(u_k(\ab))}{1 + \exp(u_k(\ab))} \right) }{-\log(p_L)}.
\end{equation}
It is straightforward to see that maximizing the geometric mean objective $\sum_{k=1}^K \lambda_k (u_k(\ab) - \log(1 + e^{u_k(\ab)}))$  is equivalent to maximizing this modified objective. In the remainder of the proof, we use this objective function for the geometric mean problem. Let 
\begin{equation}
h(u) = \frac{-\log(p_L) + \log\left( \frac{\exp(u)}{1 + \exp(u)}\right) }{-\log(p_L)}.
\end{equation} 
Observe that, if the product attribute $\ab$ is set such that none of the literals in a clause $k$ is satisfied, then $u_k(\ab) = Q_L$ and $h(u_k(\ab)) = 0$. Otherwise, $u_k(\ab) \geq Q_U$ and $h(u_k(\ab)) \geq 1 - 1/(100K)$. Moreover, $h(u) < 1$ for all $u \in \mathbb{R}$. 

Before we proceed with the proof of the theorem, we prove two lemmas. We first define $g_k(\xb)$ as in the proof of Theorem~\ref{proof:theorem_APXHard_MAX3SAT}. That is, $g_k(\xb) = 1$ if clause $k$ in the MAX-3SAT problem is satisfied by solution $\xb$ and $g_k(\xb) = 0$ otherwise. Then we establish the following relations between $g_k(\xb)$ and $h(u_k(\ab))$ for a MAX-3SAT problem solution $\xb$ and a geometric logit-based SOCPD solution $\ab$.

\begin{lemma}
	Let $\xb$ and $\ab$ defined such that $x_i = a_i$ for all $i \in \{1,\dots,n\}$. Then, we have
	\begin{equation}
	g_k(\xb) - \frac{1}{100K} \leq h(u_k(\ab)) \leq g_k(\xb).
	\end{equation}
	\label{lemma:GM_1}
\end{lemma}

\begin{proofvvm} To establish the first inequality, notice that, $h$ is an increasing function since the logarithm and logistic functions are increasing functions and $-\log(p_L)$ is a positive constant. If $g_k(\xb) = 1$, we have $u_k(\ab) \geq Q_U$, which implies that $h(u_k(\ab)) \geq h(Q_U) = 1 - 1/100K = g_k(\xb) - 1/100K$. If $g_k(\xb) = 0$, we have $u_k(\ab) = Q_L$, which implies $h(u_k(\ab)) = 0 > g_k(\xb) - 1/100K$. Therefore, for both cases, $h(u_k(\ab)) \geq g_k(\xb) - 1/100K$ holds.
	
	The second inequality also holds because $1 > h(u_k(\ab)) \geq 1 - 1/100K$ if $g_k(\xb) = 1$, and $h(u_k(\ab)) = 0$ otherwise. \hfill $\square$
\end{proofvvm}

Now, we will establish a relation between the objective functions of the geometric mean logit-based SOCPD and the MAX-3SAT problem.

\begin{lemma}
		Let $\xb$ and $\ab$ defined such that $x_i = a_i$ for all $i \in \{1,\dots,n\}$. Then, we have
	\begin{equation}
	\left\lceil \sum_{k = 1}^K h(u_k(\ab)) \right\rceil = \sum_{k = 1}^K g_k(\xb).
	\end{equation}
 	\label{lemma:GM_2}
\end{lemma}

\begin{proofvvm}
	By Lemma~\ref{lemma:GM_1}, we have
	\begin{equation}
	\sum_{k = 1}^K \left( g_k(\xb) - \frac{1}{100K}\right) = \sum_{k = 1}^K g_k(\xb) - \frac{1}{100} \leq \sum_{k = 1}^K h(u_k(\ab)) \leq \sum_{k = 1}^K g_k(\xb).
	\end{equation}
	Since $\sum_{k = 1}^K g_k(\xb)$ is an integer, this implies that 
		$\left\lceil \sum_{k = 1}^K h(u_k(\ab)) \right\rceil = \sum_{k = 1}^K g_k(\xb)$. \hfill $\square$
\end{proofvvm}

Finally, we will conclude the proof of the theorem by showing that, given an optimal solution $\ab$ to the geometric logit-based SOCPD problem, the solution $\xb$, which we obtain by setting $x_i = a_i$ for all $i \in \{1,\dots,n\}$, is an optimal solution to the MAX-3SAT problem. To verify this, we use the notation that we defined in the proof of Theorem~\ref{proof:theorem_APXHard} and follow a similar proof technique. Suppose that, $\xb$ is not the optimal solution to the MAX-3SAT problem, and there exists a solution $\tilde{\xb}$, which achieves a higher number of satisfied clauses than $\xb$. Let $\tilde{\ab}$ be the solution we obtain by setting $\tilde{a}_i = \tilde{x}_i$ for all $i \in \{1,\dots,n\}$. Then, we have

\begin{align}
	\left\lceil \sum_{k = 1}^K h(u_k(\tilde{\ab})) \right\rceil &= \sum_{k = 1}^K g_k(\tilde{\xb}) \\
	&> \sum_{k = 1}^K g_k(\xb) \\
	&= \left\lceil \sum_{k = 1}^K h(u_k(\ab)) \right\rceil
\end{align}
where the equalities follow from Lemma~\ref{lemma:GM_2} and the inequality follows from the assumption that $\sum_{k = 1}^K g_k({\tilde{\xb}}) > \sum_{k = 1}^K g_k(\xb)$. Since $\sum_{k = 1}^K g_k({\tilde{\xb}})$ and $\sum_{k = 1}^K g_k(\xb)$ are integers, this implies that $\sum_{k = 1}^K h(u_k(\tilde{\ab}))  > \sum_{k = 1}^K h(u_k(\ab))$, which contradicts the optimality of $\ab$. Therefore, $\xb$ must be the optimal solution to the MAX-3SAT problem. \Halmos

\clearpage
\pagebreak 
\section{Robust logit-based share-of-choice product design}
\label{sec:robust}

A key assumption in the logit-based SOCPD problem is that the underlying parameters that determine customer choice -- the distribution $\lambdab$ and the partworth vectors $\betab_1, \dots \betab_K$ -- are known precisely. In practice, these parameters are estimated from data (as in our numerical experiments with real data in Section~\ref{subsec:numerical_experiments_real}) and there may be errors in these estimated values; thus, these parameters are subject to uncertainty. This is important because a product that is optimized based on a single $\lambdab$ and a single collection of partworth vectors $\betab_1,\dots, \betab_K$ may yield significantly lower market share if the actual $\lambdab$ and $\betab_1,\dots, \betab_K$ values are different from the ones used in the optimization.

In this section, we consider two different robust optimization approaches to the logit-based SOCPD problem that address uncertainty in the partworth vectors $\betab_1,\dots, \betab_K$. In particular, let $\betab = (\betab_1,\dots, \betab_K) \in \Rbb^{(n+1)K}$ denote the concatenation of the partworth vectors of all $K$ customer types; we refer to $\betab$ as the \emph{grand partworth vector}. The robust logit-based SOCPD problem can then be written as the following max-min problem:
\begin{equation}
\max_{\ab \in \Acal} \min_{\betab \in \Ucal} \sum_{k=1}^K \lambda_k \cdot \sigma( \beta_{k,0} + \sum_{i=1}^n \beta_{k,i} a_i ),
\end{equation}
where $\Ucal \subseteq \Rbb^{(n+1)K}$ is an uncertainty set of possible partworth vectors. In this problem, we seek to find the product design vector $\ab$ that maximizes the worst-case share-of-choice, where the worst-case is taken over all grand partworth vectors in $\Ucal$. 

The rest of this section is organized as follows. Section~\ref{subsec:robust_product_model} presents our first approach, which assumes that $\Ucal$ is structured as a Cartesian product of smaller uncertainty sets corresponding to each customer type. Section~\ref{subsec:robust_joint_model} presents our second approach, which assumes that $\Ucal$ is structured as a budget uncertainty set. Section~\ref{subsec:robust_product_results} presents a small set of computational experiments for the first approach, while Section~\ref{subsec:robust_joint_results} presents computational results for the second approach. 

Before we continue, we note that the two approaches developed below only consider uncertainty in the grand partworth vector $\betab$ and not in the probability distribution $\lambdab$. We focus on this form of uncertainty as we believe this is the more interesting case to consider. When there is only uncertainty in $\lambdab$, the robust logit-based SOCPD problem can be written as 
\begin{equation}
\max_{\ab \in \Acal} \min_{\lambdab \in \Lambda} \sum_{k=1}^K \lambda_k \cdot \sigma( \beta_{k,0} + \sum_{i=1}^n \beta_{k,i} a_i ),
\end{equation}
where $\Lambda$ is an uncertainty set of probability mass functions supported on $[K]$. This problem can be analyzed in a straightforward fashion as the objective function is linear in $\lambdab$, and therefore the inner worst-case problem that minimizes over $\lambdab$ is a linear program. Depending on the structure of $\Lambda$, one can potentially reformulate the inner problem to eliminate the minimization over $\lambdab$ (for example, if $\Lambda$ is a polyhedron, then one can use LP duality to reformulate the problem with a finite number of additional variables and constraints) and reformulate each $\sigma(\cdot)$ term using one of the three formulations presented earlier (\modelRA, \modelP or \modelPRPT). Alternatively, one can also design a cutting plane procedure that replaces $\Lambda$ with a finite set $\hat{\Lambda}$, solves the corresponding restricted master problem, and then identifies a new $\lambdab$ to add to $\hat{\Lambda}$ by solving the worst-case problem $\min_{\lambdab \in \Lambda} \sum_{k=1}^K \lambda_k \cdot \sigma( \beta_{k,0} + \sum_{i=1}^n \beta_{k,i} a_i )$.

\subsection{Robust approach 1: product uncertainty set}
\label{subsec:robust_product_model}

In the first approach that we consider, we assume that the uncertainty set $\Ucal$ of the grand partworth vector is structured as a Cartesian product of type-specific uncertainty sets, that is,
\begin{equation}
\Ucal = \Ucal_1 \times \Ucal_2 \times \dots \times \Ucal_K,
\end{equation}
where $\Ucal_k \subseteq \Rbb^{n+1}$ is an uncertainty set governing the partworth vector $\betab_k$ of customer type $k$. Under this uncertainty set, the robust logit-based SOCPD problem can be written as 
\begin{align}
\max_{\ab \in \Acal} \min_{\betab \in \Ucal} \left\{ \sum_{k=1}^K \lambda_k \sigma( \beta_{k,0} + \sum_{i=1}^n \beta_{k,i} a_i) \right\}
\end{align}

Due to the product form of the uncertainty set, this problem admits a nice reformulation, which we now explain. In particular, we can re-write the problem as 
\begin{align}
& \max_{\ab \in \Acal} \min_{\betab \in \Ucal} \left\{ \sum_{k=1}^K \lambda_k \sigma( \beta_{k,0} + \sum_{i=1}^n \beta_{k,i} a_i) \right\} \\
& = \max_{\ab \in \Acal} \min_{\betab_1 \in \Ucal_1, \dots, \betab_K \in \Ucal_K} \left\{ \sum_{k=1}^K \lambda_k \sigma( \beta_{k,0} + \sum_{i=1}^n \beta_{k,i} a_i) \right\} \\
& = \max_{\ab \in \Acal}  \left\{ \sum_{k=1}^K \lambda_k \min_{\betab_k \in \Ucal_k} \sigma( \beta_{k,0} + \sum_{i=1}^n \beta_{k,i} a_i) \right\} ,
\end{align}
where the last step follows because, due to the product form of the uncertainty set, the minimization over the overall grand partworth vector $\betab$ decomposes into $K$ minimizations over each individual customer type's partworth vector $\betab_k$. 

From here, the problem can be further reformulated by observing that the logistic response function $\sigma(\cdot)$ is monotonic, and so the minimization over $\betab_k$ can be pushed inside of $\sigma(\cdot)$:
\begin{align}
& \max_{\ab \in \Acal} \left\{ \sum_{k=1}^K \lambda_k \min_{\betab_k \in \Ucal_k} \sigma( \beta_{k,0} + \sum_{i=1}^n \beta_{k,i} a_i) \right\}, \\
& = \max_{\ab \in \Acal} \left\{ \sum_{k=1}^K \lambda_k \sigma( \min_{\betab_k \in \Ucal_k}  \{ \beta_{k,0} + \sum_{i=1}^n \beta_{k,i} a_i \} ) \right\}.
\end{align}
Recall now from our formulation \modelP that the decision variable $w_k$ represents the linearization of $x_{k,1} \cdot u_k$. To model the inner minimization, we replace the constraint that defines $w_k$ in that formulation, which is
\begin{equation}
w_k = \beta_{k,0} x_{k,1} + \sum_{i=1}^n \beta_{k,i} y_{k,i},
\end{equation}
with the following robust constraint:
\begin{equation}
w_k \leq \beta_{k,0} x_{k,1} + \sum_{i=1}^n \beta_{k,i} y_{k,i}, \quad \forall \ \betab_k \in \Ucal_k.
\end{equation}
Formulation~\modelP thus becomes the following formulation, which we denote by \modelPRobustS:
\begin{subequations}
\begin{alignat}{2}
\modelPRobustS: \quad & \underset{\ab, \ub, \wb, \xb, \yb}{\text{maximize}} & \quad & \sum_{k=1}^K \lambda_k x_{k,1} \\
& \text{subject to} & & x_{k,1} + x_{k,1} e^{-w_k / x_{k,1} } \leq 1, \quad \forall \ k \in [K], \\
& & & w_k \leq \beta_{k,0} x_{k,1} + \sum_{i=1}^n \beta_{k,i} y_{k,i}, \quad \forall \ k \in [K], \ \betab_k \in \Ucal_k, \label{prob:PRobustS_robust}\\
& & & y_{k,i} \leq a_i, \quad \forall \ k \in [K], \ i \in [n], \\
& & & y_{k,i} \leq x_{k,1}, \quad \forall \ k \in [K], \ i \in [n], \\
& & & y_{k,i} \geq x_{k,1} + a_i - 1, \quad \forall \ k \in [K], \ i \in [n], \\
& & & y_{k,i} \geq 0, \quad \forall \ k \in [K], \ i \in [n], \\
& & & x_{k,0}, x_{k,1} \geq 0, \quad \forall \ k \in [K], \\
& & & \Cb \ab \leq \db, \\
& & & \ab \in \{0,1\}^n. 
\end{alignat}
\label{prob:PRobustS}
\end{subequations}
The key distinction between \modelPRobustS and \modelP is that constraint~\eqref{prob:PRobustS_robust} is quantified over all partworth vectors $\betab_k \in \Ucal_k$, and thus \eqref{prob:PRobustS_robust} describes a potentially uncountably infinite collection of linear inequalities. As is standard in robust optimization, if each $\Ucal_k$ admits a tractable representation, then one can re-write the constraint as 
\begin{equation}
w_k \leq \min_{\betab_k \in \Ucal_k} \{ \beta_{k,0} x_{k,1} + \sum_{i=1}^n \beta_{k,i} y_{k,i} \} \label{eq:robust_w_constraint}
\end{equation}
and reformulate the minimization problem on the right hand side of the inequality to obtain an equivalent but finite representation. For example, if $\Ucal_k$ is a polyhedron, the minimization problem is a linear program, and one can use LP duality theory to reformulate the constraint exactly using a finite number of constraints and variables. Alternatively, one can consider solving the problem using constraint generation, where one replaces $\Ucal_k$ with a finite subset $\hat{\Ucal}_k$, and solves the minimization problem on the right-hand side to identify partworth vectors at which the constraint is violated.  %

In the experiments that we will present in Section~\ref{subsec:robust_product_results}, we will assume that each $\Ucal_k$ is a continuous budget uncertainty set \citep[see][for more details]{bertsimas2004price} defined as
\begin{equation}
\Ucal_k = \{ \betab_k = \bar{\betab}_k - \hat{\betab}_k \circ \xib_k \ \mid \ \zerob \leq \xib_k \leq \oneb, \ \oneb^\top \xib_k \leq \Gamma \},
\end{equation}
where the vectors $\oneb$ and $\zerob$ are used to denote $(n+1)$-dimensional vectors of all ones and zeros, respectively, and $\circ$ denotes the component-wise product of two vectors. In this definition, the vector $\bar{\betab}_k \in \Rbb^{n+1}$ is the vector of nominal partworths and the vector $\hat{\betab}_k$ is the vector of maximum allowable deviations, where each value $\hat{\beta}_{k,i}$ represents the most that the partworth $\beta_{k,i}$ may deviate from its nominal value $\beta_{k,i}$. The value $\xi_{k,i}$ is bounded between 0 and 1 and represents the fraction of the maximum deviation of $\hat{\beta}_{k,i}$; the constraint $\sum_{i=0}^n \xi_{k,i} \leq \Gamma$ models that we only allow up to $\Gamma$ partworth values to maximally deviate from their nominal values. Note that in our uncertainty set, we only consider downward deviations, which results in the form $\bar{\betab}_k - \hat{\betab}_k \circ \xib_k$. Although budget uncertainty sets (as for example in \citealt{bertsimas2004price}) usually allow for both upward and downward deviations, in our case it is not necessary to consider upward deviations, because such deviations are never optimal for the inner minimization problem and will in general never be a part of the worst-case solution. 

For this budget uncertainty set, constraint~\eqref{eq:robust_w_constraint} can be reformulated by applying LP duality to the right-hand side of \eqref{eq:robust_w_constraint}, and introducing a new set of decision variables and constraints. In particular, it can be shown that constraint~\eqref{eq:robust_w_constraint} is equivalent to
\begin{align}
w_k & \leq \bar{\beta}_{k,0} x_{k,1} + \sum_{i=1}^n \bar{\beta}_{k,i} y_{k,i} - \Gamma q_k - \sum_{i=1}^n \tau_{k,i} \\
q_k + \tau_{k,0} & \geq \hat{\beta}_{k,0} x_{k,1}, \\
q_k + \tau_{k,i} & \geq \hat{\beta}_{k,i} y_{k,i}, \quad \forall \ i \in [n], \\
q_k & \geq 0, \\
\tau_{k,i} & \geq 0, \quad \forall \ i \in [n],
\end{align}
where $q_k$ and $\taub_k = (\tau_{k,0}, \dots, \tau_{k,n})$ are new decision variables that are added to problem~\modelPRobustS. Thus, problem~\modelPRobustS becomes
\begin{subequations}
\begin{alignat}{2}
\modelPRobustSBudget:\quad & \underset{\ab, \qb, \wb, \xb, \yb, \taub}{\text{maximize}} & \quad & \sum_{k=1}^K \lambda_k x_{k,1} \\
& \text{subject to} & & x_{k,1} + x_{k,1} e^{-w_k / x_{k,1} } \leq 1, \quad \forall \ k \in [K], \\
& & & w_k \leq \bar{\beta}_{k,0} x_{k,1} + \sum_{i=1}^n \bar{\beta}_{k,i} y_{k,i} - \Gamma q_k - \sum_{i=1}^n \tau_{k,i}, \quad \forall \ k \in [K], \\
& & & q_k + \tau_{k,0} \geq \hat{\beta}_{k,0} x_{k,1}, \quad \forall \ k \in [K] \\
& & & q_k + \tau_{k,i} \geq \hat{\beta}_{k,i} y_{k,i}, \quad \forall \ k \in [K],\ i \in [n], \\
& & & y_{k,i} \leq a_i, \quad \forall \ k \in [K], \ i \in [n], \\
& & & y_{k,i} \leq x_{k,1}, \quad \forall \ k \in [K], \ i \in [n], \\
& & & y_{k,i} \geq x_{k,1} + a_i - 1, \quad \forall \ k \in [K], \ i \in [n], \\
& & & y_{k,i} \geq 0, \quad \forall \ k \in [K], \ i \in [n], \\
& & & x_{k,0}, x_{k,1} \geq 0, \quad \forall \ k \in [K], \\
& & & \Cb \ab \leq \db, \\
& & & \ab \in \{0,1\}^n, \\
& & & q_k  \geq 0, \quad \forall \ k \in [K], \\
& & & \tau_{k,i} \geq 0, \quad \forall \ k \in [K], \ i \in [n],
\end{alignat}
\label{prob:PRobustS}
\end{subequations}

\subsection{Robust approach 2: joint uncertainty set}
\label{subsec:robust_joint_model}

In our second approach, let our uncertainty set $\Ucal$ of partworth vectors be defined as 
\begin{equation}
\Ucal = \left\{ \betab = \bar{\betab} - \hat{\betab} \circ \Xib \quad \vline \quad \sum_{k=1}^K \sum_{i=0}^n \xi_{k,i} \leq \Gamma, \ \Xib \in \{0,1\}^{(n+1)K} \right\}.
\end{equation}
In the above definition, $\bar{\betab} = (\bar{\betab}_1, \dots, \bar{\betab}_K) \in \Rbb^{K(n+1)}$ is the vector of nominal values of the partworths, where $\bar{\betab}_k = (\bar{\beta}_{k,0}, \bar{\beta}_{k,1}, \dots, \bar{\beta}_{k,n})$ is the vector containing the nominal partworths and the nominal intercept for customer type $k$. The vector $\hat{\betab} = (\hat{\betab}_1, \dots, \hat{\betab}_K) \in \Rbb^{K(n+1)}$ is the vector of maximal allowed deviations of the partworth parameters, where $\hat{\betab}_k = (\hat{\beta}_{k,0}, \dots, \hat{\beta}_{k,n})$ represents the vector of maximal allowed deviations of each parameter (i.e., $\hat{\beta}_{k,i}$ is the most that the partworth $\beta_{k,i}$ is allowed to deviate from its nominal value $\bar{\beta}_{k,i}$). The vector $\Xib$ is defined as $\Xib = (\xib_1, \dots, \xib_K) \in \{0,1\}^{(n+1)K}$, where each $\xib_k = (\xi_{k,0}, \dots, \xi_{k,n}) \in \{0,1\}^{n+1}$ is the vector of binary variables indicating whether partworth $\beta_{k,j}$ is deviating from the nominal value $\bar{\beta}_{k,j}$ ($\xi_{k,j} = 1$) or not ($\xi_{k,j} = 0$). We refer to $\xib_k$ as the perturbation pattern of customer type $k$, and $\Xib$ as the grand perturbation pattern. 

The uncertainty set $\Ucal$ represents the set of all vectors of partworth vectors where at most $\Gamma$ parameters are equal to $\bar{\beta}_{k,i} - \hat{\beta}_{k,i}$ and the rest are equal to their nominal value $\bar{\beta}_{k,i}$. The idea of this uncertainty set is that while each partworth parameter $\beta_{k,j}$ may deviate from its nominal value, we expect that in the worst case, there should not be too many such parameters deviating from their nominal value. Note that unlike the product uncertainty set in the previous section, there is no limit on how many deviations can occur for each customer type, and so an admissible grand partworth vector $\betab$ from $\Ucal$ may be such that all of the deviations occur for a small subset of the $K$ customer types, with the partworths for the other customer types unperturbed. Additionally, $\Ucal$ is a discrete uncertainty set, whereas the uncertainty set of the previous section may be discrete or continuous. The motivation for this choice is tractability; we shall discuss this in more detail shortly. 

The corresponding robust logit-based SOCPD problem is then
\begin{equation}
\max_{\ab \in \Acal} \min_{\betab \in \Ucal} \left\{ \sum_{k=1}^K \lambda_k \sigma( \beta_{k,0} + \sum_{i=1}^n \beta_{k,i} a_i) \right\}, \label{prob:robust_joint_master_abstract}
\end{equation}
where the goal is to find the product design vector $\ab$ that maximizes the worst-case logit-based share-of-choice, where the worst-case is taken over all partworth vectors $\betab$ belonging to $\Ucal$. Note that in this model, we are again assuming that the nominal values of the customer type probabilities $\lambda_1,\dots, \lambda_K$ are not subject to uncertainty. 

Note that this problem is rather difficult to solve, because the inner worst-case problem, 
\begin{equation}
\min_{\betab \in \Ucal} \left\{ \sum_{k=1}^K \lambda_k \sigma( \beta_{k,0} + \sum_{i=1}^n \beta_{k,i} a_i) \right\} \label{prob:inner_worst_case_joint}
\end{equation}
is a binary nonlinear optimization problem, similarly to the nominal logit-based SOCPD problem. Typically in robust optimization, the inner worst-case problem is a tractable optimization problem that can be reformulated using duality. For example, if the objective is linear in the uncertain parameter, and the uncertainty set is polyhedral, then the inner worst-case problem is a linear program, and one can apply linear programming duality to reformulate the inner worst-case problem using a finite number of additional variables and constraints. In our setting, such an approach is not applicable due to the nature of this inner problem. 

Instead, what we can hope to do is to solve the overall robust problem~\eqref{prob:robust_joint_master_abstract} using delayed constraint generation. In this approach, we first reformulate the problem in epigraph form:
\begin{subequations}
\begin{alignat}{2}
& \underset{\ab, \theta}{\text{maximize}} & \quad & \theta \\
& \text{subject to} & & \theta \leq \sum_{k=1}^K \lambda_k \sigma( \beta_{k,0} + \sum_{i=1}^n \beta_{k,i} a_i), \quad \forall \ \betab \in \Ucal, \label{prob:robust_joint_master_epigraph_main_constraint} \\
& & & \ab \in \Acal. 
\end{alignat}
\label{prob:robust_joint_master_epigraph}
\end{subequations}
Now, instead of solving problem~\eqref{prob:robust_joint_master_epigraph} with all possible $\betab$ enumerated, we start with constraint~\eqref{prob:robust_joint_master_epigraph_main_constraint} enforced for only a finite subset $\hat{\Ucal} \subseteq \Ucal$. We then solve this restricted master problem to obtain a solution $\ab$. With this solution $\ab$ in hand, we now solve the following separation problem:
\begin{subequations}
\begin{alignat}{2}
& \underset{\Xib}{\text{minimize}} & \quad & \sum_{k=1}^K \lambda_k \sigma( \bar{u}_k - \hat{\beta}_{k,0} \xi_{k,0} - \sum_{i=1}^n \hat{\beta}_{k,i} a_i \xi_{k,i}) \\
& \text{subject to} & & \sum_{k=1}^K \sum_{i=0}^n \xi_{k,i} \leq \Gamma, \\
& & & \xi_{k,i} \in \{0,1\}, \quad \forall \ k \in [K], \ i \in \{0,1,\dots,n\},
\end{alignat}
\label{prob:robust_joint_separation_abstract}
\end{subequations}
where $\bar{u}_k$ is defined as $\bar{u}_k = \bar{\beta}_{k,0} + \sum_{i=1}^n \bar{\beta}_{k,i} a_i$, which is the utility of $\ab$ using the nominal partworth values for customer type $k$. 

Although this problem is challenging, we can reformulate it as a mixed-integer convex program using the same type of technique as we used to obtain formulation \modelP. In particular, as in formulation \modelP, let $\pi_{k,1}$ and $\pi_{k,0}$ denote the purchase probability of the product and the no-purchase probability, respectively, for customer type $k$; let $u_k$ denote the utility of product $k$; let $h_k$ denote the linearization of $u_k \cdot \pi_{k,0}$; and let $z_{k,i}$ denote the linearization of $\pi_{k,0} \cdot \xi_{k,i}$. With these definitions, this separation problem~\eqref{prob:robust_joint_separation_abstract} can be re-written as the following mixed-integer exponential cone program:
\begin{subequations}
\begin{alignat}{2}
& \underset{ \pib, \ub, \hb, \zb, \Xib}{ \text{minimize}} & \quad & \sum_{k=1}^K \lambda_k x_{k,1} \\
& \text{subject to} & & \pi_{k,0} + \pi_{k,0} \cdot e^{ \frac{h_k}{\pi_{k,0}}} \leq 1, \quad \forall \ k \in [K], \\
& & & \pi_{k,0} + \pi_{k,1} = 1, \quad \forall\ k \in [K], \\
& & & u_k = \bar{u}_k - \sum_{i=0}^n \hat{\beta}_{k,i} \cdot a_i \cdot \xi_{k,i}, \quad \forall \ k \in [K], \\
& & & h_k = \bar{u}_k \pi_{k,0} - \sum_{i=0}^n \hat{\beta}_{k,i} \cdot a_i \cdot z_{k,i}, \quad \forall \ k \in [K], \\
& & & z_{k,i} \leq \pi_{k,0},  \quad \forall \ k \in [K], i \in \{0,1,\dots,n\}, \\
& & & z_{k,i} \leq \xi_{k,i}, \quad \forall \ k \in [K], i \in \{0,1,\dots,n\}, \\
& & & z_{k,i} \geq \pi_{k,0} + \xi_{k,i} - 1,  \quad \forall \ k \in [K], i \in \{0,1,\dots,n\}, \\
& & & \sum_{k=1}^K \sum_{i=0}^n \xi_{k,i} \leq \Gamma, \quad \forall\ k \in [K], \\
& & & \sum_{k=1}^K \sum_{i=0}^n z_{k,i} \leq \Gamma \cdot \pi_{k,0},  \quad \forall \ k \in [K], \\
& & & \sum_{k=1}^K \sum_{i=0}^n (\xi_{k,i} - z_{k,i}) \leq \Gamma \cdot \pi_{k,1},  \quad \forall \ k \in [K], \\ 
& & & \xi_{k,i} \in \{0,1\}, \quad \forall \ k \in [K], i \in \{0,1,\dots,n\}, \\
& & & z_{k,i} \geq 0,  \quad \forall \ k \in [K], i \in \{0,1,\dots,n\}.
\end{alignat}
\label{prob:robust_joint_separation}
\end{subequations}
By solving this problem, we obtain the solution $\Xib$; to obtain the corresponding $\betab$ vector, we simply calculate it as $\betab = \bar{\betab} - \hat{\betab} \circ \Xib$. We add the corresponding constraint to the master problem~\eqref{prob:robust_joint_master_epigraph}, and solve the master problem again. 

We make several important remarks about the restricted master problem. First, note that the restricted master problem can be formulated as a mixed-integer exponential cone program that is very similar to the nominal problem. The main difference is that the variables $x_{k,0}, x_{k,1}, u_k, w_k, y_{k,i}$ are now additionally indexed by $\xib$. In particular, $x_{k,\xib,0}$ and $x_{k,\xib,1}$ are the choice probabilities for the no-purchase option and the product for customer type $k$ when its partworth vector deviates according to the perturbation pattern $\xib$. Similarly, $u_{k, \xib}$ is the utility of the product for customer type $k$ with the perturbation pattern $\xib$; $w_{k, \xib}$ is the linearization of $u_{k,\xib} \cdot x_{k,\xib,1}$; $y_{k,\xib,i}$ is the linearization of $a_i \cdot x_{k,\xib,1}$. Each combination of $k$ and $\xib$ requires analogs of the constraints~\eqref{prob:P_perspective} to \eqref{prob:P_x_definition} of \modelP and in particular, requires one exponential cone.

Second, and related to the previous point, is that across different worst-case realizations of the grand partworth vector $\betab$ that are generated from $\Ucal$, the same deviation pattern $\xib$ could appear multiple times for the same customer type. This implies that the same variable $x_{k,\xib,1}$ that represents the purchase probability for customer type $k$ under perturbation pattern $\xib$ could appear in multiple instances of the epigraph constraint~\eqref{prob:robust_joint_master_epigraph_main_constraint}. This is important because it allows for efficiency in terms of how many exponential cones are used to model the $x_{k,\xib,1}$ variables. A naive implementation of constraint generation would introduce a new exponential cone for each $x$ variable that appears in constraint~\eqref{prob:robust_joint_master_epigraph_main_constraint}, resulting in $K \cdot M$ exponential cones after $M$ worst-case realizations are generated. By being careful about whether a perturbation pattern has been generated previously, one can reuse variables that have already been introduced, and reduce how many new exponential cones get added to the master problem. 

Lastly, we alluded earlier that the choice of a discrete budget uncertainty set, as opposed to a continuous budget uncertainty set is motivated by tractability, and after laying out the overall constraint generation procedure, it should become clear why a discrete uncertainty set may be easier to work with. In particular, the worst-case problem~\eqref{prob:robust_joint_separation_abstract} can be formulated exactly as the mixed-integer exponential cone program~\eqref{prob:robust_joint_separation}. If one were to consider a continuous uncertainty set $\Ucal$, then the worst-case problem~\eqref{prob:inner_worst_case_joint} would be a continuous, non-convex problem, and in such a situation, it is not clear how one can solve such a problem to global optimality (in order to implement a constraint generation/cutting plane method), or how one can otherwise tractably reformulate the overall robust problem.

\subsection{Numerical experiments with product uncertainty robust approach} 
\label{subsec:robust_product_results}

In this section, we present a small set of numerical experiments to demonstrate the value of the robust approach using the product budget uncertainty set described in Section~\ref{subsec:robust_product_model}. We consider the synthetic instances from Section~\ref{subsec:numerical_experiments_synthetic} with $n = 30$, and $K \in \{10, 20\}$, and the scale factor parameter $c$ fixed to $5$. 

We set up the budget uncertainty set $\Ucal_k$ of each customer type $k$ as follows. We use the value of each term $\beta_{k,i}$ as the nominal value $\bar{\beta}_{k,i}$ in our uncertainty set. For the intercept, we assume that there is no uncertainty, and set $\hat{\beta}_{k,0} = 0$. For each attribute, we assume that $\hat{\beta}_{k,i} = c' \cdot | \bar{\beta}_{k,i}|$, where $c' \in \{0.1, 0.2\}$ is a parameter that will be tested. We vary the budget $\Gamma$ in the set $\{1,2,3,4,5,6,7\}$. Note that in general, there are $Kn$ partworth parameters, not counting the intercept; setting the budget as $m$ means that at most $m$ out of $n$ parameters deviate from their nominal values, for a total of $Km$ out of $Kn$ parameters over all $K$ customer types. 

For each $(n,K, c')$ combination, we solve formulation~\modelPRobustSBudget for each of the 20 synthetic instances, and record the objective value. We solve the formulation using Mosek and impose a time limit of one hour. In addition, for each $n$ and $K$, we also compute the worst-case share-of-choice of the nominal product vector by solving $\min_{\betab \in \Ucal} \left\{ \sum_{k=1}^K \lambda_k \sigma( \beta_{k,0} + \sum_{i=1}^n \beta_{k,i} a_i) \right\}$. 

To compare the robust and nominal approaches, we compute two different metrics:
\begin{enumerate}
	\item \emph{Worst-case loss}: The worst-case loss (WCL) is defined as 
	\begin{equation}
	\WCL = \frac{ F(\ab^N, \bar{\betab}) - \min_{\betab \in \Ucal} F(\ab^N, \betab) }{F(\ab^N, \bar{\betab}) } \times 100\%,
	\end{equation}
	where $\ab^N$ is the nominal product design (i.e., $\ab^N \in \arg \max_{\ab \in \Acal} F(\ab, \bar{\betab})$), and $F(\ab, \betab) = \sum_{k=1}^K \lambda_k \sigma( \beta_{k,0} + \sum_{i=1}^n \beta_{k,i} a_i)$ is the share-of-choice of the product vector $\ab$ under the overall partworth vector $\betab = (\betab_1, \dots, \betab_K)$. In words, it is the percentage reduction in the share-of-choice of the nominal solution $\ab^N$, relative to the nominal share-of-choice, under the worst-case realization in $\Ucal$. A high value of $\WCL$ implies that the worst-case performance of the nominal product design deteriorates significantly when the realized partworths differ from their nominal values. 
	\item \emph{Relative improvement}: The relative improvement (RI) is defined as 
	\begin{equation}
	\RI = \frac{ \min_{\betab \in \Ucal} F(\ab^R, \betab) - \min_{\betab \in \Ucal} F(\ab^N, \betab)}{\min_{\betab \in \Ucal} F(\ab^N, \betab)} \times 100\%, \label{eq:RI_definition}
	\end{equation}
	where $\ab^N$ is the nominal product design and $\ab^R$ is the robust product design (i.e., $\ab^R \in \arg \max_{\ab \in \Acal} \min_{\betab \in \Ucal} F(\ab, \betab)$). In words, it is the improvement in worst-case share-of-choice performance of the robust product design $\ab^R$ relative to the nominal product design $\ab^N$. A high value of $\RI$ implies that the robust design delivers better performance under uncertainty than the nominal design.
\end{enumerate}

Table~\ref{table:robust_product} reports the average $\WCL$ and $\RI$ for each of the $(n, K, c')$ combinations, as well as the average computation time of \modelPRobustSBudget. From this table, we can see that the deterioration of the nominal solution when exposed to the worst grand partworth vector $\betab$ from $\Ucal$ can be large (as high 33\% when $K = 20$, $c' = 0.20$ and $\Gamma = 5$), and that the robust solution improves on the nominal solution significantly in terms of worst-case share-of-choice (as much as 28\%). Additionally, the computation times of this approach are reasonable; while they are larger than the nominal formulation (see Table~\ref{table:R1_synthetic_gap_time} in Section~\ref{subsec:numerical_experiments_synthetic}), they are no more than about 20 minutes in the largest case. 

\begin{table} \centering
\begin{tabular}{lrrrrrr} \toprule
$c'$ & $n$ & $K$ & $\Gamma$ & $\WCL$ (\%) & $\RI$ (\%) & Time (s) \\ \midrule
0.10 &  30 &  10 &   1 & 0.96 & 0.00 & 20.90 \\ 
  0.10 &  30 &  10 &   2 & 2.15 & 0.06 & 25.89 \\ 
  0.10 &  30 &  10 &   3 & 3.55 & 0.38 & 30.36 \\ 
  0.10 &  30 &  10 &   4 & 5.14 & 1.13 & 31.10 \\ 
  0.10 &  30 &  10 &   5 & 6.78 & 2.06 & 32.16 \\[0.25em]
  0.10 &  30 &  20 &   1 & 1.99 & 0.15 & 461.28 \\ 
  0.10 &  30 &  20 &   2 & 4.35 & 0.84 & 616.34 \\ 
  0.10 &  30 &  20 &   3 & 6.96 & 1.84 & 755.08 \\ 
  0.10 &  30 &  20 &   4 & 9.70 & 3.42 & 872.96 \\ 
  0.10 &  30 &  20 &   5 & 12.50 & 5.41 & 956.01 \\[0.5em]
  0.20 &  30 &  10 &   1 & 2.28 & 0.07 & 26.43 \\ 
  0.20 &  30 &  10 &   2 & 6.00 & 1.59 & 28.91 \\ 
  0.20 &  30 &  10 &   3 & 11.17 & 5.15 & 36.13 \\ 
  0.20 &  30 &  10 &   4 & 17.61 & 11.66 & 41.62 \\ 
  0.20 &  30 &  10 &   5 & 24.41 & 20.31 & 50.46 \\[0.25em]
  0.20 &  30 &  20 &   1 & 4.55 & 0.91 & 591.00 \\ 
  0.20 &  30 &  20 &   2 & 10.87 & 4.15 & 780.24 \\ 
  0.20 &  30 &  20 &   3 & 18.38 & 10.50 & 901.27 \\ 
  0.20 &  30 &  20 &   4 & 26.26 & 18.59 & 1158.55 \\ 
  0.20 &  30 &  20 &   5 & 33.90 & 27.94 & 1294.17 \\  \bottomrule
\end{tabular}
\caption{Performance of robust solutions using the product uncertainty set approach (Section~\ref{subsec:robust_product_model}) on synthetic data instances. \label{table:robust_product}}
\end{table}

\subsection{Numerical experiments with joint uncertainty robust approach}
\label{subsec:robust_joint_results}

In this second set of numerical experiments, we aim to understand the value of the robust approach using the joint budget uncertainty set described in Section~\ref{subsec:robust_joint_model}. We again consider the synthetic instances from Section~\ref{subsec:numerical_experiments_synthetic} with $n = 30$ and $K \in \{10, 20\}$. We set the scale factor parameter $c$ to  $5$. 

To set up the joint budget uncertainty set $\Ucal$, we proceed as follows. We use the value of each term $\beta_{k,i}$ as the nominal value $\bar{\beta}_{k,i}$ in our uncertainty set. For the intercept, we assume that there is no uncertainty, and set $\hat{\beta}_{k,0} = 0$. For each attribute, we assume that $\hat{\beta}_{k,i} = c' \cdot | \bar{\beta}_{k,i}|$, where $c' \in \{0.1, 0.2\}$ is a parameter that will be tested. Lastly, we vary the budget $\Gamma$ in the set $\{1 \cdot K, 2 \cdot K, \dots, 5 \cdot K\}$. Note that in general, there are $K(n+1)$ utility parameters ($n$ partworths plus one intercept, for each customer type); setting the budget as $m \cdot K$ can be interpreted as anticipating up to $m$ out of $n$ partworths of each segment (on average) to vary from their nominal values.

For each of the 20 instances corresponding to each combination $(n, K, c) \in \{30\} \times \{10, 20\} \times \{5\}$, we solve the robust formulation~\eqref{prob:robust_joint_master_abstract} for each $\Gamma$ and each $c'$. We apply the constraint generation method described in Section~\ref{subsec:robust_joint_model}. Due to the significantly greater computational requirement of the robust problem described in Section~\ref{subsec:robust_joint_model} compared to the nominal problem (formulation~\modelP), we deemed it necessary to impose time limits on several components of the overall method. In particular, we impose a time limit of one hour on the overall procedure, with a time limit of 600 seconds for each solve of the restricted master problem, and 5 seconds for each solve of the subproblem. If the subproblem fails to identify a violated constraint within the 5 second time limit, it is solved again with a longer time limit of 120 seconds. If this second solve results in a violated constraint, the procedure continues; if it does not produce a violated constraint, the procedure terminates. 

For each of the same 20 instances corresponding to each combination $(n, K, c) \in \{30\} \times \{10, 20\} \times \{5\}$, for each value of $c' \in \{0.1, 0.2\}$ and for each value of $\Gamma$, we also compute the worst-case objective value of the nominal solution. We do this by solving the separation problem~\eqref{prob:robust_joint_separation} at the nominal product vector $\ab$. To make this worst-case objective value consistent with our robust procedure, we again impose a computation time limit of 120 seconds. We again compute the $\WCL$ and $\RI$ as in the experiments of the previous section. 

Table~\ref{table:robust_joint} reports the average of $\WCL$ and $\RI$ for each combination of $(n,K,c)$, $c'$ and $\Gamma$. In addition, it also reports the average computation time. From this table, we can see that in general, the $\WCL$ of the nominal solution can be substantial. For example, in the case where $c' = 0.1$ (i.e., each partworth deviates from its nominal value by at most 10\%), the $\WCL$ ranges from 4.31\% to 24.22\%. When $c' = 0.2$, it can be as high as 53.88\%. On the other hand, the robust solution can significantly outperform the nominal solution in terms of worst-case share-of-choice, as shown in the high values of $\RI$ (for example, with $K = 10$, $c' = 0.2$, $\Gamma = 40$, the $\RI$ is over 15\%). 

\begin{table}
\centering
\begin{tabular}{lrrrrrr} \toprule
$c'$ & $n$ & $K$ & $\Gamma$ & $\WCL$ (\%) & $\RI$ (\%) & Time (s) \\ \midrule
0.10 &  30 &  10 &  10 & 4.31 & 0.19 & 777.40 \\ 
  0.10 &  30 &  10 &  20 & 7.50 & 0.70 & 1612.02 \\ 
  0.10 &  30 &  10 &  30 & 9.87 & 2.79 & 2257.32 \\ 
  0.10 &  30 &  10 &  40 & 11.60 & 4.39 & 2333.21 \\ 
  0.10 &  30 &  10 &  50 & 13.10 & 5.02 & 2105.17 \\[0.25em] 
  0.10 &  30 &  20 &  20 & 8.80 & -0.97 & 3823.48 \\ 
  0.10 &  30 &  20 &  40 & 14.58 & -0.56 & 3989.73 \\ 
  0.10 &  30 &  20 &  60 & 18.79 & 2.93 & 4021.12 \\ 
  0.10 &  30 &  20 &  80 & 21.90 & 4.91 & 3999.85 \\ 
  0.10 &  30 &  20 & 100 & 24.22 & 9.41 & 3999.27 \\[0.5em]
  0.20 &  30 &  10 &  10 & 12.29 & -0.16 & 3368.03 \\ 
  0.20 &  30 &  10 &  20 & 21.76 & 3.13 & 3854.17 \\ 
  0.20 &  30 &  10 &  30 & 29.75 & 8.48 & 3846.31 \\ 
  0.20 &  30 &  10 &  40 & 36.45 & 15.88 & 3682.21 \\ 
  0.20 &  30 &  10 &  50 & 41.75 & 33.67 & 3644.28 \\[0.25em] 
  0.20 &  30 &  20 &  20 & 18.00 & -1.86 & 3928.02 \\ 
  0.20 &  30 &  20 &  40 & 30.58 & 0.99 & 3973.20 \\ 
  0.20 &  30 &  20 &  60 & 40.20 & 6.63 & 3911.56 \\ 
  0.20 &  30 &  20 &  80 & 48.14 & 18.73 & 3938.09 \\ 
  0.20 &  30 &  20 & 100 & 53.88 & 27.28 & 3907.19 \\ \bottomrule
\end{tabular}
\caption{Performance of robust solutions using the joint uncertainty set approach (Section~\ref{subsec:robust_joint_model}) on synthetic data instances. \label{table:robust_joint}}
\end{table}

In these results, we note that in a few cases the average $\RI$ is negative. Note that by the definition of $\RI$ in equation~\eqref{eq:RI_definition}, this cannot happen if $\ab^R$ \emph{exactly} solves the robust problem, and all worst-case share-of-choice objective values are computed \emph{exactly}. This is entirely an artifact due to the computation time limits that were applied when solving the robust problem and to evaluate the worst-case share-of-choice. In particular, due to the overall time limit of one hour, it is possible to have a suboptimal solution to the robust problem; as a result, even if one could perfectly compute $\min_{\betab \in \Ucal} F(\cdot, \betab)$ for such a solution and the nominal solution, it is possible that the nominal solution outperforms it in terms of worst-case objective value, leading to a negative $\RI$. In addition, it is also possible that the worst-case objective value of the nominal solution, $\min_{\betab \in \Ucal} F(\ab^N, \betab)$, is over-estimated, which can happen if problem~\eqref{prob:robust_joint_separation} is terminated early with a suboptimal solution; this in turn could result in a negative $\RI$ as well. 

Lastly, with regard to computation times, we note that the computation times for this approach are large, and in particular much larger than for the product uncertainty set approach tested in Section~\ref{subsec:robust_product_results}, which involved solving a single finite mixed-integer exponential cone program with the same number of exponential cones as the nominal problem~\modelP. (Note that in some cases, the computation time is higher than one hour, as the global one hour time limit was reached in the middle of a solve of either the restricted master problem or the worst-case separation problem.) From a tractability standpoint, these preliminary results suggest that the product uncertainty set approach is preferable to the joint uncertainty set approach.

\clearpage

\section{Approximation algorithm for fixed $K$ and fixed partworth magnitude}
\label{sec:approximation_algorithm}

In this section, we develop an approximation algorithm for solving the logit-based SOCPD problem. Section~\ref{subsec:approximation_algorithm_DP} presents a dynamic programming approach for solving the logit-based SOCPD problem in the case when all the partworth parameters $\{ \beta_{k,i} \}$ are integer-valued. Using this dynamic programming approach, Section~\ref{subsec:approximation_algorithm_FPTAS} presents an approximation algorithm for solving the logit-based SOCPD problem in general (where the partworth parameters need not be integer-valued). Lastly, Section~\ref{proof:theorem_algo_is_FPTAS} presents the proof of the approximation guarantee in Section~\ref{subsec:approximation_algorithm_FPTAS} (Theorem~\ref{theorem:algo_is_FPTAS}).

\subsection{Dynamic programming approach when partworths are integer}
\label{subsec:approximation_algorithm_DP}

We describe here a dynamic programming approach for solving the logit-based SOCPD problem when the partworths $\beta_{k,0}, \dots, \beta_{k,n}$ take integer values and when the set of product attribute vectors is unconstrained, i.e., $\Acal = \{0,1\}^n$. When we treat the number of customer types $K$ as a constant, this approach yields a pseudo-polynomial time algorithm for solving the logit-based SOCPD problem.

Recall that the logit-based SOCPD problem is
\begin{equation*}
\max_{\ab \in \Acal \equiv \{0,1\}^n} \sum_{k=1}^K \lambda_k \cdot \sigma( u_k(\ab)),
\end{equation*}
where for convenience, we use $\sigma(\cdot)$ to denote the logistic response function, i.e., $\sigma(u) = e^u / (1 + e^u)$. Let $F(\ab)$ denote the above objective function. %

Suppose that all of the utility parameters -- $\beta_{k,0}, \dots, \beta_{k,n}$ -- are integer valued. For each $k \in [K]$ and $i \in \{1,\dots,n+1\}$, define $u_{k,i,\max}$ and $u_{k,i,\min}$ as 
\begin{align*}
u_{k,i,\max} & = \beta_{k,0} + \sum_{j=1}^{i-1} ( \beta_{k,j})_+, \\
u_{k,i,\min} & = \beta_{k,0} + \sum_{j=1}^{i-1} ( \beta_{k,j})_-,
\end{align*}
where $(\cdot)_+ = \max\{0,\cdot\}$, $(\cdot)_- = \min\{0, \cdot\}$, and the sum is defined to be zero when the range of summation is empty (i.e., when $i = 1$). In words, $u_{k,i,\min}$ is the lowest possible value that $u_k(\cdot)$ can take when we are allowed to set $a_1,\dots, a_{i-1}$ arbitrarily, but $a_i, \dots, a_n$ are fixed to zero. Similarly, $u_{k,i,\max}$ is the largest possible value that $u_k(\cdot)$ can take when we fix $a_{i},\dots, a_n$ to zero. Let $u_{\max} = \max_{k \in [K]} u_{k,n+1,\max}$ and $u_{\min} = \min_{k \in [K]} u_{k,n+1, \min}$ be the largest and smallest possible utility values, respectively, attainable from setting all $n$ attributes over all $K$ customer types.  Finally, define $B$ and $b$ as 
\begin{align}
B & = \max_{k \in [K], i \in \{0,1,\dots,n\}} ( \beta_{k,i} )_+, \\
b & = \min_{k \in [K], i \in \{0,1,\dots,n\}} ( \beta_{k,i} )_-,
\end{align}
where $B$ and $b$ represent the largest positive and smallest negative partworth parameters, respectively, over all $k$ and $i$. Consequently, 
\begin{align}
u_{\max} & \leq (n+1) B, \\
u_{\min} & \geq (n+1) b,
\end{align}
and thus $u_{\max} - u_{\min} \leq (n+1) (B - b)$. In what follows, we will refer to $B - b$ as the \emph{partworth range parameter}. 

Let $\Vcal_{i,k}$ be the set of possible integer utilities between $u_{k,i,\min}$ and $u_{k,i,\max}$: 
\begin{equation}
\Vcal_{i,k} = \{ u_{k,i,\min}, u_{k,i,\min}+1, \dots, u_{k,i,\max} \}.
\end{equation}

Consider the following dynamic program, defined using the value functions $J_1,\dots, J_{n+1}$. For $i =1,\dots, n+1$, let $J_i: \Vcal_{i,1} \times \dots \times \Vcal_{i,K} \to \mathbb{R}$ be a function that satisfies the following recursion:
\begin{align*}
J_i(v_1,\dots, v_K) & = \max\{ J_{i+1}(v_1, \dots, v_K), J_{i+1}(v_1 + \beta_{1,i}, \dots, v_K + \beta_{K,i}) \}, \\
&  \forall i \in [n], (v_1,\dots, v_K) \in \prod_{k=1}^K \Vcal_{i,k},
\end{align*}
where we use the notation $[N] = \{1,\dots, N\}$, with the terminal conditions
\begin{align*}
J_{n+1}( v_1,\dots, v_K) = \sum_{k=1}^K \lambda_k \cdot \sigma( v_k ), \quad \forall (v_1,\dots, v_K) \in \prod_{k=1}^K \Vcal_{n+1,k}
\end{align*}
Observe that by solving this dynamic program, the value of $J_1(\beta_{1,0}, \dots, \beta_{K,0})$ yields the exact optimal value of the logit-based SOCPD problem. The optimal solution can be obtained by taking the greedy action with respect to the optimal value function. 

Note also that the running time of computing all of the values of $J$ using the DP recursion is $\sum_{i=1}^{n+1} \prod_{k=1}^K | \Vcal_{i,k} | = O( (n+1) (u_{\max} - u_{\min})^K ) ) = O( (n+1)^{K+1} (B - b)$, and the time to find the optimal $\ab$ by identifying the greedy action is $O(n)$. Thus, if $K$ is treated as a constant, then we can solve the problem in  pseudopolynomial time, i.e., in time that is polynomial in the magnitude of the inputs, as represented by the partworth range parameter $B - b$, and in the number of attributes $n$. We shall leverage this method in the next section to develop an approximation algorithm that solves the general logit-based SOCPD problem (i.e., where the partworths are no longer integer), in the regime where the number of customer types $K$ and the partworth range parameter $B - b$ are both treated as constants.

\subsection{Approximation algorithm for general partworths}
\label{subsec:approximation_algorithm_FPTAS}

Using the dynamic programming method developed in the previous section, we now develop an approximation algorithm for the general case when the partworths $\{ \beta_{k,i} \}$ are no longer integer. The overall strategy that we will take to construct our approximation algorithm is to discretize the utility parameters $\beta_{k,0}, \dots, \beta_{k,n}$ of each customer type. In particular, suppose that we are given a number $R > 0$, which will serve as a discretization parameter. Consider discretizing the partworths according to $R$:
\begin{equation*}
\tilde{\beta}_{k,j} = \left \lfloor \frac{ \beta_{k,j}}{R} \right \rfloor  , \quad \forall \ k \in [K], j \in \{0,1,\dots,n\}.
\end{equation*}
Define also the discretized utility function $\tilde{u}_k(\cdot)$ as
\begin{equation*}
\tilde{u}_k(\ab) = \tilde{\beta}_{k,0} + \sum_{j=1}^n \tilde{\beta}_{k,j} a_j.
\end{equation*}
Note that by multiplying $\tilde{u}_k(\ab)$ by $R$, we approximately obtain $u_k(\ab)$; that is, $R \tilde{u}_k(\ab) \approx u_k(\ab)$. Our goal now is to solve the discretized problem 
\begin{align*}
& \max_{\ab \in \{0,1\}^n} \hat{F}(\ab) \\
& \equiv \max_{\ab \in \{0,1\}^n} \sum_{k=1}^K \lambda_k \sigma( R \cdot \tilde{u}_k(\ab)),
\end{align*}
where $\hat{F}: \{0,1\}^n \to \mathbb{R}$ denotes the discretized logit-based share-of-choice objective. We will now describe a dynamic programming approach for solving this problem, which will turn out to be the FPTAS that we seek.

As with the DP approach in Section~\ref{subsec:approximation_algorithm_DP}, let us compute bounds $\tilde{u}_{k,i,\max}$ and $\tilde{u}_{k,i,\min}$ as
\begin{align*}
\tilde{u}_{k,i,\max} & = \tilde{\beta}_{k,0} + \sum_{j=1}^{i-1} (\tilde{\beta}_{k,j})_+, \\
\tilde{u}_{k,i,\min} & = \tilde{\beta}_{k,0} + \sum_{j=1}^{i-1} (\tilde{\beta}_{k,j})_-.
\end{align*}
Let us also define $\tilde{u}_{\max} = \max_{k \in [K]} \tilde{u}_{k,n+1,\max}$, $\tilde{u}_{\min} = \min_{k \in [K]} \tilde{u}_{k,n+1,\min}$. Note that in relation to $B$ and $b$ defined in Section~\ref{subsec:approximation_algorithm_DP}, $\tilde{u}_{\max}$ and $\tilde{u}_{\min}$ can be bounded as follows:
\begin{align*}
\tilde{u}_{\max} & \leq (n+1) \left \lfloor \frac{B}{R} \right \rfloor \\
& \leq (n+1) \cdot \frac{B}{R}, \\
\tilde{u}_{\min} & \geq (n+1) \left \lfloor \frac{b}{R} \right \rfloor \\
& \geq (n+1) \left( \frac{b}{R} - 1 \right).
\end{align*}

Finally, let $\tilde{\Vcal}_{i,k}$ be defined as 
\begin{equation}
\tilde{\Vcal}_{i,k} = \{ \tilde{u}_{k,i,\min}, \tilde{u}_{k,i,\min} + 1, \dots, \tilde{u}_{k,i,\max}  \}.
\end{equation}

Note that the discretized problem, 
\begin{align*}
\max_{\ab \in \{0,1\}^n} \hat{F}(\ab) \equiv \max_{\ab \in \{0,1\}^n} \sum_{k=1}^K \lambda_k \sigma( R \cdot \tilde{u}_k(\ab))
\end{align*}
can again be solved by dynamic programming, which we now describe. Define the value functions $J_1,\dots,J_{n+1}$, where for each $i = 1,\dots, n+1$ the function $J_1: \prod_{k=1}^K \tilde{\Vcal}_{i,k} \to \mathbb{R}$ satisfies the following recursion:
\begin{align}
J_i(v_1,\dots, v_K) & = \max\{ J_{i+1}(v_1, \dots, v_K), J_{i+1}(v_1 + \tilde{\beta}_{1,i}, \dots, v_K + \tilde{\beta}_{K,i}) \}, \nonumber \\
&  \forall i \in [n], (v_1,\dots, v_K) \in \prod_{k=1}^K \tilde{\Vcal}_{i,k}, \label{eq:Fhat_DP_recursion}
\end{align}
with the terminal values defined as
\begin{align}
J_{n+1}(v_1,\dots, v_K) & = \sum_{k=1}^K \lambda_k \sigma(R \cdot v_k), \nonumber \\
& \forall (v_1,\dots, v_K) \in \prod_{k=1}^K \tilde{\Vcal}_{i,k}. \label{eq:Fhat_DP_terminal}
\end{align}
The value of $J_1(\cdot)$ at $(\tilde{\beta}_{1,0}, \dots, \tilde{\beta}_{K,0})$, i.e., $J_1(\tilde{\beta}_{1,0}, \dots, \tilde{\beta}_{K,0})$, is exactly the optimal value of $\max_{\ab \in \{0,1\}^n} \hat{F}(\ab)$. In addition, the number of steps to calculate all of the value functions $J_1,\dots, J_{n+1}$ is bounded by
\begin{align*}
& \sum_{i=1}^{n+1} \prod_{k=1}^K | \tilde{\Vcal}_{i,k} | \\
& = \sum_{i=1}^{n+1} \prod_{k=1}^K ( \tilde{u}_{k,i,\max} - \tilde{u}_{k,i,\min}) \\
& \leq \sum_{i=1}^{n+1} (\tilde{u}_{\max} - \tilde{u}_{\min})^K \\
& = (n+1) (\tilde{u}_{\max} - \tilde{u}_{\min})^K \\ 
& = (n+1) \left( (n+1) \left( \frac{B}{R} - \frac{b}{R} + 1\right) \right)^K \\
& = (n+1)^{K+1} \left( \frac{B}{R} - \frac{b}{R} + 1 \right)^K  
\end{align*}
which implies that the computation time of the DP is $O( (n+1)^{K+1} \left( \frac{B}{R} - \frac{b}{R} + 1 \right)^K)$. 

An optimal solution $\hat{\ab} \in \arg \max_{\ab \in \{0,1\}^n} \hat{F}(\ab)$ can be obtained using Algorithm~\ref{algorithm:DP_greedy} below, which requires $O(n)$ steps.
\begin{algorithm}
\caption{Greedy algorithm for obtaining optimal solution from DP value function $J(\cdot)$. \label{algorithm:DP_greedy}}
\begin{algorithmic}
\STATE Initialize $a_1 \geq 0, \dots, a_n \gets 0$.
\STATE Initialize $v_1 \gets \tilde{\beta}_{1,0}, \dots, v_K \gets \tilde{\beta}_{K,0}$, $i \gets 1$.
\WHILE{ $i \leq n$ }
	\IF{$J(i+1, v_1,\dots, v_K) < J(i+1, v_1 + \tilde{\beta}_{1,i},\dots, v_K + \tilde{\beta}_{K,1})$}
		\STATE Set $a_i \gets 1$.
	\ELSE
		\STATE Set $a_i \gets 0$.
	\ENDIF 
	\FOR{ $k = 1,\dots, K$}
		\STATE Set $v_k \gets v_k + \tilde{\beta}_{k,i} a_i$.
	\ENDFOR
\ENDWHILE
\RETURN Solution $\ab$.
\end{algorithmic}
\end{algorithm}

We can now define our approximation algorithm for the logit-based SOCPD, which is presented below as Algorithm~\ref{algorithm:FPTAS}. 
\begin{algorithm}
\caption{Approximation algorithm for logit-based SOCPD problem. \label{algorithm:FPTAS}}
\begin{algorithmic}
\STATE Set $R \gets \epsilon / [ (n+1)K ]$
\STATE Solve $\max_{\ab \in \{0,1\}^n} \hat{F}(\ab)$ using the dynamic program \eqref{eq:Fhat_DP_recursion} - \eqref{eq:Fhat_DP_terminal}. 
\STATE Using Algorithm~\ref{algorithm:DP_greedy}, obtain $\hat{\ab} \in \arg \max_{\ab \in \{0,1\}^n} \hat{F}(\ab)$. 
\RETURN Approximate solution $\hat{\ab}$. 
\end{algorithmic}
\end{algorithm}

The following theorem establishes the performance guarantee and runtime of Algorithm~\ref{algorithm:FPTAS}. 
\begin{theorem}
For any $\epsilon \in (0,1)$, Algorithm~\ref{algorithm:FPTAS} returns a $(1 - \epsilon)$-optimal solution the logit-based SOCPD problem, in running time $O \left( (n+1)^{2K + 1} K^K \left(\frac{1}{\epsilon} \right)^K (B - b + 1)^K \right)$. 
\label{theorem:algo_is_FPTAS}
\end{theorem}

Note that Algorithm~\ref{algorithm:FPTAS} is not a fully polynomial time approximation scheme (FPTAS) in general, because the runtime has an exponential dependence on the number of customer types $K$, and has a polynomial dependence on the partworth range parameter $B - b$, which implies an exponential dependence on the input size (i.e., the number of bits required to encode the problem). When $K$ and $B - b$ are treated as constant quantities, i.e., $K = O(1)$ and $(B - b) = O(1)$, then the runtime is polynomial in $n$ and $1/\epsilon$, and Algorithm~\ref{algorithm:FPTAS} could be be viewed as a fully polynomial time approximation scheme.

We now discuss the dependence of the runtime on $K$ and $B - b$ in more detail, and how these dependencies compare to practical settings in which one would need to solve the logit-based SOCPD problem:\\

\noindent \textbf{Dependence on $K$:} With regard to the dependence on $K$, we note that the assumption of a constant number of customer types is common in the assortment optimization literature that develops approximation algorithms. For example, the paper of \cite{desir2022capacitated} develops an approximation algorithm for the mixture MNL assortment problem that has an exponential dependence on $K$, but becomes an FPTAS when $K = O(1)$. %

In our context, when the customer population is represented using a latent-class MNL model that is estimated from a conjoint dataset, it is typically not necessary to consider huge values of $K$. In Section~\ref{appendix:numerics_tuning_K} of the ecompanion, we provide evidence that the appropriate value of $K$ is small by considering the Akaike information criterion (AIC), the Bayesian information criterion (BIC) and the consistent AIC (CAIC), which are commonly used measures in marketing science to tune $K$. Note however that for hierarchical Bayesian models, it is more conventional to assume that customer types correspond to individual respondents, or to multiple draws from the posterior distribution of each respondent's partworth vector. In these types of models $K$ could be large and would scale with the number of respondents in the conjoint dataset (as is the case in our experiments in Section~\ref{subsec:numerical_experiments_real}). \\

\noindent \textbf{Dependence on $B - b$:} With regard to the dependence on $(B - b)$, we make the following remarks. First, virtually all product design instances will use partworths that are estimated in some way from conjoint data. This implies that the partworths are subject to internal storage constraints from how the corresponding estimation algorithm is implemented on a computer. For example, our implementation of the EM algorithm in Julia represents all continuous quantities using Julia's native 64 bit floating point type, {\tt Float64}, which naturally imposes a limit on $(B - b)$. We do not envision a practical setting where one would need to consider partworths that would result in values of $B - b$ so large that one would need to go beyond the precision of standard floating point representations. %

Second, when estimating latent-class MNL models, it is sometimes the case that as one considers larger and larger values of $K$, the LC-MNL model that is returned is such that there are some segments with low segment probabilities $\lambda_k$ and partworth vectors $\betab_k$ containing entries with large magnitudes (either positive or negative). This would consequently lead to a high overall value of $B - b$. This occurs because as one increases $K$ starting with $K = 1$, the EM algorithm initially identifies segments with large $\lambda_k$'s in the data that greatly contribute to the fit, but beyond a certain number of segments, the subsequent segments that are identified have small $\lambda_k$'s and contribute less and less to the fit. In Section~\ref{appendix:partworth_distributions}, we illustrate this for our four data sets by plotting the distributions of the magnitudes of the partworth parameters; in general, as $K$ increases, the distributions become more diffuse. 

In practice, however, such LC-MNL models, where $K$ is large and where some segments have large partworth magnitudes and small segment probabilities, are unlikely to be justified from a data fitting perspective. As we noted above, the four data sets we consider in Section~\ref{subsec:numerical_experiments_real} are all well explained by values of $K$ between 5 and 30 (see Section~\ref{appendix:numerics_tuning_K}). Besides model selection measures, another important consideration is whether the segments produced make sense; on this point, \cite{bucklin1992brand} note (see page 206): 
\begin{quote}
\emph{
Though improvement in model fit is a critical first cut at determining the number of segments, a managerially relevant criterion is interpretability. If additional insight is not provided by increasing the number of response segments, including additional segments may not be worthwhile even though they improve the model fit. }
\end{quote}

Third, with regard to partworth magnitudes for mixture MNL models that arise in hierarchical Bayesian models, we note that it is common to assume that the mean of the second-stage prior for the mean of the partworth vector mixture distribution is zero (see equation~\eqref{eq:HB_second_stage_barbetab} in our specification of our hierarchical Bayesian model in Section~\ref{appendix:numerics_HB_specification}, which is implemented in the R package {\tt bayesm}). This naturally has the effect that the mean of the mixture distribution of the $\betab$ is shrunk towards zero, and consequently that the mean of the posterior distribution of each respondent's partworth vector is also shrunk towards zero. As a result, the posterior mean estimate of each respondent's partworth vector should generally not have large magnitudes. In Section~\ref{appendix:partworth_distributions}, we confirm that this is indeed the case by plotting the distributions of the partworths for the four hierarchical Bayesian models we consider in Section~\ref{subsec:numerical_experiments_real}. \\

We conclude this section by noting that we did investigate whether it is possible to obtain a stronger approximation algorithm that would match what is known for the mixture MNL assortment problem (i.e., an approximation algorithm that is an FPTAS when $K = O(1)$, and that does not have an exponential dependence on input size). However, we were not successful in establishing such an algorithm. The design of such an approximation algorithm for this problem remains an open research question.

\subsection{Proof of Theorem~\ref{theorem:algo_is_FPTAS}}
\label{proof:theorem_algo_is_FPTAS}

To prove Theorem~\ref{theorem:algo_is_FPTAS}, we first establish several auxiliary results. The first is Lemma~\ref{lemma:discretization_gap}, which states that the discretized utility function $R \tilde{u}_k$ underapproximates the true utility function $u_k(\ab)$, and that the gap between this discretized utility function and the true utility function is uniformly bounded by $(n+1)R$. 

\begin{lemma}
For any $\ab \in \{0,1\}^n$, we have that $u_k(\ab) \geq R \tilde{u}_k(\ab)$ and that $u_k(\ab) - R \tilde{u}_k(\ab) \leq (n+1)R$. 
\label{lemma:discretization_gap}
\end{lemma}
\begin{proofvvm}
For the inequality $u_k(\ab) \geq R \tilde{u}_k(\ab)$, observe that
\begin{align*}
R \tilde{u}_k(\ab) & = R \tilde{\beta}_{k,0} + \sum_{j=1}^n R \tilde{\beta}_{k,j} a_j \\
& = R \cdot \left \lfloor \frac{ \beta_{k,0}}{R} \right \rfloor + \sum_{j=1}^n R \cdot \left \lfloor \frac{ \beta_{k,j}}{R} \right \rfloor \cdot a_j \\
& \leq R \cdot \frac{ \beta_{k,0}}{R} + \sum_{j=1}^n R \cdot \frac{ \beta_{k,j}}{R} \cdot a_j \\
& = u_k(\ab),
\end{align*}
where the inequality follows because $\lfloor x \rfloor \leq x$. To see the second part of the lemma, observe that
\begin{align*}
u_k(\ab) - R \tilde{u}_k(\ab) & = \left( \beta_{k,0} - R \cdot \left \lfloor \frac{ \beta_{k,0}}{R} \right \rfloor \right) + \sum_{j=1}^n \left( \beta_{k,j} - R \cdot \left \lfloor \frac{ \beta_{k,j}}{R} \right \rfloor \right) \cdot a_j \\
& \leq R + \sum_{j=1}^n R \cdot a_j \\
& \leq (n+1) R,
\end{align*}
where the first inequality follows because for any $x$ and any positive $R$, we have $\lfloor x / R \rfloor \geq x / R - 1$, which implies that $x - R \lfloor x / R \rfloor \leq x - R \cdot (x/R - 1) = R$. \Halmos
\end{proofvvm}

A straightforward consequence of this lemma is that the discretized share-of-choice function $\hat{F}(\ab)$ always underapproximates the true share-of-choice function $F(\ab)$, which is captured in the next lemma.
\begin{lemma}
For all $\ab \in \{0,1\}^n$, $\hat{F}(\ab) \leq F(\ab)$.
\label{lemma:Fhat_underapproximates_F}
\end{lemma}
\begin{proofvvm}
We have
\begin{align*}
\hat{F}(\ab) & = \sum_{k=1}^K \lambda_k \cdot \sigma( R \cdot \tilde{u}_k(\ab)) \\
& \leq \sum_{k=1}^K \lambda_k \cdot \sigma( u_k(\ab)) \\
& = F(\ab),
\end{align*}
where the inequality follows by the first part of Lemma~\ref{lemma:discretization_gap} and the monotonicity of $\sigma(\cdot)$. \Halmos
\end{proofvvm}

Armed with Lemma~\ref{lemma:discretization_gap} and Lemma~\ref{lemma:Fhat_underapproximates_F}, we can prove the following guarantee on the quality (in terms of relative gap) of the solution of the discretized problem. 
\begin{lemma}
Let $\hat{\ab}$ be an optimal solution of $\max_{\ab} \hat{F}(\ab)$, and $\ab^*$ be an optimal solution of $\max_{\ab} F(\ab)$. Then
\begin{align*}
\frac{F(\ab^*) - F(\hat{\ab})}{F(\ab^*)} \leq K \cdot (n+1) \cdot R. 
\end{align*}
\label{lemma:relative_gap}
\end{lemma}
\begin{proofvvm}
We have:
\begin{align*}
\frac{F(\ab^*) - F(\hat{\ab})}{F(\ab^*)}  & \leq \frac{F(\ab^*) - \hat{F}(\ab^*) + \hat{F}(\ab^*) - \hat{F}(\hat{\ab}) + \hat{F}(\hat{\ab}) - F(\hat{\ab})}{F(\ab^*)} \\
& \leq \frac{F(\ab^*) - \hat{F}(\ab^*)}{F(\ab^*)} \\
& = \frac{ \sum_{k=1}^K \lambda_k \sigma( u_k(\ab^*)) - \sum_{k=1}^K \lambda_k \sigma( R \tilde{u}_k(\ab^*))}{ \sum_{k=1}^K \lambda_k \sigma(u_k(\ab^*))} \\
& = \sum_{k=1}^K \lambda_k \cdot \frac{ \sigma( u_k(\ab^*)) - \sigma( R \tilde{u}_k(\ab^*)) }{\sum_{k=1}^K \lambda_k \sigma(u_k(\ab^*))} \\
& \leq \sum_{k=1}^K \lambda_k \cdot \frac{ \sigma( u_k(\ab^*)) - \sigma( R \tilde{u}_k(\ab^*)) }{ \lambda_k \cdot \sigma(u_k(\ab^*))} \\
& = \sum_{k=1}^K \frac{ \sigma( u_k(\ab^*)) - \sigma( R \tilde{u}_k(\ab^*)) }{ \sigma(u_k(\ab^*))} \\
& = \sum_{k=1}^K \left( 1 - \frac{1 + e^{-u_k(\ab^*)}}{1 + e^{-R \tilde{u}_k(\ab^*)}} \right) \\
& = \sum_{k=1}^K \frac{e^{-R \tilde{u}_k(\ab^*)} - e^{-u_k(\ab^*)}}{1 + e^{-R \tilde{u}_k(\ab^*)}} \\ 
& = \sum_{k=1}^K \frac{1 - e^{R \tilde{u}_k(\ab^*) -u_k(\ab^*)}}{e^{R \tilde{u}_k(\ab^*)} + 1 } \\
& \leq  \sum_{k=1}^K \left( 1 - e^{- (u_k(\ab^*) - R \tilde{u}_k(\ab^*))} \right) \\
& \leq \sum_{k=1}^K ( u_k(\ab^*) - R \tilde{u}_k(\ab^*) ) \\
& \leq K \cdot (n+1) \cdot R 
\end{align*}
where the first step follows by algebra; the second step follows because $\hat{F}(\ab^*) - \hat{F}(\hat{\ab}) \leq 0$ (this is true because $\hat{\ab}$ is an optimal solution of $\max_{\ab} \hat{F}(\ab)$) and $\hat{F}(\hat{\ab}) - F(\hat{\ab}) \leq 0$ (this follows by Lemma~\ref{lemma:Fhat_underapproximates_F}); the third and fourth step follow by algebra; the fifth step follows because $\sum_{k=1}^K \lambda_k \sigma(u_k(\ab^*)) \geq \lambda_{k'} \sigma(u_{k'}(\ab^*))$ for any $k'$; the sixth, seventh, eighth and ninth steps follow by algebra; the tenth step follows by the fact that the denominator $e^{R \tilde{u}_k(\ab^*)} + 1$ is lower bounded by 1; the eleventh step follows because $1 - e^{-x} \leq x$ for any real $x$; and the final step follows by the second part of Lemma~\ref{lemma:discretization_gap}. \Halmos
\end{proofvvm}

With all of these results established, we now finally verify Theorem~\ref{theorem:algo_is_FPTAS}. Let $\hat{\ab}$ be the solution produced by Algorithm~\ref{algorithm:FPTAS}. The solution $\hat{\ab}$ produced by Algorithm~\ref{algorithm:FPTAS} solves the approximate problem $\max_{\ab \in \{0,1\}^n} \hat{F}(\ab)$. By Lemma~\ref{lemma:relative_gap}, this solution is a $(1 - K(n+1)R)$-optimal solution; for $R = \epsilon / (K(n+1))$, we thus have that it is a $(1 - \epsilon)$-optimal solution. 

With regard to the running time, the running time of the DP recursion in equations \eqref{eq:Fhat_DP_recursion} - \eqref{eq:Fhat_DP_terminal} is upper bounded by $( (n+1)^{K+1} \left( \frac{B}{R} - \frac{b}{R} + 1 \right)^K)$. Since $R = \epsilon / (K(n+1))$, we have the running time of the DP is upper bounded as follows:
\begin{align*}
&  (n+1)^{K+1} \left( B \cdot \frac{1}{\epsilon} \cdot K(n+1) - b \cdot \frac{1}{\epsilon} \cdot K(n+1) + 1 \right)^K \\
& \leq (n+1)^{K+1} \left( B \cdot \frac{1}{\epsilon} \cdot K(n+1) - b \cdot \frac{1}{\epsilon} \cdot K(n+1) + \frac{1}{\epsilon} \cdot K(n+1) \right)^K \\ 
& = (n+1)^{2K + 1} K^K \left(\frac{1}{\epsilon} \right)^K (B - b + 1)^K,
\end{align*}
where the inequality follows from the fact that $K \geq 1$, $n \geq1$ and $0 \leq \epsilon \leq 1$. Additionally, the number of steps to run Algorithm~\ref{algorithm:DP_greedy} after computing the value function is $O(n)$. Therefore, the overall complexity is $O \left( (n+1)^{2K + 1} K^K \left(\frac{1}{\epsilon} \right)^K (B - b + 1)^K \right)$. \Halmos

\end{document}